\documentclass[amsppt,11pt]{amsart}
\usepackage[backref]{hyperref}
\usepackage{curves}
\usepackage{amsmath}

    \usepackage{pdfsync}

\newtheorem {theorem}{Theorem}[section]

\newtheorem {proposition}[theorem]{Proposition}
\newtheorem {lemma}[theorem]{Lemma}

\newcounter{conjecture}
\newcounter{remark}

\newcommand{\eqnsection}{
   \renewcommand{\theequation}{\thesection.\arabic{equation}}
   \makeatletter
   \csname @addtoreset\endcsname{equation}{section}
   \makeatother}

\def \be{\begin{equation}}
\def \ee{\end{equation}} 
\def \bt{\begin{theorem}}
\def \et{\end{theorem}}
\def \bea{\begin{eqnarray}}
\def \eea{\end{eqnarray}}
\def \bas{\begin{eqnarray*}}
\def \eas{\end{eqnarray*}}

\def \bl{\begin{lemma}} 
\def \el{\end{lemma}}

\def \al{\alpha}
\def \bb{\beta}
\def \ga{\gamma}

\def \de{\delta}
\def \De{\Delta}
\def \ep{\epsilon}

\newcommand{\eps}{\varepsilon}

\def \la{\lambda}

\def \om{\omega}

\def \si{\sigma}

\def \th{\theta}

\def \ze{\zeta}

\def \ff{\infty}
\def \wh{\widehat}
\def \wt{\widetilde}

\newcommand{\ls}[1]
   {\dimen0=\fontdimen6\the\font \lineskip=#1\dimen0
\advance\lineskip.5\fontdimen5\the\font \advance\lineskip-\dimen0
\lineskiplimit=.9\lineskip \baselineskip=\lineskip
\advance\baselineskip\dimen0 \normallineskip\lineskip
\normallineskiplimit\lineskiplimit \normalbaselineskip\baselineskip
\ignorespaces }

\def \R{{\bf R}}

\def \Z{{\bf Z}}
\def \S{{\bf S}}

\def \AA{{\mathcal A}}
\def \BB{{\mathcal B}}

\def \HH{{\mathcal H}}
\def \II{{\mathcal I}}

\def \KK{{\mathcal K}}
\def \LL{{\mathcal L}}
\def \MM{{\mathcal M}}

\def \PP{{\mathcal P}}

\def \WW{{\mathcal W}}

\def \({\left(}
\def \){\right)}

\def \lc{\left\{}
\def \rc{\right\}}

\def \nn{\nonumber}

\def \bc{\begin{center} }
\def \ec{\end{center} }

\eqnsection

\newcommand{\beq}[1]{\begin{equation}\label{#1}}
\newcommand{\eeq}{\end{equation}}

\newcommand{\beqn}[1]{\begin{eqnarray}\label{#1}}
\newcommand{\eeqn}{\end{eqnarray}}

\def\squarebox#1{\hbox to #1{\hfill\vbox to #1{\vfill}}}
\renewcommand{\qed}{\hspace*{\fill}
            \vbox{\hrule\hbox{\vrule\squarebox{.667em}\vrule}\hrule}\smallskip}

\newcommand{\beaa}{\begin{eqnarray*}}
\newcommand{\eeaa}{\end{eqnarray*}}

\newcommand{\Pbm}{\mathbb{P}}
\newcommand{\Ebm}{\mathbb{E}}
\newcommand{\Pgw}{P^{\mbox{\fontsize{.03in}{1in}\rm GW}}}

\newcommand{\BM}{W}

\newcommand{\dwa}{d^1_{\mbox{\fontsize{.03in}{1in}\rm Wa}}}

 \date{March 27, 2022}

\title[Tightness for Thick Points in two dimensions ]
{Tightness for Thick Points  in two dimensions }

 \author[Jay Rosen]
   {Jay Rosen}

\thanks{Jay Rosen was partially  supported by the Simons Foundation}

\begin{document}

\begin{abstract}
\noindent
Let $W_{t}$ be Brownian motion in the plane started at the origin  and let $  \theta$ be the first exit time of the unit disk $D_{1}$. Let
\[
 \mu_{ \theta } (   x,\epsilon) =\frac{1}{\pi\epsilon^{ 2} }\int_{0}^{  \th }1_{\{   B(   x,\epsilon)\}}(  W_{t})\,dt,\]
and set $\mu^{ \ast}_{ \theta } (\epsilon)=\sup_{x\in D_{1}}\mu_{ \theta } (   x,\epsilon)$. We show that
\[
\sqrt{\mu^{ \ast}_{ \theta } (\epsilon)}-\sqrt{2/\pi}\( \log \epsilon^{-1}- \frac{1}{2}\log\log \epsilon^{-1}\) \]
is tight.

\end{abstract}

\maketitle

\section{Introduction}

Let $W_{t}$ be Brownian motion in the plane started at the origin  and let $  \th$ be the first exit time of the unit disk $D_{1}$. In \cite{DPRZ} we showed that
\begin{equation}
\lim_{\ep\to 0}\sup_{x\in D_{1}}\frac{1}{\ep^{ 2}\log^{ 2}(   \ep)}\int_{0}^{  \th }1_{\{   B(   x,\ep)\}}(  W_{t})\,dt=2,\quad a.s., \label{tp.0}
\end{equation}
where $ B(   x,\ep)$ is the ball of radius $\ep$ centered at $x$. The integral above is the occupation measure of $ B(   x,\ep)$, and points $x$ with large occupation measure are referred to as thick points. Taking square roots we can write this as
\begin{equation}
\lim_{\ep\to 0}\frac{1}{\log (  \ep^{ -1})}\sqrt{\sup_{x\in D_{1}}\frac{1}{\pi\ep^{ 2} }\int_{0}^{  \th }1_{\{   B(   x,\ep)\}}(  W_{t})\,dt}=\sqrt{2/\pi},\quad a.s.
 \label{tp.01}
\end{equation}
Let
\begin{equation}
 \mu_{ \th } (   x,\ep) =\frac{1}{\pi\ep^{ 2} }\int_{0}^{  \th }1_{\{   B(   x,\ep)\}}(  W_{t})\,dt,\label{plane.1}
\end{equation}
and set $\mu^{ \ast}_{ \th } (\ep)=\sup_{x\in D_{1}}\mu_{ \th } (   x,\ep)$.
Then (\ref{tp.01}) says that   $\sqrt{\mu^{ \ast}_{ \th } (\ep)}\sim \sqrt{2/\pi}\log   \ep^{ -1} $, as $\ep\to 0$. In this paper we obtain   more detailed asymptotics. Let 
\begin{equation}
m_{\ep}=\sqrt{2/\pi}\( \log \ep^{-1}- \frac{1}{2}\log\log \ep^{-1}\).\label{121}
\end{equation}
We will say that the thick points in $D_{1}$ are tight if $\sqrt{\mu^{ \ast}_{ \th } (\ep)}-m_{\ep} $ is a tight family of random variables. That is, 
\begin{equation}
\lim_{K\to\ff} \varlimsup_{\ep\to 0} \Pbm\( \Big|\sqrt{\mu^{ \ast}_{ \th } (\ep)}-m_{\ep}\Big| >K\)=0.\label{0.01gplane}
\end{equation}

\bt\label{theo-tight0plane}
The thick points  in $D_{1}$ are tight.
\et

 In fact we obtain the following improvement on the right tail of (\ref{0.01gplane}), where $c^{ \ast}: =2\sqrt{2}$.

\bt\label{theo-planetight}
On $D_{1}$, for some $0< C, C' ,z_{0}<\ff$ and all $z\geq z_{0}$,
\begin{equation}
 \varlimsup_{\ep\to 0}  \Pbm\(\sqrt{\mu^{ \ast}_{ \th } (\ep)} -m_{\ep} \geq z\)\leq C z e^{ -  c^{ \ast} \sqrt{\pi}\,  z },\label{goal.1plane}
\end{equation}
and
\begin{equation}
\varliminf_{\ep\to 0}   \Pbm\(\sqrt{\mu^{ \ast}_{ \th } (\ep)} -m_{\ep} \geq z\)\geq  C'ze^{ -c^{ \ast}\sqrt{\pi}\, z}.\label{goal.1aplane}
\end{equation}
\et

It follows from Brownian scaling that Theorems \ref{theo-tight0plane} and \ref{theo-planetight}
hold if $D_{1}$ is replaced by any disc centered at the origin.

For reasons  of symmetry it is easier to work on the sphere $\S^{ 2}$, and derive our results for  thick points  in $D_{1}$ from results for thick points on   $\S^{ 2}$. We use  $ B_{d}(   x,r)$ for  the ball centered at $x$  of radius $r$, in the spherical metric $d$.

Let $X_{t}$ be Brownian motion on $\S^{ 2}$ started at some point $v$ (the `South Pole'). For some (small) $r^{ \ast}$ let $ \tau$ be the first hitting time of $\partial B_{d}(  v,r^{ \ast})$ (the   `Antarctic Circle').  Let $\om_{\ep}=2\pi (   1-\cos \ep)$, the area of $B_{d}(   x,\ep)$, and set
\begin{equation}
\bar\mu_{ \tau } (   x,\ep) = \frac{1}{\om_{\ep}}\int_{0}^{ \tau }1_{\{   B_{d}(   x,\ep)\}}(  X_{t})\,dt,\label{tp.1}
\end{equation}

With $\bar\mu_{ \tau, \ep }^{\ast} =\sup_{x\in \S^{ 2}} \bar\mu_{ \tau } (   x,\ep)$ we will say that the thick points on $\S^{2}$ are tight if $\sqrt{\bar\mu_{ \tau, \ep }^{\ast} }-m_{\ep} $ is a tight family of random variables. 

\bt\label{theo-tight0}
The thick points on $\S^{2}$ are tight.
\et

As in Theorem \ref{theo-planetight} we obtain the following improvement for the right tail.

\bt\label{theo-spheretight}
On $\S^{2}$, for some $0< C, C' ,z_{0}<\ff$ and all $z\geq z_{0}$,
\begin{equation}
 \varlimsup_{\ep\to 0}  \Pbm\(\sqrt{\bar\mu_{ \tau, \ep }^{\ast} } -m_{\ep} \geq z\)\leq C ze^{ -c^{ \ast} \sqrt{\pi}\, z},\label{goal.1asto}
\end{equation}
and
\begin{equation}
\varliminf_{\ep\to 0}   \Pbm\(\sqrt{\bar\mu_{ \tau, \ep }^{\ast} } -m_{\ep} \geq z\)\geq  C' z e^{ -c^{ \ast} \sqrt{\pi}\,  z }.\label{goal.1astow}
\end{equation}
\et
 
 Theorems \ref{theo-tight0} and \ref{theo-spheretight} are stated and first proven for $r^{ \ast}$ sufficiently small.
In  Section \ref{sec-euclid} we show that they  hold for any $0<r^{\ast}<\pi$.

In analogy with  \cite{DPRZ}, rather than work directly with  occupation measures,
we work with  excursion counts. To define this let
  $h_{l}=2\arctan ( r_{0}e^{ -l}/2)$ with $r_{0}$ small. For some $d_{0}\leq 1/1000$ let $F_{l}$ be the centers of a $d_{0}h_{l}$ covering of $\S^{2}$.

 Let $\mathcal{T}_{x, l}^{\tau }$   be the number of excursions from  $\partial B_{d}\left(x,h_{l-1}\right)$ to $ \partial B_{d}\left(x,h_{l}\right)$  prior to $ \tau$. 
  We will  obtain the following result for $\sup_{x\in F_{L}} \mathcal{T}_{x,L}^{\tau}$.

\bt\label{theo-excovertight}On $\S^{2}$,  for some $0<z_{0},C, C' <\ff$, all $L $ large and all $z_{0}\leq z\leq \log L$,
\begin{equation}
 C'  z e^{ -2z } \leq  \Pbm\(\sup_{x\in F_{L}}\sqrt{2  \mathcal{T}_{x,L}^{\tau}} - (  2 L- \log L)\geq z\)\leq C  ze^{ -2z }.\label{goal.1astos}
\end{equation}
Equivalently
\begin{equation}
 C'   ze^{ -2z } \leq  \Pbm\( \sup_{x\in F_{L}}  \mathcal{T}_{x,L}^{\tau} \geq 2L\(L-  \log L+ z\)   \)\leq Cz  e^{ -2z }.\label{goal.1ast}
\end{equation}
\et

Since $L\sim  \log h_{L}^{-1} $, Theorem \ref{theo-excovertight} is then suggestive of Theorem \ref{theo-spheretight} if we knew that on average the occupation measure of  $ B_{d}\left(x,h_{L}\right)$ during an excursion from $\partial B_{d}\left(x,h_{L}\right)$ to $ \partial B_{d}\left(x,h_{L-1}\right)$  was `about' $h^{ 2}_{L} $. While this is basically known for our choices of $h_{L},h_{L-1}$, see \cite[Lemma 6.2]{DPRZ}, it is more delicate  to get the precision   necessary to show the equivalence of Theorem \ref{theo-excovertight} with  (\ref{goal.1asto}).

We now write (\ref{goal.1astos}) in a more convenient form.
Set
\begin{equation}
\rho_{L} =2-\frac{\log L}{L}.\label{eq:defofts}
\end{equation}
We will prove  the following version of  Theorem \ref{theo-excovertight}.
\bt\label{theo-tight}On $\S^{2}$,  for some $0<z_{0},C, C' <\ff$, all $L $ large and all $z_{0}\leq z\leq \log L$,
\begin{equation}
 C'  z e^{ -2z } \leq  \Pbm\left[\sup_{x\in F_{L}}\sqrt{2\mathcal{T}_{x,L}^{\tau}}\,\geq \rho_{L}L+z \right]\leq C z e^{ -2z}.\label{goal.1}
\end{equation}
\et

\subsection{Background}  This paper is based in many ways on my work \cite{BRZ} with Belius and Zeituni on tightness for the cover time of $\S^{2}$. The general approach is similar, and whenever results of that paper could be used directly I did so. However, the mathematics often necessitated different arguments.

The family \[\{\bar\mu_{ \tau } (   x,\ep); x\in B_{d}(  v,r^{ \ast}), \ep>0\} \] is associated with a second order Gaussian chaos
 $H  (   x,\ep) $,  $x\in B_{d}(  v,r^{ \ast}), \ep>0$ by an isomorphism
 theorem of Dynkin \cite{D}.    Intuitively, 
\begin{equation}
H  (   x,\ep)=\int_{B_{d}( x,\ep)} G^{2}_{y} \,dm(y)- E\(\int_{B_{d}( x,\ep)} G^{2}_{y} \,dm(y)\)\label{back.1}
\end{equation}
 where  $G_{x}$ is the mean zero Gaussian process with covariance $u(x,y)$,  the Green's function for $  B_{d}(  v,r^{ \ast})$ and  $m$ denotes the standard surface measure on $S^{2}$. Since $u(x,x)=\ff$ for all $x$, (\ref{back.1}) is not a priori well defined.  Nevertheless, this would suggest that there is a close relationship between $\bar\mu_{ \tau, \ep }^{\ast} =\sup_{x\in R^{ 2}} \bar\mu_{ \tau } (   x,\ep)$ and the supremum of Gaussian fields. For details on $H$ and the  isomorphism
 theorem see \cite[Section 2]{MR}.
 
For related work see \cite{Abe, J1,J2, J3, J4}.
 
 \subsection{Open problems}
\label{sec-open}
\hspace{.3 in}

 1.  Based on the analogy with the extrema of Branching random walks 
    and log-correlated Gaussian fields, one expects that Theorem 
    \ref{theo-tight0plane} should be replaced by the statement 
    \begin{eqnarray*} 
    &&\mbox{\em The sequence of random variables 
      $ \sqrt{\mu^{\ast}_{\th}(\eps)}- m_{\eps} $ converges }\\
    && \mbox{\em  
    in distribution
  to a randomly shifted Gumbel random variable.}
\end{eqnarray*}
A key step in proving such convergence would be 
the improvement of the tail estimates in Theorems  \ref{theo-planetight} and \ref{theo-spheretight} for $z$
large,
which in turn would require a corresponding improvement of Theorem
\ref{theo-tight}.

2. In \cite{DPRZ} we also proved a conjecture of Erd\H{o}s and  Taylor
concerning the number   $L^{ \ast}_n$ of visits to the most visited site for
simple random walk in
$\Z^2$ up to step $n$. It was  shown there that
\begin{equation}
\lim_{n\to \ff}\frac{L^{ \ast}_n}{ (\log n)^2} =1/\pi
\hspace{.1in}\mbox{a.s.}\label{0.0}
\end{equation}
  The approach in that paper was to first prove (\ref{tp.0}) for
planar  Brownian motion and then to use strong approximation. 
Subsequently, in
\cite{R}, we presented a purely random walk method 
to prove (\ref{0.0}) for simple random walk. A natural problem is to prove tightness for $\sqrt{L^{ \ast}_n}$.

3. Following \cite{DPRZ} we analyzed thick points for several other process. See \cite{DPRZtransient}, \cite{DPRZspatial} and \cite{DPRZint}. One can ask about tightness or some analog for these processes.

\subsection{Structure of the paper } 
In Section \ref{sec-upperboundex} we obtain the upper bounds for excursion counts in Theorem \ref{theo-tight},  and in Section \ref{sec-lowerboundex} we derive the lower bounds. These sections employ many of the tools developed in \cite{BRZ}. In Section \ref{sec-occmeas} we show how to go from results on excursion counts to Theorems \ref{theo-tight0} and \ref{theo-spheretight} which involve  $\bar\mu_{ \tau } (   x,\ep)$ in $\S^{2}$. Here we have to deal with a new problem for the upper bounds: $\bar\mu_{ \tau } (   x,\ep)$ in $\S^{2}$ is not in general monotone in $\ep$. This requires interpolation and a continuity estimate which are developed in Sections \ref{sec-interp} and \ref{sec-var}. 
 In the short Section \ref{sec-euclid} we derive our results on thick points for the unit disc in the plane from our results on thick points for $\S^{2}$, and use this to show that  Theorems \ref{theo-tight0} and \ref{theo-spheretight}   hold for any $0<r^{\ast}<\pi$. The last section is an Appendix containing the barrier estimates we need for Sections \ref{sec-upperboundex} and \ref{sec-lowerboundex}.


\newpage
\subsection{Index of Notation}
The following are frequently used notation, and a pointer to the location where the definition appears.

$ \begin{array}{ll}
 \mu_{ \th } (   x,\ep)&\mbox {\eqref{plane.1}}\\  
 
 m_{\ep}&\mbox {\eqref{121}}\\    
 
c^{\ast}&\mbox {\eqref{goal.1plane}}\\  
 
  \bar\mu_{ \tau } (   x,\ep), B_{d}(x,r)&\mbox {\eqref{tp.1}}\\

\rho_{L}& \mbox{(\ref{eq:defofts})}\\
\mathcal{T}_{x, l}^{\tau }&\mbox{(\ref{goal.1astos})}\\

r_l,h_l&\mbox{\eqref{dr.1}}\\
F_l&\mbox {\eqref{dr.2}}\\
 \mathcal{T}_{y,l }^{ k\to 0}& \mbox{\rm \eqref{eq:NotHitByrL}}  \\
l_{L} & \mbox{(\ref{h1}) }\\
 \alpha_{z,+}(l) & \mbox{  (\ref{eq:AlphaBarrierDef})}\\
k_{y} & \mbox{(\ref{label.1})}\\
F^{\ast}_{L} & \mbox{(\ref{label.2})}\\

  F^{m}_{L},
  \mathcal{H}_{ m,l} & \mbox{\rm \eqref{123H4}}\\

\mathcal{B}_{ m,l}& \mbox{\rm \eqref{123}}  \\ 
\mathcal{C}_{ m,l }& \mbox{\rm \eqref{ind.34}}  \\  
\mathcal{D}_{ m,l }(j) & \mbox{\rm \eqref{ind.35}}  \\    

 \mathcal{B}_{ m,l}^{ \ga ,k}& \mbox{\rm \eqref{zdek}}  \\   
 
  \mathcal{T}_{y,\wt r_{l}}^{u_{m}, r_{l-2},n}&\mbox{\rm \eqref{eq:to showcent}}\\  


  \beta_{z} \left(l\right)  &\mbox{\rm \eqref{eq:BetaDefn2}}\\

 \al_{z, -}\left(l\right)  &\mbox{\rm \eqref{14.1}}\\
 
 \mathcal{T}_{y, l}^{ 1}, \mathcal{T}_{y,l}^{1, x^{2}} &\mbox{\rm \eqref{14.1}}\\
 
 F^{0}_{L} &\mbox{\rm \eqref{F0.def}}\\


\WW_{y,k}(n)& \mbox{(\ref{eq:LBTruncatedSumz})}\\
N_{k,a}&\mbox{\rm \eqref{nkdef.2}}\\
N_{k}, I_u& \mbox{\rm \eqref{nkdef.2j}}\\
\HH_{k,a}&
\mbox{\rm \eqref{hka.1}}\\
k^+,k^{++}& \mbox{\rm \eqref{eq:DefOfKPlus}}\\

 \wh\II_{y,z}  &\mbox{\rm \eqref{nkdef.2}}\\
 
 
 N_{k,a} &\mbox{\rm \eqref{nkdef.2}}\\

N_{k} &\mbox{\rm \eqref{nkdef.2j}}\\

\II_{y,z} 
 &\mbox{\rm \eqref{eq:TruncatedSummandLBzz}}\\
 
   \HH_{k,a} &\mbox{\rm \eqref{hka.1}}\\

J_{y,k}^{\uparrow}   &\mbox{\rm \eqref{eq:JDownz3}}\\
 
\BB_{y,k,a}   &\mbox{\rm \eqref{eq:KUpz3}}\\

   \overline \MM_{x, \ep, a,b}(n)  &\mbox{\rm \eqref{occm.n}}\\

 t_{L}\left(z\right) &\mbox{\rm \eqref{clt.2}}\\

 \overline \MM_{y,\bar \ep_{y}, y_{0}, a,b}(n)  &\mbox{\rm \eqref{unif.1a}}\\

D_{\ast}  &\mbox{\rm \eqref{dast}}\\

\end{array} $

\section{Upper bounds for excursions} 
\label{sec-upperboundex}

Let
\begin{equation}
h(r)=2\arctan ( r/2)\label{hdef}
\end{equation} 
and let
\begin{equation}
r_{l}=r_{0}e^{-l},l=0,1,\ldots,\hspace{.2 in}\mbox{and }\hspace{.2 in}h_{l}=h(r_{l})\label{dr.1}
\end{equation}
for some $r_{0}<1$. We can take $r_{0}<1$ sufficiently small that for all $0\leq x\leq r_{0}$
\begin{equation}
x-x^{3}\leq h(x) \leq x \hspace{.2 in}\mbox{and}\hspace{.2 in} |h'(x)-1|\leq x^{2}.\label{dr.1c}
\end{equation}
For some $d_{0}\leq 1/1000$ let $F_{l}$ be the centers of an $d_{0}h_{l}$ covering of $S^{2}$.
We record for future
use that  
\begin{equation}
\left|F_{l}\right|\asymp cr_{l}^{-2}=cr^{-2}_{0} e^{2l},l\ge0.\label{dr.2}
\end{equation}
Recall that $\mathcal{T}_{x, l}^{\tau }$ is the number of excursions from $\partial B_{d}\left(x,h_{l-1}\right)$ to $ \partial B_{d}\left(x,h_{l}\right)$    prior to  $ \tau$. In this section we will  assume that $2r^{ \ast}\leq h_{0}$   so that  for all $y\in B_{d}(  v,r^{ \ast})$ we have 
$B_{d}(  v,r^{ \ast})\subseteq   B_{d}( y,h_{0})$.

The reason for using $h(   r)$ is due to the following result for $\S^{ 2}$, see \cite[(2.6)]{BRZ}.  For any $u_{1}<u_{2}<u_{3}$
\begin{equation}
 \Pbm^{x\in\partial B_{d}\left(0,h (u_{2}) \right)}\left[H_{\partial B_{d}\left(0,h (u_{1})\right)}<H_{\partial B_{d}\left(0,h (u_{3})\right)}\right]=\frac{\log \(\frac{u_{2}}{u_{3}}\)}{\log \(\frac{u_{1}}{u_{3}}\)}.\label{isc.5}
\end{equation}

The next Lemma provides simple bounds which will be adequate to handle points which are    close to the `South Pole` $v$.

\bl
\label{lem:NotHitByrL}For $L $ large, any $y\in B_{d}^{c}\left(v,h_{k} \right)$ and  all $|z|\leq \log L$, 
\begin{equation}
\Pbm \left[\sqrt{2\mathcal{T}_{y,L}^{\tau}}\,\geq \rho_{L}L+z \right]\leq cke^{-2L}L e^{ -2z},\label{eq:NotHitByrL}
\end{equation}
for some $c<\ff$ independent of $1\leq k\leq L-1$.

If $y\in B_{d}\left(v,h_{L-1} \right)$
\begin{equation}
\Pbm \left[\sqrt{2\mathcal{T}_{y,L}^{\tau}}\,\geq \rho_{L}L+z \right]\leq c e^{-2L}L^{2} e^{ -2z}.\label{eq:NotHitByrLbig}
\end{equation}
\el
{\bf  Proof: } 
For $k\leq l-1$, let $ \mathcal{T}_{y,l }^{ k\to 0}$ be the number of excursions from   $\partial B_{d}\left(y,h_{l-1}\right)$ to $ \partial B_{d}\left(y,h_{l}\right)$   between $H_{\partial B_{d}\left(y,h_{k}\right)}$ and $H_{\partial B_{d}\left(y,h_{0}\right)}$. We first estimate probabilities involving $ \mathcal{T}_{y,l }^{  k\to 0}$.
Using (\ref{isc.5}), an excursion from $\partial B_{d}\left(y,h_{k}\right)$ hits $B_{d}\left(y,h_{l-1} \right)$ before exiting  $B_{d}( y,h_{0})$
has probability   $ k/(l -1  )$,   and then the  probability to hit $\partial B_{d}\left(y,h_{l}\right)$ before  exiting  $B_{d}( y,h_{0})$ is $1-\frac{1}{l }$. Thus, using the strong Markov property, 
\bea
\Pbm \left[\mathcal{T}_{y, l }^{k\to 0}\geq n\right]&=&\frac{k }{l-1 } \left(1-\frac{1}{l }\right)^{n}\label{1.1}\\ & \leq &\frac{k }{l-1 } e^{-\frac{n}{l }}, \nn
\eea
for $n$ large.  Since (recall (\ref{eq:defofts}))
\bea
(\rho_{L}L+z )^{ 2}&=&(\rho_{L}L)^{2}+2z\rho_{L}L+z^{ 2}\nn\\
&=&4L^{ 2}-4L\log L+4zL+z^{ 2}-2z  \log L+\log^{2} L,\label{eq:TsDivdedByL}
\eea
it follows that for $L$ large
\begin{equation}
\Pbm\left[\sqrt{2\mathcal{T}_{y,L}^{  k\to 0}}\,\geq \rho_{L}L+z \right]\leq cke^{-2L}L e^{ -2z}.\label{eq:NotHitByrLa}
\end{equation}
(\ref{eq:NotHitByrL}) follows since for $y\in B_{d}^{c}\left(v,h_{k} \right)$ we have $\mathcal{T}_{y,L}^{\tau}\leq \mathcal{T}_{y,L}^{ k\to 0}$.

For (\ref{eq:NotHitByrLbig}) we note that for $y\in B_{d}\left(v,h_{L-1} \right)$ we have $\mathcal{T}_{y,L}^{\tau}\leq \mathcal{T}_{y,L}^{L-1\to 0}$.
\qed

\bl
\label{lem:NotHitByrLlogL}For $L $ large and  all $0\leq z\leq \log L$, 
\begin{equation}
\Pbm \left[\sup_{y\in F_{L}\cap B_{d}\left(v,h_{\log L} \right)}\sqrt{2\mathcal{T}_{y,L}^{\tau}}\,\geq \rho_{L}L+z \right]\leq c e^{ -2z},\label{eq:NotHitByrLlogL}
\end{equation}
for some $c<\ff$.
\el

{\bf  Proof: }By Lemma \ref{lem:NotHitByrL} the probability in (\ref{eq:NotHitByrLlogL}) is bounded by 
\begin{eqnarray}
&&\sum_{k=\log L}^{\ff}\Pbm \left[\sup_{y\in F_{L}\cap B_{d}\left(v,h_{k} \right)\cap B^{c}_{d}\left(v,h_{k+1} \right)}\sqrt{2\mathcal{T}_{y,L}^{\tau}}\,\geq \rho_{L}L+z \right]
\label{eq:NotHitByrLlogLa}\\
&&\leq   \sum_{k=\log L}^{\ff} |F_{L}\cap B_{d}\left(v,h_{k} \right)\cap B^{c}_{d}\left(v,h_{k+1} \right)| cke^{-2L}L e^{ -2z}\nonumber\\
&&\leq   cL e^{ -2z}\sum_{k=\log L}^{\ff}   ke^{-2k}\leq c e^{ -2z}.\nonumber
\end{eqnarray}
\qed

Thus we only need deal with $y\in B_{d}^{c}\left(v,h_{\log L} \right)$. However,  Lemma \ref{lem:NotHitByrL} would give, for example, that 
\begin{equation}
\Pbm\left[\sup_{y\in F_{L}\cap B_{d}^{c}\left(v,h_{1} \right)}\sqrt{2\mathcal{T}_{y,L}^{ \tau}}\,\geq \rho_{L}L+z \right]\leq C L e^{ -2z},\label{122}
\end{equation}
which would be disastrous if we let $L\to\ff$. To deal with this we introduce a barrier.

 Let
 \begin{equation}
l_{L} =  l\wedge (L-l) .\label{h1}
\end{equation} 
 Fix $z$ and set
\begin{equation}
\alpha_{  z, +}\left(l\right)=\alpha\left(l,L, z\right)= \rho_{L}l +z + l_{L}^{ 1/4}.\label{eq:AlphaBarrierDef}
\end{equation}

Let
\begin{equation}
k_{y}=\inf \{k\,|\, y\in B_{d}^{c}\left(v,h_{k} \right) \},\label{label.1}
\end{equation}
and
\begin{equation}
 F^{\ast}_{L}=F_{L}\cap B_{d}^{c}\left(v,h_{\log L} \right).\label{label.2}
\end{equation}

Since $\alpha_{z, +}\left(L\right)=\rho_{L}L+z$ and $k_{y}\leq \log L$ for $y\in  F^{\ast}_{L}$,  our desired upper bound will follow from the next Lemma. 
\bl
There exists $z_{0}>0$ such that for all  $z_{0}\leq z\leq \log L$  and all $L $ large
\bea
&&
 \Pbm \left[\exists y\in  F^{\ast}_{L}, l\in\left\{k_{y}+ 1,\ldots, L\right\}\mbox{ s.t. } \mathcal{T}_{y,l}^{\tau }\ge\alpha_{z,+}^{2}\left(l\right)/2 \right]\nn\\
&&\hspace{3 in}\leq cze^{ -2z }.\label{eq:SmartMarkov20}
\eea
\el

Although this formulation looks more complicated and demanding than  our desired upper bound, it will allow us to proceed level by level and to eventually use a barrier estimate. The next Lemma will be used in our proof.

\bl\label{lem-singlelevel}For $L $ large, any $y\in B_{d}^{c}\left(v,h_{k} \right)$ and  all $0\leq z\leq \log L$, 
\be 
 \Pbm \left[ \mathcal{T}_{y,l}^{\tau }\ge\alpha_{z,+}^{2}\left(l\right)/2 \right]\leq ckl e^{-2l}e^{ -2(z+l_{L}^{ 1/4}) -(z+l_{L}^{ 1/4})^{ 2}/2l },\label{1.1act0}
\ee 
for some $c<\ff$ independent of $k\geq 1$ and $l\in\left\{k+ 1,\ldots, L\right\}$.
\el
{\bf  Proof: }As in (\ref{1.1})
\begin{equation}
 \Pbm \left[ \mathcal{T}_{y,l}^{\tau }\ge\alpha_{z,+}^{2}\left(l\right)/2 \right]\le c\frac{k }{l-1 } e^{-\frac{\alpha_{z,+}^{2}\left(l\right)}{2l }}\label{1.1act}
\end{equation} 
and 
\bea
&&\alpha_{z,+}^{2}\left(l\right)= l^{2} \rho^{ 2}_{L}+2(z+l_{L}^{ 1/4})l\rho_{L}+(z+l_{L}^{ 1/4})^{ 2}\nn\\
&&=l^{2}\(4-4{\log L \over L} +{\log^{2} L \over L^{2}}\)+2(z+l_{L}^{ 1/4})l\(2-{\log L \over L}\)\nn\\
&&\hspace{3 in}+(z+l_{L}^{ 1/4})^{ 2}\label{1.1act2}
\eea
Hence 
\bea
\alpha_{z,+}^{2}\left(l\right)/2l&=&\(2l -2\frac{l }{L }\log L\)+2(z+l_{L}^{ 1/4}) +(z+l_{L}^{ 1/4})^{ 2}/2l+o_{L}(1 )\nn\\
&&\geq \(2l -2 \log l\)+2(z+l_{L}^{ 1/4}) +(z+l_{L}^{ 1/4})^{ 2}/2l+o_{L}(1),\label{1.1act3}
\eea
using the concavity of the logarithm. Our result follows. \qed

\subsection{Proof of (\ref{eq:SmartMarkov20}) for $l$ not too large}\label{sec-begena}

\begin{proposition}
\label{prop:SmartMarkov}  There exists $z_{0}>0$ such that for all  $z_{0}\leq z\leq \log L$ and all $L $ large
\bea
&& 
 \Pbm\left[\exists y\in F^{\ast}_{L}, l\in\left\{k_{y}+ 1,\ldots, L-\(4\log L\)^{4}\right\}\right.\nn\\
&&\hspace{2.5 in}\left. \mbox{ s.t. } \mathcal{T}_{y,l}^{\tau }\ge\alpha_{z,+}^{2}\left(l\right)/2 \right]\leq ce^{ -2z }.\label{eq:SmartMarkov}
\eea
\end{proposition}

{\bf  Proof: } 
 We use a packing argument. 
Let $\phi(l)=e^{.25\,l_{L}^{ 1/4}}$.  
Considering separately the case of $ l\leq L/2$ and $L/2<l\leq L-\(4\log L\)^{4}$, we see that for some $m_{0}$
\begin{equation}
l_{L}^{ 1/4}\geq  4\log l, \hspace{.2 in}m_{0}\leq l\leq L-\(4\log L\)^{4},\label{1.2bound9}
\end{equation}
so that 
\begin{equation}
{l \over \phi(l)}={l \over e^{.25\,l_{L}^{ 1/4}}}\leq 1, \hspace{.2 in}m_{0}\leq l\leq L-\(4\log L\)^{4}.\label{1.2bound}
\end{equation}
 
We define modified radii   by
\begin{equation}
r_{l-1}^{-}=\left(1-\frac{1}{\phi(l-1)}\right)r_{l-1}\mbox{ and }r_{l}^{+}=\left(1+\frac{1}{\phi(l)}\right) r_{l}\mbox{ for }l\ge 1.\label{eq:ModifiedRadiiDef}
\end{equation}
Note that
 \be
 h(r_{l+\log \phi(l)}) \overset{\eqref{dr.1c}}{\le}   r_{l+\log \phi(l)}\overset{\eqref{dr.1}}{=}  \frac{r_{l}}{\phi(l)}\label{rel.10}
= \frac{r_{l-1}}{e \,\phi(l)}\leq  \frac{r_{l-1}}{ \phi(l-1)}.  \ee
Using this and  (\ref{dr.1c}) we have for $\phi(l)$ large enough  
\be
h(r_{l-1}^{-})+{1 \over 10^{3}}h(r_{l+\log \phi(l)})\le h(r_{l-1})\mbox{ and }h(r_{l})+{1 \over 10^{3}}h(r_{l+\log \phi(l)})\le h(r_{l}^{+}).\label{rel.11}
\ee

For each $y\in \S^{ 2}$ let $y_{l}$ denote the point in $F_{l}$
closest to $y$ (breaking ties in some arbitrary way).
By the definition   of $F_{l+\log \phi (l)}$
we have
\be 
\hspace{.2 in}d\left(y,y_{l+\log \phi(l)}\right)\le {1 \over 10^{3}}h(r_{l+\log \phi(l)}), \label{eq:Distance} 
\ee
so that using (\ref{rel.11}) we see that for all $y\in \S^{ 2}$
\begin{equation}
 B_{d}\left(y,h_{l}\right)\subset B_{d}\left(y_{l+\log \phi(l)},h(r_{l}^{+})\right)  \subset B_{d}\left(y_{l+\log \phi(l)},h(r_{l-1}^{-})\right)\subset  B_{d}\left(y,h_{l-1}\right).\label{eq:UpperBoundSandwichA}
\end{equation}

Now for $k\leq l-1$ set 
\begin{equation}
r_{ k ,l}^{-}=\left(1-\frac{1 }{\phi(l)}\right)r_{ k }\mbox{ and }r_{0,l}^{+}=\left(1+\frac{1 }{\phi(l)}\right) r_{0}\mbox{ for }l\ge0.\label{eq:ModifiedRadiiDef0}
\end{equation}

As  in the proof of (\ref{rel.11}) we have
\[
h(r_{ k ,l}^{-})+{1 \over 10^{3}}h(r_{l+\log \phi(l)})\le h(r_{ k }) 
\mbox{ and } 
h(r_{0})+{1 \over 10^{3}}h(r_{l+\log \phi(l)})\le h(r_{0,l}^{+}),
\]
so that  (\ref{eq:Distance}) also
implies that 
\begin{equation}
B_{d}\left(y_{l+\log \phi(l)},h(r_{ k ,l}^{-})\right) \subset B_{d}\left(y,h(r_{ k })\right)\subset B_{d}\left(y,h(r_{0})\right)\subset  B_{d}\left(y_{l+\log \phi(l)},h(r_{0,l}^{+})\right).\label{eq:UpperBoundSandwichB}
\end{equation}

 For each $y\in F_{l+\log \phi(l)}$ let $\wt{\mathcal{T}}_{y, l}^{k \to 0 }$ be the number of excursions from  $\partial B_{d}\left(y,h(r_{l-1}^{-})\right)$ to $ \partial B_{d}\left(y,h(r_{l}^{+})\right)$ prior to the first excursion from $\partial B_{d}\left(y, h(r_{ k ,l}^{-})\right)$ to  $\partial B_{d}\left(y, h(r_{0,l}^{+})\right)$.  Then define
 \[
\wt{\mathcal{T}}_{y,l}^{k\to 0}=\wt{\mathcal{T}}_{y_{l+\log \phi(l)},l}^{ k\to 0 },\mbox{\,\ for }y\in \S^ { 2}\backslash F_{l+\log \phi(l)}\mbox{ for all } l\ge k+1.
 \]
It follows from (\ref{eq:UpperBoundSandwichA}) and (\ref{eq:UpperBoundSandwichB}) that $\wt{\mathcal{T}}_{y, l}^{k_{y}\to 0 }\ge \mathcal{T}_{y,l}^{k_{y}\to 0}\ge \mathcal{T}_{y,l}^{\tau}$ for all $l\ge k_{y}+1$. Thus

\bl
\label{lem:TraversalPackingDomination}For all $y\in \S^{ 2},l\ge k_{y} +1$
  we have that
\begin{equation}
\wt{\mathcal{T}}_{y,l}^{k_{y}\to 0}\ge \mathcal{T}_{y,l}^{\tau}.\label{eq:TraversalPackingDominationk}
\end{equation}
\el

Because of this    the probability
in (\ref{eq:SmartMarkov}) is bounded above by 
\bea
&&
{\displaystyle \sum_{k= 1}^{ \log L}}\,\,\,{\displaystyle \sum_{l=k +1}^{L-\(4\log L\)^{4}}}{\displaystyle \sum_{y\in B_{d}\left(v,h(r_{ k-1 })\right)\cap  F_{l+\log \phi(l)}}  }\Pbm \left[\wt{\mathcal{T}}_{y, l}^{k\to 0 }\ge\alpha_{z,+}^{2} \left(l\right)/2    \right] \label{2.28f}\\ & =&{\displaystyle \sum_{k= 1}^{ \log L}}\,\,\,{\displaystyle \sum_{l=k +1}^{L-\(4\log L\)^{4}}}{\displaystyle \left|B_{d}\left(v,h(r_{ k-1 })\right)\cap F_{l+\log \phi(l)}\right|}\Pbm \left[\wt{\mathcal{T}}_{y, l}^{k\to 0  }\ge\alpha_{z,+}^{2} \left(l\right)/2     \right] \nn\\
 & \leq & c {\displaystyle \sum_{k= 1}^{ \log L}}\,\,\,{\displaystyle \sum_{l=k +1}^{L-\(4\log L\)^{4}}}e^{.5\,l_{L}^{ 1/4}}e^{2(l-k)}\Pbm \left[\wt{\mathcal{T}}_{y, l}^{k\to 0  }\ge\alpha_{z,+}^{2} \left(l\right)/2     \right] ,\nn
\eea
for some arbitrary $y\in F_{l+\log \phi(l)}$.  We show below that for all $k\leq l$
\begin{equation}
\Pbm \left[\wt{\mathcal{T}}_{y, l}^{k\to 0  }\geq \alpha_{z,+}^{2}(l)/2 \right]\leq ce^{-2l-  l_{L}^{ 1/4}}e^{ -2z }, \label{goalub1.j}
\end{equation}
and since 
\begin{equation}
{\displaystyle \sum_{k= 1}^{ \ff}}\,\,\,{\displaystyle \sum_{l= 1}^{L }}e^{.5\,l_{L}^{ 1/4}}e^{2(l-k)}e^{-2l-  l_{L}^{ 1/4}}e^{ -2z }\leq ce^{ -2z },\nn
\end{equation}
this will complete the proof of (\ref{eq:SmartMarkov}). 

We now turn to the proof of (\ref{goalub1.j}).
Let
\bea
p_{l}&=& \frac{\log\left(r_{l-1}^{-}/r_{0,l}^{+}\right)}{\log\left(r_{l}^{+}/r_{0,l}^{+}\right)}
=\frac{\log\left(\left(1-\frac{1}{\phi(l-1)}\right)\left(1+\frac{1}{\phi(l)}\right)^{ -1}e^{ -(   l-1) }\right)}{\log\left(e^{ -l }\right)}\nn
\\&
=&\frac{l -1+2 /\phi(l)+O( \phi(l)^{ -2}   )}{l }=1-\frac{1-2/\phi(l)+O( \phi(l)^{- 2}   )}{l }, \label{peel.1}
\eea
and
\bea
q_{l}&=& \frac{\log\left(r_{k,l}^{-}/r_{0,l}^{+}\right)}{\log\left(r_{l-1}^{-}/r_{0,l}^{+}\right)}
=\frac{\log\left(\left(1-\frac{1}{\phi(l)}\right)\left(1+\frac{1}{\phi(l)}\right)^{ -1}e^{ -k }\right)}{\log\left(\left(1-\frac{1}{\phi(l-1)}\right)\left(1+\frac{1}{\phi(l)}\right)^{ -1}e^{ -(   l-1) }\right)}\nn
\\&
=&\frac{k+2/\phi(l)+O( \phi(l)^{ -2}   )}{l -1+2/\phi(l)+O( \phi(l)^{ -2}   )}=\frac{k+O( \phi(l)^{- 1}  )}{l -1}.\nn
\eea
Using the fact that $p_{l}<1$ together with (\ref{peel.1}) we can write
\begin{equation}
p_{l}=1-\frac{1-b_{l}/\phi(l)}{l }\label{peel.2}
\end{equation}
with $1-b_{l}/\phi(l)>0$. In addition, using (\ref{1.2bound}) and possibly increasing $m_{0}$,
\begin{equation}
\frac{lb_{l}}{\phi(l)}\leq 3,\qquad l\geq m_{0}.\label{peel.2a}
\end{equation}

Since $q_{l}$ is the probability for an excursion from $\partial B_{d}\left(y, h(r_{k,l}^{-})\right)$ to  hit $B_{d}\left(y, h(r_{l-1}^{-})\right)$ before $\partial B_{d}\left(y,h(r_{0,l}^{+})\right)$, and $p_{l}$ is  the  probability for an excursion from $B_{d}\left(y, h(r_{l-1}^{-})\right)$ to hit $\partial B_{d}\left(y, h(r_{l}^{+})\right)$ before $\partial B_{d}\left(y,h(r_{0,l}^{+})\right)$,
we see that as in (\ref{1.1}) 
\be
\Pbm \left[\wt{\mathcal{T}}_{y, l}^{k\to 0  }\geq \alpha_{z,+}^{2}(l)/2 \right]\leq \frac{ck}{l} e^{-\frac{\alpha_{z,+}^{2}(l)}{2l }\left( 1-\frac{b_{l}}{\phi(l)}  \right)}. \label{1.2axz}
\ee
By (\ref{1.1act3})
we have that
\bea
{\alpha_{z,+}^{2}(l)\over 2l }
&\geq&    \(2l -2 \log l\)+2(z+l_{L}^{ 1/4}) +z^{ 2}/2l+o_{L}(1),\label{1.3v}
\eea
so that, for $k\leq l$
\be
\Pbm \left[\wt{\mathcal{T}}_{y, l}^{k\to 0  }\geq \alpha_{z,+}^{2}(l)/2 \right]\leq c  l^{2} e^{-\(2l+2(z+l_{L}^{ 1/4}) +z^{ 2}/2l\)\left( 1-\frac{b_{l}}{\phi(l)}  \right)}. \label{1.2axzbr}
\ee
We claim that
\begin{equation}
{z^{ 2} \over 2l }\left( 1-\frac{b_{l}}{\phi(l)}  \right) -2z b_{l} /\phi(l)>0,\nn
\end{equation}
that is
\begin{equation}
z\left( 1-\frac{b_{l}}{\phi(l)}  \right) >4l b_{l} /\phi(l) \nn
\end{equation} 
for $z\geq z_{0}$ sufficiently large.  For $l>m_{0}$, this follows from  (\ref{peel.2a}), and for $l\leq m_{0}$ we can just increase $z$
further.
Thus for such $z$
\be
\Pbm \left[\wt{\mathcal{T}}_{y, l}^{k\to 0  }\geq \alpha_{z,+}^{2}(l)/2 \right]\leq cl^{2} e^{-\(2l+2l_{L}^{ 1/4}  \)\left( 1-\frac{b_{l}}{\phi(l)}  \right)}e^{-2z}. \label{1.2axzb}
\ee
For $k\leq l\leq m_{0}$ this already proves (\ref{goalub1.j}) with $c$ sufficiently large. 
For $l>m_{0}$, using  (\ref{peel.2a}) again we now have
\be
\Pbm \left[\wt{\mathcal{T}}_{y, l}^{k\to 0  }\geq \alpha_{z,+}^{2}(l)/2 \right]\leq cl^{2} e^{-\(2l+2 l_{L}^{ 1/4}  \) }e^{-2z}. \label{1.2axzbs}
\ee
and (\ref{1.2bound9}) completes the proof of (\ref{goalub1.j}).
\qed

\subsection{Proof of (\ref{eq:SmartMarkov20}) for $l$ very large}\label{sec-begen}

It suffices to show that for some small but fixed constant $\tilde{c}$ to be chosen later 
we have that for all  $L $ sufficiently  large   and all $z_{0} \leq z\leq \log L$
\be
\mathbb{P}\left[\begin{array}{c}
\exists y\in F^{\ast}_{L}\cap B_{d}\left(0,\tilde{c}h_{0}\right)\mbox{ and }k_{y}+1\le l\le L\\
\mbox{ such that }\sqrt{2\mathcal{T}_{y,l}^{ \tau}}\ge \al_{z,+}(l)
\end{array}\right]\le cze^{-2z }. \label{g17.0ab}
\ee
Here $0$, the center of $B_{d}\left(0,\tilde{c}h_{0}\right)$, is used to denote an arbitrary point in $\S^{2}$.

Now consider
\[
\mathcal{G}_{l}=\left\{ \sqrt{2\mathcal{T}_{y,l'}^{\tau}}\le \al_{z,+}(l')\mbox{ for all }l'=k_{y}+1,\ldots,l\mbox{ and }\forall y\in F^{\ast}_{L}\cap B_{d}\left(0,\tilde{c}h_{0}\right)\right\} .
\]
Let $L'=L-\left(4 \log L\right)^{4}$. With
\begin{equation}
\mathcal{H}_{l}=\left\{ \exists y\in F^{\ast}_{L}\cap B_{d}\left(0,\tilde{c}h_{0}\right)\mbox{ s.t. }\sqrt{2\mathcal{T}_{y,l}^{ \tau}}\ge \al_{z,+}(l),\,k_{y}<l\right\},
\label{123H}
\end{equation}
we will prove that for all $l>L'$
\begin{equation}
\Pbm\left[\mathcal{H}_{l} \cap\mathcal{G}_{l-2} \right]\le cz e^{-  l_{L}^{ 1/4}-2z},\label{geq: to prove Gl}
\end{equation}
so that  we have  
\bea
\Pbm\left[\mathcal{G}_{L}^{c} \right]&&\le\sum_{l=L'+1}^{L}\Pbm\left[\mathcal{G}_{l}^{c}\cap\mathcal{G}_{l-1} \right]+\Pbm\left[\mathcal{G}^{c}_{L'}  \right]\nn\\
&&\le\sum_{l=L'+1}^{L}\Pbm\left[\mathcal{H}_{l}\cap\mathcal{G}_{l-2} \right] +\Pbm\left[\mathcal{G}^{c}_{L'}  \right]\nn\\
&&\le\sum_{l=L'+1}^{L}cz e^{-  l_{L}^{ 1/4}-2z }+\Pbm\left[\mathcal{G}^{c}_{L'}  \right] \nn\le cze^{-2z },\label{cbsum}
\eea
by Proposition \ref{prop:SmartMarkov}, which  will prove (\ref{g17.0ab}).

Setting $F^{m}_{L}=F_{L}\cap B_{d}^{c}\left(v,h_{m} \right)\cap B_{d}\left(v,h_{m-1} \right)$, so that  $k_{y}=m$ for $y\in F^{m}_{L}$, and for any $l>m$
\begin{equation}
\mathcal{H}_{ m,l}=\left\{ \exists y\in F^{m}_{L}\cap B_{d}\left(0,\tilde{c}h_{0}\right)\mbox{ s.t. }\sqrt{2\mathcal{T}_{y,l}^{ \tau}}\ge \al_{z,+}(l) \right\},
\label{123H4}
\ee
we will prove that for all $l>L'$
\begin{equation}
\sum_{m=1}^{\log L}\Pbm\left[\mathcal{H}_{ m,l} \cap\mathcal{G}_{l-2} \right]\le cz e^{-  l_{L}^{ 1/4}-2z},\label{geq: to prove Gl4}
\end{equation}
 which gives (\ref{geq: to prove Gl}).

For any $l>m$ let 
\begin{equation}
\mathcal{B}_{ m,l}=\left\{ \exists x\in F^{m}_{L}\cap B_{d}\left(u_{m},\tilde{c}h_{l}\right)\mbox{ s.t. }\sqrt{2\mathcal{T}_{x,l}^{ \tau}}\ge \al_{z,+}(l) \right\},\label{123}
\end{equation}
where we can assume that $u_{m}\in F^{m}_{L}$.
By a union bound, $\Pbm\left[\mathcal{H}_{ m,l} \cap\mathcal{G}_{l-2} \right]$
is bounded above by
\be
ce^{2(l-m)}\times \Pbm\left[ \mathcal{B}_{ m,l}\cap\mathcal{G}_{l-2} \right].\label{geq: union 1}
\ee 
Hence it suffices to show that
\begin{equation}
\sum_{m=1}^{\log L}e^{-2m}\Pbm\left[\mathcal{B}_{ m,l} \cap\mathcal{G}_{l-2}  \right]\le cze^{-2l-  l_{L}^{ 1/4}-2z}.\label{nod.12}
\end{equation}

 Since $\mathcal{G}_{l-2}\subset\mathcal{C}_{ m,l}$, where 
\be
\mathcal{C}_{ m,l }=\left\{ \sqrt{2\mathcal{T}_{u_{m}, l'}^{ \tau}}\le \al_{z,+}(l')\mbox{ for all }l'=m+1,\ldots,l-2\right\},\label{ind.34}
\ee
it suffices to show that
\begin{equation}
\sum_{m=1}^{\log L}e^{-2m}\Pbm\left[\mathcal{B}_{ m,l}  \cap\mathcal{C}_ {m,l }  \right]\le cze^{-2l- l_{L}^{ 1/4}-2z}.\label{gsuf1}
\end{equation}

We show in the next Section that  for all $l\geq L-\left(4 \log L\right)^{4}$
\begin{equation}
\Pbm\left[\mathcal{B}_{ m,l}  \cap\left\{\sqrt{2\mathcal{T}_{u_{m}, l-2}^{ \tau}}\le {1 \over 2}\al_{z,+}(l-2)\right\} \right]\le ce^{-c'L^{2}}.\label{notSmall}
\end{equation}
It follows from this that with
\be
\mathcal{D}_{ m,l }(j)=\left\{ \sqrt{2\mathcal{T}_{u_{m}, l-2}^{ \tau}}\in I_{\al_{z,+}(l-2)+j}\right\},\label{ind.35}
\ee
 it suffices to show that 
\begin{eqnarray}
&&\sum_{m=1}^{\log L} e^{-2m}\sum_{j=0}^{{1 \over 2}\al_{z,+}(l-2)}\Pbm\left[\mathcal{B}_{ m,l} \cap\mathcal{C}_{ m,l } \cap\mathcal{D}_{ m,l }(-j)\right] \le cze^{-2l- l^{1/4}_{L}-2z}.\label{gcb.9dec}
\end{eqnarray}

We also show in the next Section that we can find a fixed $j_{0}$ such that for all $j_{0}\leq j\leq {1 \over 2}\al_{+,z}(l-2)$, uniformly in $1\leq m\leq \log L$ and $z_{0}\leq z\leq \log L$, for any 
$\wt{ \mathcal{C}}_{ m,l }\in  { \mathcal{F}}\(\mathcal{T}_{u_{m}, k}^{ \tau}, k=1,\ldots, l-2\)$  
\begin{equation}
\Pbm \left[\mathcal{B}_{ m,l} \,\big |\, \wt{ \mathcal{C}}_{ m,l} \cap\mathcal{D}_{ m,l }(-j)\right]      \le Ce^{-4j },\label{cb.15y}
\end{equation} 
 by taking $\tilde{c}>0$ sufficiently small.

It follows from the barrier estimate  (\ref{18.20}) that for $0\leq   j\leq {1 \over 2}\al_{z,+}(l-2)$,
\begin{equation}
\Pbm\left[ \mathcal{C}_{ m,l } \cap\mathcal{D}_{ m,l }(-j)\right]\leq  ce^{-2l-2z-2l^{1/4}_{L}+2j } \times m^{2}\left(1+z+m + l^{ 1/4}_{L}\right)\left(1+j  \right).\label{cb.decp7}
\end{equation}
Combining the last 2 displays   we can bound the left hand side of (\ref{gcb.9dec}) by
\begin{eqnarray}
&&\hspace{.3in} Cze^{-2l-2z- l^{1/4}_{L}} \sum_{j=0}^{\ff}e^{-4j 1_{\{j\geq j_{0}\}}+2j}
(1+j),\label{gcb.9dec3} 
\end{eqnarray}
which   proves (\ref{gcb.9dec}).

\subsection{Proof of the conditional barrier estimate  (\ref{cb.15y}) and  (\ref{notSmall})}\label{sec-beearly2}

We first prove that for some $j_{0}$ fixed and all $j_{0}\leq j\leq {1 \over 2}\al_{+,z}(l)$,
\begin{equation}
\Pbm\left[\mathcal{B}_{ m,l} \,\big |\, \wt{ \mathcal{C}}_{ m,l } \cap\mathcal{D}_{ m,l }(-j)\right]  \le Ce^{-4j }.\label{para.1}
\end{equation}

{\bf Proof:   } 
For each $\ga \in(0,1]$ and $y$, let $\mathcal{T}_{y,\wt r_{l}}^{ \tau }$ be the number of excursions from  $\partial B\left(y,h( \wt r_{l-1})\right)$ to $ \partial B\left(y,h(\wt r_{l})\right)$    prior to  $ \tau$,
where
\begin{equation}
\wt r_{l-1}=r_{l-1}\left(1-\ga \right) ,\hspace{.2 in}\wt r_{l}=r_{l}\left(1+\ga \right).\label{scale2}
\end{equation} 
Note that  
\be
\mathcal{T}_{y',\wt r_{l}}^{\tau}\ge \mathcal{T}_{y,l}^{\tau}\mbox{ for all } y'\mbox{ such that }d\left(y,y'\right)\le{\ga  r_{l} \over 2},\label{cent1}
\ee
since then 
\[
B_{d}\left(y,h_{l-1}\right)\supset B_{d}\left(y',h(r_{l-1}\left(1-\ga \right))\right)\supset  B_{d}\left(y',h(r_{l}\left(1+\ga \right))\right)\supset B_{d}\left(y,h_{l}\right).
\] 
 
Let
\be
 \mathcal{B}_{ m,l}^{ \ga ,k} =\left\{ \exists y\in F^{m}_{k}\cap B_{d}\left(u_{m},\tilde{c}\,h_{l}\right)\mbox{ such that }\sqrt{2\mathcal{T}_{y,\wt r_{l}}^{\tau}}\ge \al_{+,z}(l)\right\}.\label{zdek}
\ee
Note $F^{m}_{k}$ not $F^{m}_{L}$.
From now on we fix
\be
\ga =\frac{1}{\al_{+,z}(l)-j},\hspace{.2 in}\mbox{and }\hspace{.2 in}k=\log\left(2(\al_{+,z}(l)-j)\right)+l.\label{conv-k}
\ee
We will show that with these values
\begin{equation}
\Pbm\left[\mathcal{B}_{ m,l}^{  \ga ,k} \,\bigg|\, \wt{ \mathcal{C}}_{ m,l} \cap\mathcal{D}_{ m,l}(-j)\right]      \le Ce^{-4j}.\label{eq:to show}
\end{equation}
Using (\ref{cent1}) this will imply (\ref{para.1}), since  
for each $y\in F^{m}_{L}\cap B_{d}\left(u_{m},\tilde{c}h_{l}\right)$ there exists  a representative $y'\in F^{m}_{k}\cap B_{d}\left(u_{m},\tilde{c}h_{l}\right)$
such that 
\[
d\left(y,y'\right)\le r_{k}=\frac{1}{2(\al_{+,z}(l)-j)}\,\,r_{l}={\ga  r_{l} \over 2}.
\]

To show (\ref{eq:to show}), we first show that for some $c_{3}>0$
\be
\Pbm\left[\sqrt{2\mathcal{T}_{u_{m},\wt r_{l}}^{\tau}} \ge\al_{+,z}(l)-\frac{j}{2}\,\bigg|\, \wt{ \mathcal{C}}_{ m,l} \cap\mathcal{D}_{ m,l}(-j)\right]  \leq c'e^{-c_{3}j^{2}}.\label{eq:to showcent}
\ee
Let $\mathcal{T}_{y,\wt r_{l}}^{u_{m}, r_{l-2},n} $  be the number of excursions from  $\partial B\left(y,h( \wt r_{l-1})\right)$ to $ \partial B\left(y,h(\wt r_{l})\right)$  during the first $n$  excursions from  $\partial B\left(u_{m},h_{l-2}\right)$ to $ \partial B\left(u_{m},h_{l-3}\right)$. Using the Markov property we have that 
\begin{eqnarray}
&&
\Pbm\left[\sqrt{2\mathcal{T}_{u_{m},\wt r_{l}}^{\tau}} \ge\al_{+,z}(l)-\frac{j}{2}\,\bigg|\, \wt{ \mathcal{C}}_{ m,l} \cap\mathcal{D}_{ m,l}(-j)\right]\label{bu1}\\
&=& \Pbm\left[\sqrt{2\mathcal{T}_{u_{m},\wt r_{l}}^{\tau}} \ge\al_{+,z}(l)-\frac{j}{2}\,\bigg|\,\mathcal{D}_{ m,l}(- j)\right]\nn\\
&=&
\Pbm\left[\sqrt{2\mathcal{T}_{u_{m},\wt r_{l}}^{\tau}} \ge\al_{+,z}(l)-\frac{j}{2}\,\bigg|\, \sqrt{2\mathcal{T}_{u_{m},l-2}^{\tau}}\in I_{\al_{z,+}(l)-j}\right]\nn\\
\hspace{-.3 in}&=& \Pbm\left[ 
\sqrt{2\mathcal{T}_{u_{m},\wt r_{l}}^{u_{m}, r_{l-2},\mathcal{T}_{u_{m},l-2}^{\tau}}}\ge\al_{+,z}(l) -\frac{j}{2}\,\bigg|\, \sqrt{2\mathcal{T}_{u_{m},l-2}^{\tau} }\in I_{\al_{z,+}(l)-j}\right]. \nn
\end{eqnarray}

To prove (\ref{eq:to showcent}) it suffices to show that show  that uniformly  for $s\in I_{\al_{z,+}(l)-j}$
\begin{equation}
\mathbb{P}\left[ 
\sqrt{2\mathcal{T}_{u_{m},\wt r_{l}}^{u_{m}, r_{l-2},s^{2}/2}}\ge\al_{+,z}(l)-\frac{j}{2} \right]\leq c'e^{-c_{3}j^{2}}.\label{bu2}
\end{equation}

To see this, let $s=\al_{z,+}(l )-j+\ze$, where $0\leq \ze\leq 1$.
It follows from  \cite[Lemma 4.6]{BK},  that
\begin{equation}
\mathbb{P}\left[ 
\sqrt{2\mathcal{T}_{u_{m},\wt r_{l}}^{u_{m}, r_{l-2},s^{2}/2}}\ge\al_{+,z}(l)-\frac{j}{2} \right]\leq e^{-\left(\sqrt{q} (\al_{+,z}(l)-j+ \ze) -\sqrt{p}(\al_{+,z}(l)-\frac{j}{2})\right)^{2}/2}\label{bu2a}
\end{equation}
where
\bea
q&=&\Pbm^{u\in\partial B_{d}\left(u_{m},h_{l-2}\right)}\left[H_{\partial B_{d}\left(u_{m},h(r_{l}\left(1+\ga \right))\right)}<H_{\partial B_{d}\left(u_{m},h_{l-3}\right)}\right]\nn\\
&=&\frac{\log r_{l-3}-\log r_{l-2}}{\log r_{l-3}-\log (r_{l}\left(1+\ga \right))}=\frac{1}{3+O\left(\ga \right)},
\eea
and
\bea
p&=&\Pbm^{u\in\partial B_{d}\left(u_{m},h(r_{l-1}\left(1-\ga \right))\right)}\left[H_{\partial B_{d}\left(u_{m},h_{l-3}\right)}<H_{\partial B_{d}\left(u_{m},h(r_{l}\left(1+\ga \right))\right)}\right]\nn\\
&=&\frac{\log (r_{l-1}\left(1-\ga \right))-\log (r_{l}\left(1+\ga \right))}{\log r_{l-3}-\log (r_{l}\left(1+\ga \right))}=\frac{1+O\left(\ga \right)}{3+O\left(\ga \right)}.
\eea

To apply \cite[Lemma 4.6]{BK} it suffices to show  that 
\be
\al_{z,+}(l)-\frac{j}{2}\geq (\al_{+,z}(l )-j+1)\sqrt{q/p},
\nn
\ee
and the right hand side
\be
 =(\al_{+,z}(l)-j+1)(1+O\left(\ga \right))=(\al_{+,z}(l)-j)+O\left(1\right),  \nonumber
\ee
since $\ga \left(\al_{+,z}(l)-j\right)=1$. Thus we can use (\ref{bu2a}) for all $j\geq c_{3}$ for some $c_{3}<\ff$. For such $j$
we therefore have
\[
\mathbb{P}\left[ 
\sqrt{2\mathcal{T}_{u_{m},\wt r_{l}}^{u_{m}, r_{l-2},s^{2}/2}}\ge\al_{+,z}(l)-\frac{j}{2} \right]\le ce^{ -\frac{1}{6}\left(-\frac{j}{2}+\ze+O\left(\ga \left(\al_{+,z}(l)-\frac{j}{2}\right)\right)\right)^{2}}, 
\]
and since $\ga \left(\al_{+,z}(l)-j\right)=1$ and $j\leq {1 \over 2}\al_{+,z}(l)$ so that $j\leq \left(\al_{+,z}(l)-j\right)$, it follows that 
\[\ga \left(\al_{+,z}(l)-\frac{j}{2}\right)=\ga \left(\al_{+,z}(l)-j\right)+\ga\frac{j}{2}\leq 2,
\]
so that we obtain  (\ref{bu2}) for all $j\geq c_{3}$. By enlarging $c'$ we then have (\ref{bu2}) for all $j$.

We now  bound
\bea
&&\Pbm\left[\mathcal{B}_{m,l}^{  \ga ,k} \,\bigg|\, \wt{ \mathcal{C}}_{m,l} \cap\mathcal{D}_{m,l}(-j)\right] \le \Pbm\left[\sqrt{2\mathcal{T}_{u_{m},\wt r_{l}}^{\tau}} \ge\al_{+,z}(l)-\frac{j}{2}\,\bigg|\,\wt{ \mathcal{C}}_{m,l} \cap\mathcal{D}_{m,l}(-j)\right]\nn\\
&&\hspace{1 in} +\Pbm\left[\mathcal{B}_{m,l}^{ \ga ,k} \cap \Big\{\sqrt{2\mathcal{T}_{u_{m},\wt r_{l}}^{\tau}}<\al_{+,z}(l)-\frac{j}{2}\Big\}\,\bigg|\,\wt{ \mathcal{C}}_{m,l} \cap\mathcal{D}_{m,l}(-j)\right].
\eea
Because of the bound (\ref{eq:to showcent}), to prove (\ref{eq:to show}) it suffices to
show that
\begin{equation}
\Pbm\left[\mathcal{B}_{m,l}^{ \ga ,k} \cap \Big\{\sqrt{2\mathcal{T}_{u_{m},\wt r_{l}}^{\tau}}<\al_{+,z}(l)-\frac{j}{2}\Big\}\,\bigg|\,\wt{ \mathcal{C}}_{m,l} \cap\mathcal{D}_{m,l}(-j)\right]\le Ce^{-4j}.\label{eq:to show2}
\end{equation}

  Assign to each $y\in F^{m}_{l+i}\cap B_{d}\left(u_{m},\tilde{c}h_{l}\right)$
a unique ``parent'' $\tilde{y}\in F^{m}_{l+i-1}\cap B_{d}\left(u_{m},\tilde{c}h_{l}\right)$ such that
$d\left(\tilde{y},y\right)\le r_{l+i}$. In particular, for $i=1$ we set $\tilde{y}=u_{m}$. Let $q=q( \tilde{y},y  )=d\left(\tilde{y},y\right)/r_{l}$ and set
\begin{equation}
A_{i}=\Big\{\underset{y\in F^{m}_{l+i}\cap B_{d}\left(u_{m},\tilde{c}h_{l}\right)}{\sup}\left|\mathcal{T}_{y,\wt r_{l}}^{\tau}-\mathcal{T}_{\tilde{y},\wt r_{l}}^{\tau}\right|\le d_{0}j i\left(\al_{+,z}(l)-j\right)\sqrt{q}\Big\},\label{s1.1}
\end{equation}
where $d_{0}$ will be chosen later, but small enough that $d_{0}\sum_{i\ge1}ie^{-i/2}\le\frac{1}{8}$.
We now show that 
\bea
&&
\bigcap_{i=1}^{k-l}A_{i}\cap \Big\{\sqrt{2\mathcal{T}_{u_{m},\wt r_{l}}^{\tau}}<\al_{+,z}(l)-\frac{j}{2}\Big\}\label{s1.2}\\
&&\subseteq \Big\{\sqrt{2\mathcal{T}_{y,\wt r_{l}}^{\tau}}<\al_{+,z}(l),\,\forall y\in F^{m}_{k}\cap B_{d}\left(u_{m},\tilde{c}h_{l}\right)\Big\}. \nn
\eea
For this, we use $q=d\left(\tilde{y},y\right)/r_{l}\le r_{l+i}/r_{l}=e^{-i}$ for $y\in F^{m}_{l+i}$ to see that for any trajectory in the left hand side of (\ref{s1.2}) and all $y\in F^{m}_{k}\cap B_{d}\left(u_{m},\tilde{c}h(r_{l})\right)$
\[
\mathcal{T}_{y,\wt r_{l}}^{\tau}\le\left(\al_{+,z}(l)-\frac{j}{2}\right)^{2}/2+j\left(\al_{+,z}(l)-j\right)d_{0}\sum_{i\ge1}ie^{-i/2},
\]
which, since $d_{0}\sum_{i\ge1}ie^{-i/2}\le\frac{1}{8}$,
implies that
\[
\begin{array}{ccl}
\mathcal{T}_{y,\wt r_{l}}^{\tau} & \le & \left(\al_{+,z}(l)-\frac{j}{2}\right)^{2}/2+\frac{1}{8}j\left(\al_{+,z}(l)-j\right)\\
 & = & \al^{2}_{+,z}(l)/2 -\al_{+,z}(l)j/2+(\frac{j}{2})^{2}/2+\frac{1}{8}\al_{+,z}(l)j-\frac{j^{2}}{8}\\
 & < & \al^{2}_{+,z}(l)/2.
\end{array}
\]
This establishes (\ref{s1.2}) and taking complements we see that
\bea
&&
 \mathcal{B}_{m,l}^{\ga ,k} =\Big\{\exists y\in F^{m}_{k}\cap B_{d}\left(u_{m},\tilde{c}h_{l}\right)\mbox{ such that }\sqrt{2\mathcal{T}_{y,\wt r_{l}}^{\tau}}\ge\al_{+,z}(l)\Big\}\label{s1.3}\\
&&\hspace{1.5 in}\subseteq \bigcup_{i=1}^{k-l}A^{c}_{i}\cup \Big\{\sqrt{2\mathcal{T}_{u_{m},\wt r_{l}}^{\tau}}\geq \al_{+,z}(l)-\frac{j}{2}\Big\}.\nn
\eea
It follows that
\begin{eqnarray}
&& \mathcal{B}_{m,l}^{\ga ,k} \cap \Big\{\sqrt{2\mathcal{T}_{u_{m},\wt r_{l}}^{\tau}}<\al_{+,z}(l)-\frac{j}{2}\Big\}
\subseteq \bigcup_{i=1}^{k-l}A^{c}_{i}. \label{gamlab}
\end{eqnarray}

We can thus bound $\Pbm\left[\mathcal{B}_{m,l}^{ \ga ,k} \cap \Big\{\sqrt{2\mathcal{T}_{u_{m},\wt r_{l}}^{\tau}}<\al_{+,z}(l)-\frac{j}{2}\Big\}\,\bigg|\,\wt{ \mathcal{C}}_{m,l} \cap\mathcal{D}_{m,l}(-j)\right]$
by
\begin{equation}
\sum_{i=1}^{k-l}\mathbb{P}\left[
\underset{y\in F^{m}_{l+i}\cap B_{d}\left(u_{m},\tilde{c}h_{l}\right)}{\sup}\left|\mathcal{T}_{y,\wt r_{l}}^{\tau}-\mathcal{T}_{\tilde{y},\wt r_{l}}^{\tau}\right|\ge d_{0}j i\left(\al_{+,z}(l)-j\right)\sqrt{q}\,\bigg|\,\wt{ \mathcal{C}}_{m,l} \cap\mathcal{D}_{m,l}(-j)\right].\label{eq:msum}
\end{equation}
Since $\left|F^{m}_{l+i}\cap B_{d}\left(u_{m},\tilde{c}h_{l}\right)\right|\leq c e^{2i}$,
a union bound   gives
that (\ref{eq:msum}) is at most
\bea
&&\hspace{.2 in} c\sum_{i=1}^{k-l}e^{2i}\underset{y\in F^{m}_{l+i}\cap B_{d}\left(u_{m},\tilde{c}h_{l}\right)}{\sup}\nn\\
&&\hspace{1 in}\mathbb{P}\left[
\left|\mathcal{T}_{y,\wt r_{l}}^{\tau}-\mathcal{T}_{\tilde{y},\wt r_{l}}^{\tau}\right|\ge d_{0}j i\left(\al_{+,z}(l)-j\right)\sqrt{q}\,\bigg|\,\wt{ \mathcal{C}}_{m,l} \cap\mathcal{D}_{m,l}(-j)\right].\label{17.1ax}
\eea

We can write the last probability as
\begin{equation}
\mathbb{P}\left[
\left|\mathcal{T}_{y,\wt r_{l}}^{u_{m}, r_{l-2},\mathcal{T}_{u_{m},l-2}^{\tau}}-\mathcal{T}_{\tilde{y},\wt r_{l}}^{u_{m}, r_{l-2},\mathcal{T}_{0,l-2}^{\tau}}\right|\ge d_{0}j i\left(\al_{+,z}(l)-j\right)\sqrt{q}\,\bigg|\,\wt{ \mathcal{C}}_{m,l}\cap\mathcal{D}_{m,l}(-j)\right].\label{eq-thanks4lem}
\end{equation}
Using  \cite[Lemma 5.6]{BRZ} with $\th=d_{0}j i$ and $n=(\al_{+,z}(l)-j)^{2}/2$,  we find that for an appropriate choice of $d_{0}, \tilde{c}$, the last probability is bounded by $ ce^{-4(j+ i)} $.
To apply \cite[Lemma 5.6]{BRZ} we must verify several points. 

First, we need to verify that for some small $\bar c_{0}$ we have $\th\leq \bar c_{0} (n-1)$, that  is 
$  d_{0}j i\leq \bar c'_{0}(\al_{+,z}(l)-j)^{2}$. For this it suffices to note that for $j,l$ in our range  $i/(\al_{+,z}(l) -j)\leq (k-l)/(\al_{+,z}(l) -j)=(\log 2(\al_{+,z}(l) -j))/(\al_{+,z}(l) -j)$ goes to $0$ as $L\to\ff$.

Secondly,  we need to show that $  \theta \le  ((n-1) q)^{2} $.
Since we have already seen that  $ \theta\leq\bar  c_{0}(n-1)$, it suffices to show that    $(n-1) q^{2}\geq   c^{2}_{2}$ for some $c_{2}>0$, or equivalently that $\sqrt{2n}\,\,    q\geq   c'_{2}>0$. That  is,
$(\al_{+,z}(l) -j) d\left(\tilde{y},y\right)/r_{l}  \geq c'_{2}$.     Assume that $d\left(\tilde{y},y\right)\ge c_{3}r_{k}$
for a small $c_{3}>0$, so that, see (\ref{conv-k}),
\begin{equation}
(\al_{+,z}(l) -j) d\left(\tilde{y},y\right)/r_{l}\ge c_{3}(\al_{+,z}(l) -j)e^{-\left(k-l\right)}=c_{3}/2.\nn
\end{equation}
With the $F_{l}$ constructed appropriately we can indeed assume that
either  $d\left(\tilde{y},y\right)\ge c_{3}r_{k}$
for a small $c_{3}>0$, or that $y=\tilde{y}$, in which case the corresponding
term in the sum in (\ref{17.1ax}) is zero. Also, by taking $\tilde{c}=q_{0}/2$ we will have $d\left(\tilde{y},y\right)/r_{l}\leq q_{0}$.

Thus  we see
that (\ref{17.1ax}) is at most
\be
c\sum_{i=1}^{k-l}e^{2i}e^{-4(j+ i)}\le Ce^{-4j}.\label{dec1a}
\ee 
This completes the proof of (\ref{eq:to show2}).
\qed

 {\bf  Proof of  (\ref{notSmall}): } As in (\ref{cent1})
 \be
\mathcal{T}_{u_{m},\wt r_{l}}^{\tau}\ge \mathcal{T}_{y,l}^{\tau}\mbox{ for all } y\mbox{ such that }d\left(y,u_{m}\right)\le{\ga  r_{l} \over 2},\label{cent1j}
\ee
where we take $\ga $ to be some fixed small number. Hence under $\mathcal{B}_{m,l}$ we have $\sqrt{2\mathcal{T}_{u_{m},\wt r_{l}}^{ \tau}}\ge \al_{z,+}(l)$. The fact that
 for all $l\geq L-\left(4\log L\right)^{4}$
\begin{equation}
\Pbm\left[ \sqrt{2\mathcal{T}_{u_{m}, l-2}^{ \tau}}\le {1 \over 2}\al_{z,+}(l-2),\,\, \sqrt{2\mathcal{T}_{u_{m},\wt r_{l}}^{ \tau}}\ge \al_{z,+}(l) \right]\le ce^{-c'L^{2}}\label{notSmallm}
\end{equation}
then follows easily as in the proof of (\ref{bu2}). (In fact, the proof uses the same ideas but  is much easier).\qed

\section{Lower bounds for excursions} 
\label{sec-lowerboundex}

In this section we will prove the following.

\bl\label{lem-exclb}
There exist  $0<c_{1}, c_{2}<\ff$ such that for all 
  $L $ large and all $0\leq z\leq \log L$,
\begin{equation}
 \Pbm\left[\sup_{y\in F_{L}}\sqrt{2\mathcal{T}_{y,L}^{\tau}}\,\geq \rho_{L}L+z \right]\geq c_{1}{ (1+z) e^{ -2z} \over (1+z)   e^{ -2z}+c_{2} }.\label{goal.lbastv}
\end{equation} 
\el
This will immediately give the lower bounds in Theorems \ref{theo-tight} and \ref{theo-excovertight} and hence complete the proofs of those Theorems.

Note that for any $z_{0}$ it suffices to show that (\ref{goal.lbastv}) holds for all $z_{0}\leq z\leq \log L$, since by adjusting $c_{1}$ we then get (\ref{goal.lbastv}) for all $0\leq z\leq \log L$.

Let
\begin{equation}
 \beta_{z} \left(l\right)  =  \rho_{L}l+z,\label{eq:BetaDefn2}
\end{equation}
and 
\begin{equation}
\al_{z, -}\left(l\right)=\al_{z, -}\left(l,L,z\right)=\rho_{L}l+z- l_{L}^{ 1/4}.\label{14.1}
\end{equation}

  For each $k\ge1$ we define $\mathcal{T}_{y,l}^{k,m}$    be the number of excursions from  $\partial B_{d}\left(y,h_{l-1}\right)$ to $ \partial B_{d}\left(y,h_{l}\right)$  during the first $m$  excursions from  $\partial B_{d}\left(y,h_{k}\right)$ to $ \partial B_{d}\left(y,h_{k-1}\right)$. We abbreviate 
  $\mathcal{T}_{y, l}^{ 1}=\mathcal{T}_{y,l}^{1, x^{2}}$ with $x$ fixed.

Choose $r_{0}$ in (\ref{dr.1}) sufficiently small that $4h(r_{-1})\leq r^{\ast}$.
Let $ \wh r= h_{1}/20$, and with   $F^{0}:=B_{d}(v, \wh r) $ 
we set
\begin{equation}
F^{0}_{L}=F^{0}\cap F_{L},\hspace{.2 in}\mbox{so that }\hspace{.2 in}|F^{0}_{L}|\asymp e^{2L}. \label{F0.def}
\end{equation}
That is, $|F^{0}_{L}|$ does not depend on $r_{0}$.

In this section we show that
\bl\label{lem-exclb0}
There exists a $0<c<\ff$ such that for all $0<r_{0} $ sufficiently small,
  $L $ large and all $0\leq z\leq \log L$,
\begin{equation}
 \Pbm\left[\sup_{y\in F^{0}_{L}}\sqrt{2\mathcal{T}_{y,L}^{1}}\,\geq \rho_{L}L+z \right]\geq {(1+z)  e^{ -2z} \over  (1+z)  e^{ -2z}+c }.\label{goal.lbast}
\end{equation}
\el

Since the probability of $x^{2}$ excursions from $ \partial B_{d}(v, h_{1}-\wh r)$ to $ \partial B_{d}(v, h_{0}+\wh r)$ before $\tau$ is greater than $0$  and does not depend on $L$, (\ref{goal.lbast}) will imply  
Lemma \ref{lem-exclb}. We note that the $r_{0} $ used in this Lemma, and hence all $h_{l}$, are smaller than the corresponding quantities used until now. This is for notational convenience and, as can easily be seen, does not affect Lemma \ref{lem-exclb} which concerns large $L$. We could have kept the original $r_{0} $ and in place of $h_{l}$ used $h_{l+k}$ for some fixed $k$, but this would have made the notation cumbersome.

The proof of Proposition \ref{lem-exclb0} uses a modified 
second moment method and occupies the rest of this section.

We introduce the events $\II_{y,z}$, beginning with a barrier event. 
Let  
\begin{equation}
\wh\II_{y,z}=\left\{ \sqrt{2\mathcal{T}_{y,l}^{1}}\le\al_{z,-}\left(l\right)\mbox{ for }l=1,\ldots,L-1\mbox{ and }  \sqrt{2\mathcal{T}_{y,L}^{1}}\geq \rho_{L}+z\right\} ,\label{eq:TruncatedSummandLB}
\end{equation}
for $y\in  F_{L} $.  
As discussed in \cite{BRZ}, we need to augment $\wh\II_{y,z}$ by information on 
the angular increments of the excursions.
Instead of keeping track of individual 
excursions, we track the empirical measure of the increments, by comparing it
in Wasserstein distance to a reference measure. Recall that the Wasserstein
$L^1$-distance between probability measures on $\R$ is given by
\begin{equation}
  \dwa(\mu,\nu)=\inf _{\xi\in \PP^2(   \mu,\nu)}\Big \{   \int |x-y|\,d\xi(   x,y)\Big \},\label{was0}
\end{equation}
where $\PP^2(\mu,\nu)  $ denotes the set of probability measures on $\R\times\R$
with marginals $\mu,\nu$. If $\mu$
is a probability measure on $\R$ with finite support and if 
$\theta_{i},\, 1\leq i\leq n$ denotes a sequence of i.i.d 
$\mu$-distributed random variables 
then it follows from \cite[Theorem 2]{FG}  that for some $c_{0}=c_0(\mu)$
\begin{equation}
  \mbox{\rm Prob}\left\{   
  \dwa\(   \frac{1}{n} \sum_{i=1}^{ n} \de_{\theta_{i}},\mu  \)> \frac{c_{0}x}{\sqrt{n}}     \right\}\leq 2 e^{ -x^{ 2}}.\label{was1}
\end{equation}

Let $\BM_{t}$ be Brownian motion in the plane.
For each $k$ let  $\nu_{k}$ be the probability measure on $[0,2\pi]$ defined by 
 \be
\nu_{k}(   dx)=P^{   (r_{k}, 0)}\(\mbox{arg }\BM_{H_{\partial B(  0, r_{k-1})}}\in \,dx\),\label{meas45}
\ee  
where $\mbox{arg }x$ for $x\in \R^{2}$ is the argument of $x$ measured from the positive $x$-axis and $P^w$ is the law of $\BM_\cdot$ started from $w$.

Returning to $X_{t}$, our Brownian motion on the sphere, and using isothermal coordinates, see \cite[Section 2]{BRZ},
let   $0\leq \th_{k,i}\leq 2\pi$, $i=1,2,\ldots$ be  the angular increments centered at $y$,   mod $2\pi$,  from  $X_{H^{ i}_{\partial B(   y, h_{k})}}$ 
to $X_{H^{ i}_{\partial B(   y, h_{k-1})}}$,   the endpoints of the
$i$'th  excursion between $\partial B(   y, h_{k})$ and $\partial B(   y, h_{k-1})$. By the Markov property the $  \th_{k,i} $, $i=1,2,\ldots$ are independent,  and using  \cite[Section 2]{BRZ} we see that   each $  \th_{k,i} $ has distribution $\nu_{k}$. We set, for $n$ a positive integer,
\begin{equation}
\WW_{y,k}(n)=\left\{   \dwa
  \(   \frac{1}{n} \sum_{i=1}^{ n} \de_{\th_{k,i}},\nu_{k}  \) \leq  \frac{c_{0}\log (L-  k)}{2\sqrt{n}}     \right\}.\label{eq:LBTruncatedSumz}
\end{equation}

We are ready to define the good events $\II_{y,z}$. 
For $a\in \Z_{+}$ let
\begin{equation}
N_{k,a}=[(\rho_{L} k+z-a+1)^{ 2}/2].\label{nkdef.2}
\end{equation}
We set 
\begin{equation}
N_{k}=N_{k,a} \hspace{.1 in}\mbox{ if }\hspace{.1 in}\sqrt{2\mathcal{T}_{y,k}^{1}}\in I_{\rho_{L} k+z-a},\label{nkdef.2j}
\end{equation}
where $I_{s}=[s,s+1]$.
With $L_{+}= L -(500 \,log L )^{ 4}$ and $d^*$ a constant to be determined below,
let 
\begin{equation}
\II_{y,z}=\wh\II_{y,z}\cap_{k=L_{+}}^{L-d^{\ast}} \WW_{y,k}\(N_{k}\),
\label{eq:TruncatedSummandLBzz}
\end{equation}
and define the count
\begin{equation}
  J_z=\sum_{y\in F^{0}_{L}}{\bf 1}_{\II_{y,z}}.\label{eq:LBTruncatedSum}
\end{equation}
To obtain \eqref{goal.lbast}, we need a control on the
first and second moments of $J_z$, 
which is provided by the next two lemmas. In fact, \eqref{goal.lbast} will follow directly from these two Lemmas as in the proof of \cite[Proposition 4.2]{BRZ}, taking into account that $|F^{0}_{L}|$ does not depend on $r_{0}$.
Most of this section is devoted to their proof.
We emphasize that in the statements of the lemma, the implied constants are
uniform in $r_0$ smaller than a fixed small threshold.

\bl [First moment estimate]
\label{prop:LowerBoundOneProfilez} There is a large enough $d^{\ast}$, such that for all $L$ sufficiently large, all $0\leq   z\leq \log L$,
and all $y\in F^{0}_{L}$,
\begin{equation}
\Pbm  \left(\II_{y,z} \right)\asymp  (1+z) e^{-2L} e^{ -2z}.\label{eq:LowerBoundOneProfilez}
\end{equation}
\el

Let 
\begin{eqnarray}
  \label{eq-G0L} &&\\
  G_{0} & = & \left\{ \left(y,y'\right):y,y'\in F_{L}\mbox{ s.t. }d\left(y,y'\right)>2h_{0}\right\} ,\nn\\
G_{k} & = & \left\{ \left(y,y'\right):y,y'\in F_{L}\mbox{ s.t. }2h_{k}<d\left(y,y'\right)\le2h_{k-1}\right\} \mbox{ for }1\le k<L,\nn\\
G_{L} & = & \left\{ \left(y,y'\right):y,y'\in F_{L}\mbox{ s.t. }0<d\left(y,y'\right)\le2h_{L-1}\right\} \nn .
\end{eqnarray}
 
\bl     [Second moment estimate]
\label{prop:UpperBoundz}  There are large enough $d^*,c'$,
such that for all $L $ sufficiently  large,   
all $0\leq    z\leq  \log L $ 
and  all $(y, y')\in G_k $, $1\leq k\leq L$,
\begin{equation}
\Pbm\(\II_{y,z} \cap \II_{y',z}\)\leq  c'  (1+z)e^{-4L+2k}e^{ -2z  }e^{-ck_{L}^{1/4}}.\label{eq:UpperBoundz}
\end{equation}
\el

Before turning to the  proofs, we introduce some notation and
record some  simple estimates that will be
useful in calculations.
Recall  \eqref{eq:BetaDefn2}, (\ref{nkdef.2})-(\ref{nkdef.2j})
and for $a\in \Z_{+}$ let
\begin{equation}
  \HH_{k,a}=\lc  \sqrt{2\mathcal{T}_{y,k}^{1}}\in I_{\rho_{L} k+z-a} \rc=\lc   \sqrt{2\mathcal{T}_{y,k}^{1}}\in I_{\bb_{z}(k)-a}\rc.\label{hka.1}
\end{equation}
Note that on $\HH_{k,a} $ we have 
$N_{k}=N_{k,a}$.

Before proceeding we need to state a deviation inequality of 
Gaussian type for
the Galton-Watson process $T_l,l\ge0$ under $\Pgw_n$, the law of a critical Galton-Watson process with geometric offspring distribution with initial offspring $n$.
The proof is very similar to \cite[Lemma 4.6]{BK},
and is therefore omitted.
\begin{lemma}\label{lem: GW proc LD}
For all $n =1,2,3,\ldots$,
\begin{equation}
\Pgw_n\left(\left|\sqrt{2T_{l}}-\sqrt{2T_{0}}\right|\ge \theta
	\right)
	\le ce^{-\frac{\theta^{2}}{2l}},\quad \theta\geq 0.\label{eq: tail bound}
\end{equation}
\end{lemma}

In the proof of our moment estimates we will need the following. 
\bl\label{lem-moddev}
For any $  a, b\leq L/\log L $
\bea
&&\Pbm \left[\sqrt{2\mathcal{T}_{y,L}^{k,(   \bb_{z}(k)-a)^{ 2}/2}}\ge \rho_{L}L+z-b\right]\label{md.2}\\
&&\hspace{1 in}\leq  ce^{- 2(L-k)- 2(a-b)-{ ( a-b)^{2}\over 2(L-k)}}L^{2{(L-k)  \over L}}.\nn
\eea
\el
 
{\bf  Proof: } By (\ref{eq: tail bound}) we have  that for all  $\theta> n\ge1$  
\begin{equation}
\mathbb{P}\left[\mathcal{T}_{y,L}^{ k,n^{2}/2}\ge\theta^{2}/2\right]\le  c\exp\left(-{\left(\theta-n\right)^{2} \over 2(L-k)}\right).\label{11.44}
\end{equation}

We apply this with   $\theta= \rho_{L}L+z-b $ and 
\[n=  \bb_{z}(k)-a= \rho_{L}k+z -a \] so that 
\be 
\theta-n=\rho_{L}\(L-k \)   +a-b,
  \nn
\ee 
and hence 
\[
 {\(\theta-n\)^{2} \over 2(L-k)}\geq   2\(L-k \)-2{(L-k)  \over L}(\log L)
+2( a-b)+{ ( a-b)^{2}\over 2(L-k)}+O_{L}(1). 
\]
This gives (\ref{md.2}). 
 \qed

\subsection{First moment estimate}

In this subsection we prove Lemma \ref{prop:LowerBoundOneProfilez}.

 For the lower bound we have that
\begin{eqnarray}
&&\Pbm  \left[  \II_{y,z}\right]\geq P \left[\wh I_{y,z}\right]-\sum_{k=L_{+}}^{L-d^{\ast}}  \Pbm  \left[\wh \II_{y,z}\bigcap  \,W^{c}_{y,k}\(N_{k}\)  \right]
\label{ty.1a}\\
&&\geq c(1+z) e^{-2L}  e^{ -2z}-\sum_{k=L_{+}}^{L-d^{\ast}}  \Pbm  \left[\wh \II_{y,z}\bigcap  \,W^{c}_{y,k}\(N_{k}\)  \right], \nn
\end{eqnarray}
where  for $P \left[\wh I_{y,z}\right]$ we have used the barrier estimate (\ref{14.6a}) of Appendix I.
We note that
\begin{equation}
  \Pbm  \left[\wh \II_{y,z}\bigcap  \,W^{c}_{y,k}\(N_{k}\)  \right]\leq 
  \sum_{a\geq k_{L}^{ 1/4}}  \Pbm  \(\wh I_{y,z}^{k,a}\),\label{ty.1b}
\end{equation}
where
\begin{equation}
\wh I_{y,z}^{k,a}=\wh \II_{y,z}\bigcap \HH_{k,a}\,\bigcap \,W^{c}_{y,k}\(N_{k,a}\).\label{ty.1c}
\end{equation}
We show below that for all $L_{+}\leq k\leq L-d^{\ast}$ and $0\leq   z\leq \log L$,
\begin{equation}
  \sum_{a\geq k_{L}^{ 1/4}}  \Pbm  \(\wh I_{y,z}^{k,a}\)\leq c'(1+z) \(e^{-2L}  e^{ -2z}\)e^{ -c\log^{ 2} (L-  k)},\label{was1p.1}
\end{equation}
which will finish the proof of the lower bound for (\ref{eq:LowerBoundOneProfilez})  for $d^{\ast}$ sufficiently large.

Furthermore, it is easily seen using (\ref{md.2}) and the fact that $L-k\leq (500 \log L )^{ 4}$ that the sum in (\ref{was1p.1}) over $a\geq L^{3/4}$
is much smaller than the right hand side of (\ref{was1p.1}), hence it suffices to show that
\begin{equation}
  \sum_{a\geq k_{L}^{ 1/4}}^{L^{3/4}}   \Pbm  \(\wh I_{y,z}^{k,a}\)\leq c'(1+z) \(e^{-2L}  e^{ -2z}\)e^{ -c\log^{ 2} (L-  k)},\label{was1p}
\end{equation}

We now turn to the proof of (\ref{was1p}). Let
\be
J_{y,k}^{\uparrow} \label{eq:JDownz3}=\left\{ \sqrt{2\mathcal{T}_{y,l}^{ 1}}\le\rho_{L}l+z\mbox{ for }l=1,\ldots,k  \right\}, 
\ee
and
\be
\BB_{y,k,a} \label{eq:KUpz3}=\left\{ \sqrt{2\mathcal{T}_{y,L}^{ k,(\bb_{z}(k)-a)^{2}/2}}\ge \rho_{L}L+z\right\}. 
\ee
Then with 
\begin{equation}
\mathcal{K}_{k,p,a}=J_{y,k-3}^{\uparrow}\bigcap \HH_{k-3,p}\bigcap \HH_{k,a}\,\bigcap \,W^{c}_{y,k}\(N_{k,a}   \) \bigcap \BB_{y,k,a}\nn
\end{equation}
we have
\begin{equation}
\Pbm  \(\wh I_{y,z}^{k,a}\)\leq 
  \sum_{p\geq (k-3)_{L}^{ 1/4}}^{L^{3/4}}  \Pbm  \(\mathcal{K}_{k,p,a}\),\label{ty.1d}
\end{equation}
plus a term which is much smaller than  the right hand side of (\ref{was1p.1}).

Let 
\begin{equation}
  \WW^{\in x}_{y,k}(   n)=\left\{  \dwa
   \(\frac{1}{n} \sum_{i=1}^{ n} \de_{\th_{k,i}},\nu_{k}\) 
   \in \frac{c_{0}}{ \sqrt{n}} I_{ x  }     \right\},\label{eq:LBTruncatedSumz4}
\end{equation}
so that
\begin{equation}
\WW^{c}_{y,k}(N_{k,a})\subseteq 
\cup_{m=\log (L-k)}^{\ff} \WW^{\in m}_{y,k}( N_{k,a}),\label{compde.1}
\end{equation}
and consequently, setting
\begin{equation}
  \LL_{k,m,p,a}=\KK_{k,p,a}\cap \WW_{y,k}^{\in m}(N_{k,a}),  \label{extra.1}
\end{equation} 
we have
\begin{equation}
\Pbm \left(\KK_{k,p,a}\right)\leq 
\sum_{m=\log  (L-k)}^{\ff}   \Pbm \left(\LL_{k,m, p,a}\right).\label{wasd.1}
\end{equation}
 
Let
\begin{equation}
\LL'_{k,m, p,a}=:J_{y,k-3}^{\uparrow}\bigcap \HH_{k-3,p}\bigcap \HH_{k,a}\,\bigcap \,W^{\in m}_{y,k}\(N_{k,a}   \).\nn
\end{equation}
To prove (\ref{was1p.1}) it suffices to prove that for all $m\geq \log (L-k)$,
\begin{equation}
  \sum_{a\geq k_{L}^{ 1/4}}^{L^{3/4}}  \sum_{p\geq (k-3)_{L}^{ 1/4}}^{L^{3/4}} \Pbm  \(\BB_{y,k,a}\cap \LL'_{k,m, p,a}\)\leq c'(1+z) \(e^{-2L}  e^{ -2z}\)e^{ -cm^{ 2}  }. \label{was1p.1z}
\end{equation} 

\bl
\begin{eqnarray}
&&\Pbm \(\LL'_{k,m, p,a}\)= \Pbm \(J_{y,k-3}^{\uparrow}\bigcap \HH_{k-3,p}\bigcap \HH_{k,a}\,\bigcap \,W^{\in m}_{y,k}\(N_{k,a}   \) \)
\label{extra.6}\\
&&  \leq C(1+z)\(   1+p\)  e^{ -2k-2(z -p)}e^{-c(p-a)^{2}}e^{-m^{2}}. \nonumber
\end{eqnarray}
\el

{\bf  Proof: }
By (\ref{was1})
\begin{equation}
\Pbm  \(W^{\in m}_{y,k}\(N_{k,a}   \)\,|\,\HH_{k,a} \)\leq e^{-m^{2}}.\label{nty.2}
\end{equation}
By (\ref{md.2})
\begin{equation}
\Pbm  \(\HH_{k,a} \,|\,\HH_{k-3,p}\)\leq ce^{-c(p-a)^{2}},\label{nty.3}
\end{equation}
and by (\ref{genlow.1}) we see that
\begin{equation}
\Pbm  \(J_{y,k-3}^{\uparrow}\bigcap \HH_{k-3,p}\) \leq C(1+z)\(   1+p\)  e^{ -2k-2(z -p)}.\label{nty.4}
\end{equation}
\qed

It follows as in the proof of \cite[Lemma 4.7]{BRZ} that for some $M_{0}<\ff$
\begin{eqnarray}
&& \Pbm  \(\BB_{y,k,a}\cap \LL'_{k,m, p,a}\)\leq \Pbm \lc\sqrt{2\mathcal{T}_{y,L}^{k, (\bb_{z}(k)-a)^{2}/2 } }\ge \rho_{L}L+z-M_{0}m\rc
\label{was11.00}\\
&& \hspace{2 in} \times \Pbm \(\LL'_{k,m, p,a}\)+ e^{-4L}.\nonumber
\end{eqnarray}
We note that by (\ref{md.2})
\begin{equation}
\Pbm \lc\sqrt{2\mathcal{T}_{y,L}^{k, (\bb_{z}(k)-a)^{2}/2 } }\ge \rho_{L}L+z-M_{0}m\rc\leq  ce^{- 2(L-k)- 2(a-M_{0}m)-{ (a-M_{0}m)^{2}\over 2(L-k)}},\label{nty.1}
\end{equation}

Putting this all together with (\ref{extra.6}), and using $|a-p|\leq 1+(p-a)^{2}$ we find that
\begin{equation}
 \Pbm  \(\BB_{y,k,a}\cap \LL'_{k,m, p,a}\)\leq C(1+z) e^{ -2L-2z}e^{ -m^{ 2}/2}   \(   1+p\)e^{-c(p-a)^{2}}
e^{ -{ a^{2}\over 2(L-k)}}.     \label{nty.5}
\end{equation}
Summing first over $p$ and then over $a$ it is easy to see, using a fraction of the exponent $m^{ 2}/2$, that (\ref{was1p.1z}) holds for all $m\geq \log (L-k).$
This completes the proof of the lower bound in (\ref{eq:LowerBoundOneProfilez}) . 

Since $  \II_{y,z}\subseteq \wh \II_{y,z}$ the upper bound in (\ref{eq:LowerBoundOneProfilez}) follows from the barrier estimate (\ref{14.6}) of Appendix I. \qed

\subsection{Second moment estimate: branching in the bulk}

We prove the second moment estimate for $y, y'\in F^{0}_{L}$  with
 \[2h_{k-1}<d(   y,y')\leq 2h_{k-2}.\]

In this subsection we prove Lemma \ref{prop:UpperBoundz} for
\begin{equation}
(500 \log L )^{ 4}<k\le L-(500 \log L )^{ 4}.\label{eb.00}
\end{equation}

We  need to ``give ourselves a bit
of space'' , and we therefore define
\begin{equation}
k^{+}=k+\lceil100\log L\rceil.\label{eq:DefOfkPlus}
\end{equation}
Let 
\[\wh\II_{y,z; k\pm 3}=\left\{ \sqrt{2\mathcal{T}_{y,l}^{ 1}}\le\rho_{L}l+z;\,l=1,\ldots, k-4, k+4,\ldots,L-1\right\}\]
\be\hspace{3in}\cap \left\{ \sqrt{\mathcal{T}_{y,L}^{1}}\ge \rho_{L}L+z\right\} ,\label{eq:JDownside}
\ee 
where we have skipped the barrier condition for $k-3,\ldots, k+3$.
To obtain
the two point bound for the range (\ref{eb.00}) we will
bound the probability of
\begin{equation}
\wh\II_{y,z; k\pm 3}\cap \left\{\sqrt{2\mathcal{T}_{y',L}^{k^{+},\al_{z,-}^{2}(k^{+})/2}}\ge \rho_{L}L+z\right\},\label{eq:SmallKBiggerEvent+}
\end{equation}
which   contains the event $\II_{y,z}\cap \II_{y',z}$.  

Let $\mathcal{G}^{y'}$ denote the $\si$-algebra generated by the excursions from $\partial B_{d}(y', h_{k-1})$ to $\partial B_{d}(y',h(r_{k^{+}}))$. Note that $\wh\II_{y,z; k\pm 3}\in \mathcal{G}^{y'}$. 
Since \[\left\{\sqrt{2\mathcal{T}_{y',L}^{k^{+},\al_{z,-}^{2}(k^{+})/2}}\ge \rho_{L}L+z\right\}\] is measurable with respect to the first $\al_{z,-}^{2}(k^{+})$ excursions from $\partial B_{d}(y',h(r_{k^{+}}))$ to $\partial B_{d}(y',h(r_{k^{+}-1}))$,
it follows from the basic ideas in  \cite[sub-section 6.2]{BK} that
\begin{eqnarray}
&&\Pbm \left[\wh\II_{y,z; k\pm 3},\,\sqrt{2\mathcal{T}_{y',L}^{k^{+},\al_{z,-}^{2}(k^{+})/2}}\ge \rho_{L}L+z\right]
\label{125}\\
&&\leq c  \Pbm \left[\wh\II_{y,z; k\pm 3}\right]P \left[ \sqrt{2\mathcal{T}_{y',L}^{k^{+},\al_{z,-}^{2}(k^{+})/2}}\ge \rho_{L}L+z\right]. \nonumber
\end{eqnarray}

By Lemma \ref{prop:BarrierSecGWProp}
\begin{equation}
\Pbm \left[\wh\II_{y,z; k\pm 3}\right]\leq  c(1+z)e^{-2L}  e^{ -2z}.\label{Ipm.3}
\end{equation}

Using  (\ref{md.2}) for the last term in  (\ref{125}) together with  the fact that in the range (\ref{eb.00}) we have $(k^{+})_{L}^{ 1/4}\geq 500\log L$, we find that (\ref{125}) is bounded by
\bea
&& c(1+z)e^{-2L}  e^{ -2z}e^{- 2(L-k^{+})- 2(k^{+})_{L}^{ 1/4}}L^{2},\label{md.1+}\\
&&\leq c(1+z)e^{-2L}L^{ 20 2} e^{- 2(L-k)-   2(k^{+})_{L}^{ 1/4}} e^{ -2z}.
\nn\\
&&\leq c(1+z)e^{-2L}e^{- 2(L-k)-  k_{L}^{ 1/4}} e^{ -2z}.
\nn
\eea

\subsection{Second moment estimate: early branching}

In this subsection we prove Lemma \ref{prop:UpperBoundz} for \[1\le k<(500 \log L )^{ 4}.\]


Since we no longer have  $k_{L}^{ 1/4}\geq  \log L$ we will have to use barrier estimates to control the factors of $L$ such as arise in the first line of (\ref{md.1+}). On the other hand, since the number of excursions at lower levels is not so great we don't need such a large separation. Let
\begin{equation}
\wt k =k+\lceil100\log k\rceil, \hspace{.2 in}  k_{z} =k+\lceil100\log z\rceil.\label{eq:DefOfKPlus}
\end{equation}

  For $v\in\left\{ y,y'\right\} $
\bea
J_{v,s,\wt k}^{\downarrow}&=&\left\{ \sqrt{2\mathcal{T}_{v,l}^{\wt k,s^{2}/2}}\le\rho_{L}l+z \mbox{ for }l=\wt k+1,\ldots,L-1;\right.\nn\\
&&\left.\hspace{2 in}\,\, \sqrt{2\mathcal{T}_{v,L}^{ \wt k,s^{2}/2}}\ge \rho_{L}L+z \right\} ,\label{eq:JDown}\nn
\eea
with the barrier condition applied only for $l\ge \wt k$.

We first consider the case where $z\leq 100k$. Then
\be 
 \Pbm \(I_{y,z}\cap I_{y',z}\) 
  \leq  \sum_{n=1}^{\al_{z,-}\left(\wt k\right)}\Pbm \(J_{y,n,\wt k}^{\downarrow} \cap \wh\II_{y',z; k\pm 3}\). \label{newbreak.1}
\ee

Let $\mathcal{G}^{y}$ denote the $\si$-algebra generated by the excursions from $\partial B_{d}(y,h_{k-1})$ to $\partial B_{d}(y,h(r_{\wt k}))$. Note that $\wh\II_{y',z; k\pm 3} \in \mathcal{G}^{y}$. Since, under our assumption that $z\leq 100k$, the number of excursions from $\partial B_{d}(y,h_{k-1})$ to $\partial B_{d}(y,h(r_{\wt k}))$ is dominated by $n=O(k^{2})$, it follows as before that
\bea
&&
 \Pbm\(J_{y,n,\wt k}^{\downarrow}\cap \wh\II_{y',z; k\pm 3} \)\leq cP \(J_{y,n,\wt k}^{\downarrow} \)    \Pbm \( \wh\II_{y',z; k\pm 3} \).   \label{eb.30}
\eea

By the barrier estimate (\ref{14.8}) of Appendix I, with $n=\bb_{z}(\wt k)-t$
\begin{equation}
 \Pbm \(J_{y,n,\wt k}^{\downarrow} \) \leq  ct\,n^{1/2}e^{- 2(L-\wt k)- 2t   }
\leq  ck^{202} e^{- 2(L-  k)- 2 k_{L}^{ 1/4}}.\label{eb.30c}
\end{equation}
where the last step followed from the fact that  $k_{L}^{ 1/4}\leq\wt k_{L}^{ 1/4}\leq t\leq \bb_{z}(\wt k)\leq ck$. Since, under our assumption that $z\leq 100k$, the number of terms in  (\ref{newbreak.1}) is $\leq ck^{2}$, and using (\ref{Ipm.3}), we find that (\ref{newbreak.1}) is bounded by
\bea
&&ck^{204} e^{- 2(L-  k)- 2k_{L}^{ 1/4}} (1+z) e^{-2L }   e^{ -2z}\label{eb.35}\\
&\leq &c (1+z) e^{-4L+2k- k_{L}^{ 1/4}}   e^{ -2z}.\nn
\eea 
 
 Thus we can assume that
 \begin{equation}
z\geq 100 k.\label{break.6}
\end{equation}

We have
\begin{eqnarray}
&& \Pbm \(I_{y,z}\cap I_{y',z}\)
\label{newbreak.2}\\
&&=\sum_{n=1}^{\al_{z,-}\left(k_{z}\right)}   1_{\{n=\bb_{z}(k_{z})-t;\,t\geq z/2\}}     \Pbm \(    \Big    \{\sqrt{2 \mathcal{T}_{y,k_{z}}^{1}}=n \Big \}              \cap  I_{y,z}\cap I_{y',z}\)   \nonumber\\
&&+\sum_{n=1}^{\al_{z,-}\left(k_{z}\right)}   1_{\{n=\bb_{z}(k_{z})-t;\,t< z/2\}}     \Pbm \(     \Big   \{\sqrt{ 2\mathcal{T}_{y,k_{z}}^{1}}=n \Big \}              \cap  I_{y,z}\cap I_{y',z}\)  \nonumber
\end{eqnarray}

Since in the above sums $n\leq cz $ in view of (\ref{break.6}), we can bound the first sum in (\ref{newbreak.2}) by
\bea
&&
\sum_{n=1}^{\al_{z,-}\left(k_{z}\right)}   1_{\{n=\bb_{z}(k_{z})-t;\,t\geq z/2\}}  \Pbm \(J_{y,n,k_{z}}^{\downarrow} \cap \wh\II_{y',z; k\pm 3}\)\label{newbreak.3}\\
&&\leq c \sum_{n=1}^{\al_{z,-}\left(k_{z}\right)}   1_{\{n=\bb_{z}(k_{z})-t;\,t\geq z/2\}}    \Pbm \(J_{y,n,k_{z}}^{\downarrow} \)\Pbm \( \wh\II_{y',z; k\pm 3}\)        \nn\\
&&\leq c \sum_{n=1}^{\al_{z,-}\left(k_{z}\right)}   1_{\{n=\bb_{z}(k_{z})-t;\,t\geq z/2\}}    \Pbm \(J_{y,n,k_{z}}^{\downarrow} \)(1+z) e^{-2L }   e^{ -2z},        \nn
\eea
as before. Instead of (\ref{eb.30c}) we now have
\begin{equation}
 \Pbm \(J_{y,n,k_{z}}^{\downarrow} \) \leq  ct\,z^{1/2}e^{- 2(L-k_{z})- 2t   }
\leq  cz^{202} e^{- 2(L-  k)- z/2},\label{eb.30cn}
\end{equation}
where the last inequality used $t\geq z/2$.
In view of (\ref{break.6}) and the fact that the number of terms in the sum is $\leq cz$, this gives the desired bound for the first sum in (\ref{newbreak.2}).

Note next  that if $t< z/2$ then we must have $n=\bb_{z}(k_{z})-t\geq z/2$, (but we still have $n\leq cz $ by (\ref{break.6})). Thus we can bound the second sum in (\ref{newbreak.2}) by
\bea
&&
\sum_{n,n'=1}^{\al_{z,-}\left(k_{z}\right)}   1_{\{n\geq z/2\}}  \Pbm \(    \Big \{ \sqrt{2 \mathcal{T}_{y,k_{z}}^{1}}=n\Big\}              \cap J_{y,n,k_{z}}^{\downarrow} \cap J_{y',n',k_{z}}^{\downarrow}\)\label{newbreak.4}\\
&&\leq c \sum_{n,n'=1}^{\al_{z,-}\left(k_{z}\right)}  1_{\{n\geq z/2\}}    \Pbm \(      \Big \{ \sqrt{2 \mathcal{T}_{y,k_{z}}^{1}}=n \Big \}              \cap J_{y,n,k_{z}}^{\downarrow} \)\Pbm \( J_{y',n',k_{z}}^{\downarrow}\).      \nn
\eea
as before. Then by the Markov property, this is bounded by
\be
 c  \sum_{n,n'=1}^{\al_{z,-}\left(k_{z}\right)}  1_{\{n\geq z/2\}}    \Pbm \(     \Big  \{ \sqrt{ 2\mathcal{T}_{y,k_{z}}^{1}}=n \Big \}   \)  \Pbm \(   J_{y,n,k_{z}}^{\downarrow} \) \Pbm \( J_{y',n',k_{z}}^{\downarrow}\).
\label{newbreak.5}
\ee
By (\ref{eq: tail bound}) with $n\geq z/2$ and then (\ref{break.6})
\begin{equation}
\Pbm \(     \Big  \{ \sqrt{ 2\mathcal{T}_{y,k_{z}}^{1}}=n \Big \}   \) \leq e^{-z^{2}/4k_{z}}\leq e^{-10z  },\label{newbreak.6}
\end{equation}
while now, instead of (\ref{eb.30cn}), we use
\begin{equation}
 \Pbm \(J_{y,n,k_{z}}^{\downarrow} \) \leq  ct\,z^{1/2}e^{- 2(L-k_{z})- 2t   }
\leq  cz^{202} e^{- 2(L-  k)},\label{eb.30cnp}
\end{equation}
and a similar bound for $\Pbm \( J_{y',n',k_{z}}^{\downarrow}\).$
Thus (\ref{newbreak.5}) is bounded by
\begin{equation}
c  \sum_{n,n'=1}^{\al_{z,-}\left(k_{z}\right)} e^{-10z  }
 z^{404} e^{- 4(L-  k)}\leq ce^{-10z  }
 z^{408} e^{ 2 k}e^{- 4L+2 k}.    \nn
\end{equation}
In view of (\ref{break.6}), this gives the desired bound for the second sum in (\ref{newbreak.2}).

\subsection{Second moment estimate: late branching} In this subsection we prove Lemma \ref{prop:UpperBoundz} for $L- (500 \log L )^{ 4}\le k<L-1$.

Consider first the case \[L- (500 \log L )^{ 4}\le k<L-d^{\ast}.\]
 We will bound
the probability of
\begin{equation}
\mathcal{A}= \left\{\sqrt{2\mathcal{T}_{y,L}^{k,\al_{z,-}^{2}(k)/2}}\ge \rho_{L}L+z\right\}\cap \WW_{y,k}\(N_{k}\) \cap \wh\II_{y',z; k\pm 3},\label{eq:LargeKBiggerEvent}
\end{equation}
(which  contains the event $\II_{y,z}\cap \II_{y',z}$).

It follows as in the proof of \cite[Lemma 4.7]{BRZ} that for some $M_{0}<\ff$
\begin{eqnarray}
&&\Pbm\(\mathcal{A}\)\leq \Pbm \lc\sqrt{2\mathcal{T}_{y,L}^{k, \al _{z,-}^{2}(k)/2 } }\ge \rho_{L}L+z-M_{0}\log (L-k)\rc
\label{was11.0}\\
&& \hspace{2 in} \times \Pbm \(\wh\II_{y',z; k\pm 3}\)+ e^{-4L}.\nonumber
\end{eqnarray}

Using (\ref{md.2}) and (\ref{Ipm.3}) this shows that
\be
 \Pbm\(\mathcal{A}\)\leq  ce^{- 2(L-k)+2M_{0}\log (L-k)-2k_{L}^{ 1/4} }(1+z)e^{-2L}  e^{ -2z}+ e^{-4L}
\label{was11.0a}
\ee
By taking $d^{\ast}$ sufficiently large we will have $M_{0}\log (L-k)\leq k_{L}^{ 1/4} /2$, which then gives (\ref{eq:LowerBoundOneProfilez}).\qed

For $L-d^{\ast}\le k<L-1$ we simply bound the term $\Pbm\left(\II_{y,z}\cap \II_{y',z}\right) $ by $\Pbm\left(\II_{y,z}\right)$
and obtain from (\ref{eq:LowerBoundOneProfilez}) the following upper
bound
\begin{equation}
\Pbm\left(\II_{y,z}\cap \II_{y',z}\right) \leq c(1+z)e^{-2L}e^{ -2z}\leq c_{d^{\ast}}(1+z) e^{-(4L-2k)-ck^{1/4}_{L}}e^{ -2z}.\label{eq:LowK}
\end{equation}
\qed

\section{Excursion counts and occupation measure on $\S^{ 2}$.}\label{sec-occmeas}

In this section we prove Theorems \ref{theo-tight0} and \ref{theo-spheretight}.

For $0<\ep<a<b<\pi$, let  
$ \MM_{x,\ep, a,b}(n)$ be the total     occupation measure of $B_{d}(   x,\ep)$ until the end of the first $n$ excursions   from 
$\partial B_{d}\left(x,a\right)$ to $\partial B_{d}\left(x,b\right)$. With $\om_{\ep}=2\pi (1-\cos (\ep))$, the area of $ B_{d}\left(x,\ep\right)$, let 
\begin{equation}
\overline \MM_{x, \ep, a,b}(n)={1 \over \om_{\ep}}\MM_{x,\ep, a,b}(n).\label{occm.n}
\end{equation}
 In particular, when starting from   $\partial B_{d}(   x,a)$,
\begin{equation}
\overline \MM_{x,\ep, a,b}(1)={1 \over \om_{\ep}}\int_{0}^{ H_{ \partial B_{d}(   x,b)}} 1_{\{B_{d}(   x,\ep)\}}(  X_{t} )\,dt.\label{occm.1}
\end{equation}
The following Lemma is proven in Section \ref{sec-dev}.

\bl\label{lem-clt}
For some $c >0$, uniformly in $x\in  S^{2}$, and $h_{k}/100\leq \ep\leq h_{k}$,
\begin{equation}
\Pbm\(\overline \MM_{x,\ep, h_{k}, h_{k-1}}(n)\leq {1 \over \pi}(1-\de) \,\,n\)\leq e^{-c\de^{2}n} \label{clt.1}
\end{equation}
and
\begin{equation}
\Pbm\(\overline \MM_{x,\ep, h_{k}, h_{k-1}}(n)\geq {1 \over \pi}(1+\de) \,\,n\)\leq e^{-c\de^{2}n} \label{clt.1a}
\end{equation}
\el

Recall $ \bar \mu_{\tau}\(y,\ep\)$ from (\ref{tp.1}) and set
\begin{equation}
 t_{L}\left(z\right)=2L\left(L-\log L+ z\right).\label{clt.2}
\end{equation}

\bl\label{lem-time}
We can find $0<c,c',z_{0}<\ff$ such that for $L $ large,  all $z_{0}\leq z\leq \log L$, and all 
$h_{L}/100\leq \ep_{y}\leq h_{L}$
\begin{equation}
c' ze^{-2z}\leq \Pbm\( \exists y\in F_{L} \mbox{ s.t.  } \bar \mu_{\tau}\(y,\ep_{y}\)\geq   {1 \over \pi}  t_{L}\left(z \right)\)\leq 
cze^{-2z}.\label{time.1}
\end{equation}
\el

For the sphere, it suffices to take $\ep_{y}=\ep$  independent of $y$. The present formulation is needed for the plane, as we will see in Section \ref{sec-euclid}.
To clarify the connection with (\ref{goal.1asto}) we note that for some $0<c_{\ast}=c_{\ast}(r_{0})<\ff$,
\begin{equation}
\( m_{h_{L}}+z\)^{2}= {1 \over \pi}  t_{L}\left(\sqrt{2\pi}z +c_{\ast}+o_{L}(1)\right). \label{con1.9}
\end{equation}

\subsection{The upper bound}

We first show that, with $F^{+}_{L}=F_{L}\cap B_{d} \left(v,h_{\log L} \right)$,
\begin{equation}
\PP_{1}=: \Pbm\( \exists y\in F^{+}_{L} \mbox{ s.t.  } \bar \mu_{\tau}\(y,\ep_{y}\)\geq   {1 \over \pi}  t_{L}\left(z \right)\)\leq 
cze^{-2z}.\label{timem.2}
\end{equation} 
If 
\begin{equation}
\wh \AA_{L,z}= \Bigg\{ \sup_{y\in F^{+}_{L} }\sqrt{2\mathcal{T}_{y,L}^{\tau}}\,\geq \rho_{L}L+z \Bigg\}, \nn
\end{equation}
then by (\ref{eq:NotHitByrLlogL})
\begin{eqnarray}
 \PP_{1}&\leq &\Pbm\(\wh \AA_{L,z}\) 
 + \Pbm\(\wh \AA^{c}_{L,z}, \exists y\in F^{+}_{L}  \mbox{ s.t.  } \bar \mu_{\tau}\(y,\ep_{y}\)\geq   {1 \over \pi}  t_{L}\left(z \right)\) \nonumber\\
& \leq& c e^{ -2z}+ \Pbm\(\wh \AA^{c}_{L,z}, \exists y\in F^{+}_{L} \mbox{ s.t.  } \bar \mu_{\tau}\(y,\ep_{y}\)\geq   {1 \over \pi}  t_{L}\left(z \right)\). \nonumber
\end{eqnarray}
Recalling   the notation $F^{m}_{L}=F_{L}\cap B_{d}^{c}\left(v,h_{m} \right)\cap B_{d}\left(v,h_{m-1} \right)$, we then bound   
\begin{eqnarray}
&&\Pbm\(\wh \AA^{c}_{L,z}, \exists y\in F^{+}_{L}  \mbox{ s.t.  } \bar \mu_{\tau}\(y,\ep_{y}\)\geq   {1 \over \pi}  t_{L}\left(z \right)\)
\label{eq:NotHitByrLlogLag}\\
&&\leq   \sum_{m=\log L}^{L} ce^{2(L-m)}\nn\\
&&\hspace{.5 in}\sup_{y\in F^{m}_{L} } \Pbm\(\sqrt{2\mathcal{T}_{y,L}^{\tau}}\,\leq \rho_{L}L+z,\,\bar \mu_{\tau}\(y,\ep_{y}\)\geq   {1 \over \pi}  t_{L}\left(z \right)\).\nonumber
\end{eqnarray}
We can  write
 \begin{eqnarray}
&&\Pbm\( \sqrt{2\mathcal{T}_{y,L}^{\tau}}\,\leq \rho_{L}L+z,\,\bar \mu_{\tau}\(y,\ep_{y}\)\geq   {1 \over \pi}  t_{L}\left(z \right)\)
\label{36.13n}\\
&& =\sum_{j=1}^{z+ML^{1/2}}P\( \,\sqrt{2\mathcal{T}_{y,L}^{\tau }}\in I_{\rho_{L}L+z -j} \mbox{ and } \bar \mu_{\tau}\(y,\ep_{y}\)\geq   {1 \over \pi}  t_{L}\left(z \right)\)  \nonumber\\
&& + \Pbm\(\,\sqrt{2\mathcal{T}_{y,L}^{\tau }}\le \rho_{L}L-ML^{1/2} \mbox{ and } \bar \mu_{\tau}\(y,\ep_{y}\)\geq   {1 \over \pi}  t_{L}\left(z \right)\).   \nonumber
\end{eqnarray}

\bl\label{lem-Mbounds}
For all $y\in F^{m}_{L}$,  $\log L \leq m\leq L$         and $j\leq z+M L^{1/2}$
\bea
&&P\( \,\sqrt{2\mathcal{T}_{y,L}^{\tau }}\in I_{\rho_{L}L+z -j} \mbox{ and } \bar \mu_{\tau}\(y,\ep_{y}\)\geq   {1 \over \pi}  t_{L}\left(z \right)\)\label{Mbound.1}\\
&&\leq cme^{-2L}L e^{ -2(z-j)}e^{-c'j^{2}},\nn
\eea
and
\bea
&&
\Pbm\(\,\sqrt{2\mathcal{T}_{y,L}^{\tau }}\le \rho_{L}L-ML^{1/2} \mbox{ and } \bar \mu_{\tau}\(y,\ep_{y}\)\geq   {1 \over \pi}  t_{L}\left(z \right)\)\label{Mbound.2}\\
&&\leq cme^{-2L}L e^{ -2(z-j)}e^{-4L}.\nn
\eea
\el

{\bf  Proof of Lemma     \ref{lem-Mbounds}:}
By (\ref{eq:TsDivdedByL})
\be 
  \(  \rho_{L}L+z-j\)^{ 2}/2\leq t_{L}\left(z-j +2M^{ 2}\right)
\label{calc.w1} 
\ee
for all  $j\leq z+M L^{1/2}$. Hence for such $j$
\begin{eqnarray}
 &&\Pbm\(\sqrt{2\mathcal{T}_{y,L}^{\tau }}\in I_{\rho_{L}L+z -j} \mbox{ and } \bar \mu_{\tau}\(y,\ep_{y}\)\geq   {1 \over \pi}  t_{L}\left(z \right)\)
\nn\\
&&  \leq \Pbm\( \sqrt{2\mathcal{T}_{y,L}^{\tau }}\in I_{\rho_{L}L+z -j} ,\,\,  \overline \MM_{y, \ep_{y},h_{L},h_{L-1}}\(t_{L}\left(z-j+2M^{ 2} \right) \)\geq   {1 \over \pi}  t_{L}\left(z \right)\). \nonumber
\end{eqnarray}
Using the Markov property and then (\ref{eq:NotHitByrL}), we have for $ y\in F^{m}_{L}$ this is
\begin{eqnarray}
&& =   \Pbm\( \sqrt{2\mathcal{T}_{y,L}^{\tau }}\in I_{\rho_{L}L+z -j} \)\Pbm\( \overline \MM_{y, \ep_{y},h_{L},h_{L-1}}\(t_{L}\left(z-j+2M^{ 2} \right) \)\geq   {1 \over \pi}  t_{L}\left(z \right)   \)\nonumber \\
&&\leq cme^{-2L}L e^{ -2(z-j)}\Pbm\( \overline \MM_{y, \ep_{y},h_{L},h_{L-1}}\(t_{L}\left(z-j+2M^{ 2} \right) \)\geq   {1 \over \pi}  t_{L}\left(z \right)   \).\nonumber
\end{eqnarray}

Consider first the case of $4M^{ 2} \leq j$.
We now apply (\ref{clt.1a}) with \[n=t_{L}\left(z-j+2M^{ 2} \right)=t_{L}\left(z \right)-2(j-2M^{ 2})L \sim L^{2}\] and   \[\de=2(j-2M^{ 2})L/t_{L}\left(z-j+2M^{ 2} \right)<<1 \] for $4M^{ 2} \leq j\leq z+M L^{1/2}$ to see that 
\bea
&&\hspace{1 in}
\Pbm\( \overline \MM_{y, \ep_{y},h_{L},h_{L-1}}\(t_{L}\left(z-j+2M^{ 2} \right) \)\geq   {1 \over \pi}  t_{L}\left(z \right)  \)  \label{ez.3}\\
&&=\Pbm \( \overline \MM_{y, \ep_{y},h_{L},h_{L-1}}\(t_{L}\left(z-j+2M^{ 2} \right) \)\right.\nn\\
&&\left.\hspace{1.5 in}\geq   {1 \over \pi} \( t_{L}\left(z-j+2M^{ 2}  \right)+2(j-2M^{ 2})L\)  \)\nn\\
&&=\Pbm\( \overline \MM_{y, \ep_{y},h_{L},h_{L-1}}\(t_{L}\left(z-j+2M^{ 2} \right) \)\right.\nn\\
&&\left.\hspace{1 in}\geq   {1 \over \pi} \( 1+{2(j-2M^{ 2})L \over t_{L}\left(z-j+2M^{ 2}  \right)}\)t_{L}\left(z-j+2M^{ 2}  \right)  \)\nn\\
&&\leq e^{-c{\((j-2M^{ 2})L\)^{2} \over t_{L}\left(z-j +2M^{ 2} \right)}}\leq e^{-c'j^{2}}.\nn
\eea

Similarly, for (\ref{Mbound.2}) we use
\begin{eqnarray}
&&\Pbm\(\sqrt{2\mathcal{T}_{y,L}^{\tau }}\le \rho_{L}L-ML^{1/2}  \mbox{ and } \bar \mu_{\tau}\(y,\ep_{y}\)\geq   {1 \over \pi}  t_{L}\left(z \right)\)
\label{36.12e}\\
&&  \leq \Pbm\( \overline \MM_{y, \ep_{y},h_{L},h_{L-1}}\(t_{L}\left(-ML^{1/2}+2M^{ 2} \right) \)\geq   {1 \over \pi}  t_{L}\left(z \right)\)\leq  e^{-4L}\nonumber
\end{eqnarray} 
by (\ref{ez.3}) for $M$ sufficiently large.

For $j<4M^{ 2} $ we simply bound the probability in (\ref{ez.3}) by $1$ which we can bound by 
$Ce^{-c'j^{2}}$ for $C$ sufficiently large.
 \qed

Then using  (\ref{eq:NotHitByrLlogLag}) and Lemma     \ref{lem-Mbounds} we see that
\begin{eqnarray}
&&\hspace{.2 in}\Pbm\(\wh\AA^{c}_{L,z}, \exists y\in F^{+}_{L} \mbox{ s.t.  } \bar \mu_{\tau}\(y,\ep_{y}\)\geq   {1 \over \pi}  t_{L}\left(z \right)\)
\label{36.12f}\\
&& \leq  C\sum_{m=\log L}^{  L} cmLe^{2(L-m)}\sum_{j=1}^{z+ML^{1/2}}e^{-2L} e^{ -2(z-j) }   e^{-c'j^{2}}\nn\\
&&+  \sum_{m=\log L}^{  L} cmLe^{2(L-m)}e^{-4L}.\nonumber
\end{eqnarray}
This is easily seen to be bounded by the right hand side of (\ref{timem.2}).

Recalling the notation $ F^{\ast}_{L}=F_{L}\cap B_{d}^{c}\left(v,h_{\log L} \right)$
from (\ref{label.2}), to complete the proof of the upper bound it remains to show that
\begin{equation}
\PP_{2} =:\Pbm\( \exists y\in F^{\ast}_{L} \mbox{ s.t.  } \bar \mu_{\tau}\(y,\ep_{y}\)\geq   {1 \over \pi}  t_{L}\left(z \right)\)\leq 
cze^{-2z}.\label{timem.1}
\end{equation}

Note that with $k_{y}$ as in (\ref{label.1}), if
\begin{equation} 
\AA_{L,z}= \Big\{\exists y\in F^{\ast}_{L}, l\in\left\{ k_{y}+1,\ldots, L\right\}\mbox{ s.t. } \mathcal{T}_{y,l}^{\tau }\ge\alpha_{z,+}^{2}\left(l\right)/2 \Big\},\label{cru.1}
\end{equation}
then
\begin{eqnarray}
 \PP_{2} &\leq & \Pbm\(\AA_{L,z}\)
+ \Pbm\(\AA^{c}_{L,z}, \exists y\in F^{\ast}_{L} \mbox{ s.t.  } \bar \mu_{\tau}\(y,\ep_{y}\)\geq   {1 \over \pi}  t_{L}\left(z \right)\)  \nonumber\label{cru.2}\\
 &\leq & cze^{ -c^{ \ast}z }+ \Pbm\(\AA^{c}_{L,z}, \exists y\in F^{\ast}_{L} \mbox{ s.t.  } \bar \mu_{\tau}\(y,\ep_{y}\)\geq   {1 \over \pi}  t_{L}\left(z \right)\),   \nonumber
\end{eqnarray}
by (\ref{eq:SmartMarkov20}). Recalling again the notation $F^{m}_{L}=F_{L}\cap B_{d}^{c}\left(v,h_{m} \right)\cap B_{d}\left(v,h_{m-1} \right)$, we have that 
\begin{eqnarray}
&&\Pbm\(\AA^{c}_{L,z}, \exists y\in F^{\ast}_{L} \mbox{ s.t.  } \bar \mu_{\tau}\(y,\ep_{y}\)\geq   {1 \over \pi}  t_{L}\left(z \right)\)
\label{36.11a}\\
&& =\sum_{m=1}^{\log L}\Pbm\(\AA^{c}_{L,z}, \exists y\in F^{m}_{L} \mbox{ s.t.  } \bar \mu_{\tau}\(y,\ep_{y}\)\geq   {1 \over \pi}  t_{L}\left(z \right)\)  \nonumber
\end{eqnarray}

Since 
\begin{equation}
\AA^{c}_{L,z} =  \Big\{ \sup_{y\in F^{\ast}_{L}}\mathcal{T}_{y,l}^{\tau }\le\alpha_{z,+}^{2}\left(l\right)/2, k_{y}+1\leq l\leq  L \Big\}\label{36.11}
\end{equation}
and    $k_{y}=m$ for $y\in F^{m}_{L}$,   we see that
\begin{eqnarray}
&&\hspace{.2 in}\Pbm\(\AA^{c}_{L,z}, \exists y\in F^{m}_{L} \mbox{ s.t.  } \bar \mu_{\tau}\(y,\ep_{y}\)\geq   {1 \over \pi}  t_{L}\left(z \right)\)
\label{36.12}\\
&& \leq   ce^{2(L-m)}\nn\\
&&\sup_{y\in F^{m}_{L} }\Pbm\(\mathcal{T}_{y,l}^{\tau }\le\alpha_{z,+}^{2}\left(l\right)/2, m+1\leq l\leq  L  \mbox{ and } \bar \mu_{\tau}\(y,h_{L}\)\geq   {1 \over \pi}  t_{L}\left(z \right)\).\nonumber
\end{eqnarray}
With 
\begin{equation}
\BB^{ y}_{L,m, z}=   \Big\{ \mathcal{T}_{y,l}^{\tau }\le\alpha^2_{z,+}\left(l\right)/2, m+1\leq l\leq  L-1 \Big\}\label{cru.3}
\end{equation}
we have for $y\in F^{m}_{L}$,
\begin{eqnarray}
&&\hspace{.3in}\Pbm\(\mathcal{T}_{y,l}^{\tau }\le\alpha_{z,+}^{2}\left(l\right)/2, m+1\leq l\leq  L  \mbox{ and } \bar \mu_{\tau}\(y,\ep_{y}\)\geq   {1 \over \pi}  t_{L}\left(z \right)\)
\label{36.13}\\
&& =\sum_{j=1}^{z+ML^{1/2}}P\(\BB^{ y}_{L,m,z},\,\sqrt{2\mathcal{T}_{y,L}^{\tau }}\in I_{\alpha_{z,+}\left(L\right) -j} \mbox{ and } \bar \mu_{\tau}\(y,\ep_{y}\)\geq   {1 \over \pi}  t_{L}\left(z \right)\)  \nonumber\\
&& + \Pbm\(\BB^{ y}_{L,m,z},\,\sqrt{2\mathcal{T}_{y,L}^{\tau }}\le \alpha_{z,+}\left(L\right)-z-ML^{1/2} \mbox{ and } \bar \mu_{\tau}\(y,\ep_{y}\)\geq   {1 \over \pi}  t_{L}\left(z \right)\).  \nonumber
\end{eqnarray}
Here, $M\geq 1$ is a fixed constant to be chosen shortly.

Recalling, see (\ref{eq:AlphaBarrierDef}), that $ \alpha_{z,+}\left(L\right)=\rho_{L}L+z$, and using (\ref{calc.w1}) we see that 
\begin{eqnarray}
&&\(   \alpha_{z,+}\left(L\right) -j\)^{ 2}/2 \leq t_{L}\left(z-j +2M^{ 2}\right)\nn
\end{eqnarray}
for all  $j\leq z+M L^{1/2}$.
It follows that for such $j$
\begin{eqnarray}
&&P\(\BB^{ y}_{L,m,z},\,\sqrt{2\mathcal{T}_{y,L}^{\tau }}\in I_{\alpha_{z,+}\left(L\right) -j} \mbox{ and } \bar \mu_{\tau}\(y,\ep_{y}\)\geq   {1 \over \pi}  t_{L}\left(z \right)\)
\nn\\
&&  \leq \Pbm\(\BB^{ y}_{L,m,z},\,\sqrt{2\mathcal{T}_{y,L}^{\tau }}\in I_{\alpha_{z,+}\left(L\right) -j},  \overline \MM_{y,\ep_{y}, h_{L},h_{L-1}}\(t_{L}\left(z-j+2M^{ 2} \right) \)\geq   {1 \over \pi}  t_{L}\left(z \right)\) \nonumber\\
&& =   \Pbm\(\sqrt{2\mathcal{T}_{y,l}^{\tau }}\le\alpha_{z,+}\left(l\right), m+1\leq l\leq  L-1,\,\sqrt{2\mathcal{T}_{y,L}^{\tau }}\in I_{\alpha_{z,+}\left(L\right) -j}\)\nn\\
&&\hspace{1.8 in}   \times \Pbm\( \overline \MM_{y,\ep_{y}, h_{L},h_{L-1}}\(t_{L}\left(z-j +2M^{ 2}\right) \)\geq   {1 \over \pi}  t_{L}\left(z \right)   \)\nonumber,\label{cru.5}
\end{eqnarray}
by the Markov property. Using the barrier estimate (\ref{18.20}) of Appendix I, and recalling that $m=k_{y}<\log L$, this is bounded by 
\be 
 ce^{-2L} e^{ -2(z-j) } \times m^{2}j\left(z+m  \right)    \Pbm\( \overline \MM_{y, \ep_{y},h_{L},h_{L-1}}\(t_{L}\left(z-j +2M^{ 2}\right) \)\geq   {1 \over \pi}  t_{L}\left(z \right)   \).\label{36.14}
\ee

The rest of the proof of  (\ref{timem.1}) follows as in  the proof of (\ref{timem.2}). This completes the proof of the upper bound in Lemma \ref{lem-time}.

We now remove the restriction that $y\in F_{L}$ in the upper bound, subject to a continuity restriction on $\ep_{y}$.
\bl\label{lem-unif}
We can find $0<c, C, z_{0}<\ff$ such that for $L $ large,  all $z_{0}\leq z\leq \log L$, and all 
$h_{L}/20\leq \ep_{y}\leq h_{L+1} $ such that $|\ep_{y}-\ep_{y'}|\leq C\,d(y,y')/L$ for all $y,y'\in \S^{ 2}$,
\begin{equation}
 \Pbm\( \exists y \mbox{ s.t.  } \bar \mu_{\tau}\(y,\ep_{y}\)\geq   {1 \over \pi}  t_{L}\left(z \right)\)\leq 
cze^{-2z}.\label{unif.1}
\end{equation}
\el
 {\bf   Proof of Lemma \ref{lem-unif}: } 
 Let $F'_{L}$ be the centers of a  $\frac{d_{0}}{ L}h_{L}$ covering of $\S^{ 2}$ which contains $F_{L}$. For any $y\in \S^{ 2}$ we can find $y'\in F'_{L}$ such that $d(y,y')\leq \frac{d_{0}}{ L}h_{L}$, so that by our assumptions
 $|\ep_{y}-\ep_{y'}|\leq C\,\frac{d_{0}}{ L^{2}}h_{L}$. If we set $\bar\ep_{y}=\(1+\frac{1}{ L}\)\ep_{y}$ for all  $y\in \S^{ 2}$  we see that for $L$ large $h_{L}/30\leq \bar \ep_{y}\leq 2h_{L+1} $ and $|\bar \ep_{y}-\bar \ep_{y'}|\leq  \,\frac{d_{0}}{ L }h_{L}$.
  It follows from Lemma \ref{lem-interp} below that it suffices to prove that 
 \begin{equation}
 \Pbm\( \exists y\in F'_{L}  \mbox{ s.t.  } \bar \mu_{\tau}\(y,\bar \ep_{y}\)\geq   {1 \over \pi}  t_{L}\left(z \right)\)\leq 
cze^{-2z}.\label{unif.1L}
\end{equation}

 For $0<\ep<a<b<\pi$, let  
$ \MM_{y,\bar \ep_{y}, y_{0}, a,b}(n)$ be the total     occupation measure of $B_{d}(   y,\bar \ep_{y})$ during the first $n$ excursions   from 
$\partial B_{d}\left(y_{0},a\right)$ to $\partial B_{d}\left(y_{0},b\right)$. With $\om_{\ep}=2\pi (1-\cos (\ep))$, the area of $ B_{d}\left(y,\ep\right)$, let 
\begin{equation}
\overline \MM_{y,\bar \ep_{y}, y_{0}, a,b}(n)={1 \over \om_{\bar \ep_{y}}}\MM_{y,\bar \ep_{y}, y_{0}, a,b}(n).\label{unif.1a}
\end{equation}
For $y_{0}\in F_{L}$ let
 \begin{equation}
D_{y_{0}}=\{   y\in F'_{L}\,|\, d\(   y,y_{0}\)\leq d_{0}h_{L}/2\}.\label{unif.2}
 \end{equation}
Following the proof of the upper bound for Lemma \ref{lem-time}, to prove (\ref{unif.1L}) it suffices to show that 
\begin{equation}
\Pbm\(\sup_{y\in D_{y_{0}}} \overline \MM_{y, \bar \ep_{y}, y_{0}, h_{L}, h_{L-1}}\(t_{L}\left(z-j+2M^{ 2} \right) \)\geq   {1 \over \pi}  t_{L}\left(z \right)   \)\leq ce^{-c'j^{2}}\label{unif.3}
\end{equation}
for $j\leq z+ML^{1/2} $ sufficiently large. Setting $\ep=\sup_{y\in D_{y_{0}}} \bar \ep_{y}$ and using our condition on $| \bar \ep_{y}- \bar \ep_{y'}| $ to control the denominator in (\ref{unif.1a}), we see that it suffices to show that
\begin{equation}
\Pbm\(\sup_{y\in D_{y_{0}}} \overline \MM_{y, \ep, y_{0}, h_{L},h_{L-1}}\(t_{L}\left(z-j+2M^{ 2} \right) \)\geq   {1 \over \pi}  t_{L}\left(z-M^{2}/2 \right)   \)\leq ce^{-c'j^{2}}.\label{unif.3a}
\end{equation}

 Abbreviating $Y^{(n)}_{y}= \overline \MM_{y, \ep, y_{0},h_{L},h_{L-1}}\(n \)$ where $n=t_{L}\left(z-j+2M^{ 2} \right)$ we have that 
\begin{eqnarray}
&&\Pbm\(\sup_{y\in D_{  y_{0}}} \overline \MM_{y, \ep, y_{0}, h_{L},h_{L-1}}\(t_{L}\left(z-j+2M^{ 2} \right) \)\geq   {1 \over \pi}  t_{L}\left(z-M^{2}/2   \right)   \)
\label{unif.6}\\
&&\leq \Pbm\( \overline \MM_{ y_{0},\ep, h_{L},h_{L-1}}\(t_{L}\left(z-j+2M^{ 2} \right) \)\geq   {1 \over \pi}  t_{L}\left(z  -j/2 -M^{2}/2 \right)   \)\nn\\
&&+\Pbm\( \sup_{y\in D_{  y_{0}}} |Y^{(n)}_{ y}-Y^{(n)}_{ y_{0}}|\geq   jL/2   \).\nn
\end{eqnarray}
As in the proof of Lemma     \ref{lem-Mbounds}, the first term on the right hand side is bounded by $ce^{-c'j^{2}}$
for $j\leq z+ML^{1/2} $ sufficiently large. We then bound
\bea
&&
\Pbm\( \sup_{y\in D_{  y_{0}}} |Y^{(n)}_{ y}-Y^{(n)}_{ y_{0}}|\geq   jL/2   \) \label{unif.20}\\
&&\leq \sum_{l=1}^{ \log_{2} L}\Pbm\( \sup_{y,y'\in D_{  y_{0}},\,d(y,y')\approx 2^{ -l}d_{0}h_{L}} |Y^{(n)}_{ y}-Y^{(n)}_{ y'}|\geq   jL/2 l^{ 2}  \) \nn\\
&&\leq \sum_{l=1}^{ \log_{2} L}2^{ 2l}\sup_{y,y'\in D_{  y_{0}},\,d(y,y')\approx 2^{ -l}d_{0}h_{L}} \Pbm\(   |Y^{(n)}_{ y}-Y^{(n)}_{ y'}|\geq   jL/2 l^{ 2}  \). \nn
\eea

It follows from Lemma \ref{lem-conta} with $n=t_{L}\left(z-j+2M^{ 2} \right)\sim 2L^{2}$ as above and $\th=j/2^{3/2}l^{ 2}, \bar d(y,y')=2^{ -l}d_{0}$ that for some $C_{0}>0$
\bea
&&
2^{ 2l}\sup_{y,y'\in D_{  y_{0}},\,d(y,y')\approx 2^{ -l}d_{0}h_{L}} \Pbm\(   |Y^{(n)}_{ y}-Y^{(n)}_{ y'}|\geq   jL/2 l^{ 2}  \)\label{unif.21}\\
&&\leq 2^{ 2l}\exp \( -C_{0} j^{ 2}2^{ l/2}/8 d^{1/2}_{0} l^{ 4}  \)\nn
\eea
whose sum over $l$ is bounded by $ce^{-c'j^{2}}$. In order to apply Lemma \ref{lem-conta} we have to verify that $\th\leq \sqrt{\bar d(y,y') n}/2$. In our situation this means that
$j/2^{3/2}l^{ 2}\leq 2^{ -l/2}d^{1/2}_{0}L/2$, for all $j\leq 2ML^{1/2}$. Thus we have to verify that 
$2^{1/2}M2^{ l/2}/l^{ 2}\leq d^{1/2}_{0} L^{1/2}$, which follows from the fact that $l\leq \log_{2} L$, $L$ is large and $d_{0}, M$ are fixed.

\subsection{The lower bound}

Recall the notation    $\mathcal{T}_{y, l}^{ 1}=\mathcal{T}_{y,l}^{x^{2},1}$ from the beginning of Section 3. Let $\tau_{y }$ be the time needed to complete $x^{2}$ excursions from  $\partial B_{d}\left(y,h_{ 1}\right)$ to $ \partial B_{d}\left(y,h_{0}\right)$, and set
\begin{equation}
\bar\mu_{ \tau_{y } } (   y,\ep) = \frac{1}{\om_{\ep}}\int_{0}^{ \tau_{y } }1_{\{   B_{d}(   y,\ep)\}}(  X_{t})\,dt.\label{tp.1y}
\end{equation}
Recall $F^{0}_{L}$ from (\ref{F0.def}).
We will prove the following analogue of Lemma \ref{lem-exclb0}.
\bl\label{lem-exclb0t}
There exists a $0<c<\ff$ such that for all $0<r_{0} $ sufficiently small,
  $L $ large,  all $0\leq z\leq \log L$,  and all 
$h_{L}/100\leq \ep_{y}\leq h_{L}$
\begin{equation}
 \Pbm\left[\sup_{y\in F^{0}_{L}} \bar \mu_{\tau_{y }}\(y,\ep_{y}\)\geq   {1 \over \pi}  t_{L}\left(z \right)\right]\geq {  (1+z)e^{ -2z} \over   (1+z) e^{ -2z}+c}.\label{goal.lbastt}
\end{equation}
\el
As before, the lower bound in (\ref{time.1}) will follow from this, and hence combined with (\ref{unif.1}) we see that for some $0<z_{0}$, and all $z_{0}\leq z\leq \log L$
\begin{equation}
 c'ze^{-2z}\leq\Pbm\( \exists y \mbox{ s.t.  } \bar \mu_{\tau}\(y,\ep_{y}\)\geq   {1 \over \pi}  t_{L}\left(z \right)\)\leq 
cze^{-2z}.\label{unif.1m}
\end{equation}
Combined with (\ref{con1.9}) it is easy to check that this implies Theorem \ref{theo-spheretight}.

To prove (\ref{goal.lbastt}) set
\begin{equation}
\wt \II_{y,z+d} = \II_{y,z+d}  
\,\cap  \{ \bar \mu_{\tau_{y }}\(y,\ep_{y}\)\geq   {1 \over \pi}  t_{L}\left(z \right) \} \label{ezlb.1}
\end{equation}
for some $d<\ff$ to be chosen shortly. We use the second moment method used in the proof of 
Lemma \ref{lem-exclb0}. Indeed, since $\wt \II_{y,z+d} \subseteq   \II_{y,z+d} $ all upper bounds needed follow from those used in the proof of Lemma \ref{lem-exclb0}, and it only remains to prove the appropriate lower bound for $\wt \II_{y,z+d}$.

As in (\ref{ty.1a})-(\ref{ty.1b}) we have
\begin{eqnarray}
&& \Pbm\(\wt \II_{y,z+d}\)\geq   \Pbm\(\wh \II_{y,z+d},\,\bar \mu_{\tau_{y }}\(y,h_{L}\)\geq   {1 \over \pi}  t_{L}\left(z \right)\) 
\label{ezlb.2a}\\
&& \hspace{1 in} -\sum_{k=L_{+}}^{L-d^{\ast}}\sum_{a\geq k_{L}^{ 1/4}} \Pbm  \left[\wh \II_{y,z+d}\bigcap \HH_{k,a}\,\bigcap \,W^{c}_{y,k}\(N_{k,a}\) \right].\nonumber
\end{eqnarray}
Using the Markov property and then the barrier estimate (\ref{14.6a}) of Appendix I,
\begin{eqnarray}
&& \hspace{.3 in} \Pbm\(\wh \II_{y,z+d},\,\bar \mu_{\tau_{y }}\(y,\ep_{y}\)\geq   {1 \over \pi}  t_{L}\left(z \right)\)
\label{ezlb.2}\\
&&\geq \Pbm\( \wh   \II_{y,z+d} \mbox{ and }
\,\overline \MM_{y, \ep_{y}, h_{L},h_{L-1}}\(t_{L}\left(z+d \right) \)\geq   {1 \over \pi}  t_{L}\left(z \right)\)   \nonumber\\
&&=\Pbm\(\wh \II_{y,z+d}\) \Pbm\(  
\,\overline \MM_{y, \ep_{y}, h_{L},h_{L-1}}\(t_{L}\left(z+d \right) \)\geq   {1 \over \pi}  t_{L}\left(z \right)\)  \nonumber\\
&&\geq \bar c   (1+z) e^{-2L} e^{-2(z+d )}\nn\\
&&\hspace{.5 in} \Pbm\( \overline \MM_{y, \ep_{y}, h_{L},h_{L-1}}\(t_{L}\left(z+d  \right) \)\geq   {1 \over \pi} \( t_{L}\left(z+d  \right)-dL\)  \)\nn\\
&&=\bar c   (1+z) e^{-2L} e^{-2(z+d )}\nn\\
&&\hspace{.5 in}\Pbm\(\overline \MM_{y, \ep_{y}, h_{L},h_{L-1}}\(t_{L}\left(z+d  \right) \)\geq   {1 \over \pi} \( 1-{dL \over t_{L}\left(z+d \right)}\)t_{L}\left(z+d \right)  \)\nn\\
&&\geq \bar c   (1+z) e^{-2L} e^{-2(z+d )}(1-e^{-c''{\(dL\)^{2} \over t_{L}\left(z+d \right)}}),\nn
\eea
where the last line used  (\ref{clt.1}). It should be clear from the structure of $t_{L}\left(z+d \right)$ that we can choose some $d<\ff$ so that $e^{-c''{\(dL\)^{2} \over t_{L}\left(z+d \right)}}\leq 1/2$ uniformly in  $0\leq z\leq \log L$.  Finally, after fixing such a $d$, we can show as in the proof of the first moment estimate in Section 3.1, that for   $d^{\ast}$ large enough, the last line in (\ref{ezlb.2a}) is much smaller than the last line of (\ref{ezlb.2}).
\qed

\subsection{The left tail}\label{sec-lowertail}

\bl\label{lem-exclb0tl}
There exists a $0<c<\ff$ such that for all $0<r_{0} $ sufficiently small,
  $L $ large,  all $0\leq z\leq \log L$,   and all 
$h_{L}/100\leq \ep_{y}\leq h_{L}$
\begin{equation}
 \Pbm\left[\sup_{y\in F^{0}_{L}} \bar \mu_{\tau_{y }}\(y, \ep_{y}\)\geq   {1 \over \pi}  t_{L}\left(-z \right)\right]\geq {  e^{ 2z} \over   e^{ 2z}+c}.\label{goal.lbasttl}
\end{equation}
\el
This will   complete the proof of Theorem \ref{theo-tight0} since,  as discussed right after the statement of  Lemma \ref{lem-exclb0}, the probability of completing  $x^{2}$ excursions from  $\partial B_{d}\left(y,h_{ 1}\right)$ to $ \partial B_{d}\left(y,h_{0}\right)$ before time $\tau$ for  all $y\in F^{0}_{L}$ is a strictly positive function of $r_{0}$ which goes to $1$ as $r_{0}\to 0$.  

The proof of Lemma \ref{lem-exclb0tl} is very similar to our proof of the  lower bound on the right tail, except we now have to change the upper barrier to allow for negative $z$. 
  Fix $|z|\leq \log L$.   We fix    $ \wh x>0$ once and for all.   
We abbreviate, 
\begin{equation}
\wh\beta_{z} \left(l\right)=f_{\wh x, \rho_{L}L+z}\left(l;L\right)= \wh x\(   1-\frac{l}{L}\)+ \(\rho_{L}l+z\frac{l}{L}\), \label{27.2}
\end{equation}
and     
\begin{equation}
\wh\gamma_{z, -}\left(l\right)=\wh\gamma_{z, -}\left(l,L,z\right)= \wh\beta_{z }\left(l\right)- l_{L}^{ 1/4}.\label{27.4}
\end{equation}
The barrier estimates needed are given in Lemma \ref{lem-NegBarrier}. We point out that the factors $ (1+z)$ which appear on the right hand side of (\ref{goal.lbastt}) but not (\ref{goal.lbasttl}) come from the difference in  the initial points of the barriers. \qed

  \section{  Interpolation}\label{sec-interp}

 Recall, (\ref{tp.1}), that
\begin{equation}
\bar\mu_{ \tau } (   y,\ep_{y}) = \frac{1}{\om_{\ep_{y}}}\int_{0}^{ \tau }1_{\{   B_{d}(   y,\ep_{y})\}}(  X_{t})\,dt,\label{int.1}
\end{equation}
where $\om_{\ep_{y}}=2\pi (   1-\cos \ep_{y})$, the area of $B_{d}(   y,\ep_{y})$, and, (\ref{clt.2}),
\begin{equation}
t_{L}\left(z\right)=2L\left(L-\log L+z\right).\label{clt.20}
\end{equation}

\bl\label{lem-interp} Assume that $d(y,y')\leq  a{h_{L} \over  L} $,  $|\ep_{y}-\ep_{y'}|\leq b \, {h_{L} \over  L} $, and $h_{L}/30\leq \ep_{y}, \ep_{y'}\leq 2h_{L+1}$.    We can find a $ d_{1} <\ff$  such that for all $L$ large and $z\leq \log L$, if
\be
\bar \mu_{\tau }\(y',\ep_{y'}\)\geq   {1 \over \pi}  t_{L}\left(z \right),\label{timet.90}
\ee
then for any  $c_{1}\geq 30(a+b) $,
\be
\bar \mu_{\tau }\(y,\(1+ c_{1}/L  \)\ep_{y}\)\geq   {1 \over \pi}  t_{L}\left(z -d_{1}\right).\label{timet.9}
\ee

\el

 {\bf  Proof: }
Under our assumptions,   for any $z\in   B_{d}(   y',\ep_{y'})$ we have
  $d(z,y)\leq d(z,y')+ d(y,y') \leq  \ep_{y'}+ a{h_{L} \over  L} \leq \(1+ {c_{1} \over L}  \)\ep_{y} $    so that
  \begin{equation}
    B_{d}(   y',\ep_{y'})\subseteq B_{d}\(   y,    \(1+ c_{1}/L  \)\ep_{y}\).\label{timet.2}
  \end{equation}
It follows that
  \begin{eqnarray}
\bar\mu_{ \tau } ( y',\ep_{y'})   &=& \frac{1}{\om_{\ep_{y'}}}\int_{0}^{ \tau }1_{\{   B_{d}(  y',\ep_{y'})\}}(  X_{t})\,dt
  \label{timet.3}\\
  &\leq &  \frac{1}{\om_{\ep_{y'}}}\int_{0}^{ \tau }1_{\{   B_{d}(   y, \(1+ c_{1}/L  \)\ep_{y})\}}(  X_{t})\,dt \nonumber\\
  &=&  \frac{\om_{ \(1+ c_{1}/L  \)\ep_{y}}}{\om_{\ep_{y'}}}\bar\mu_{ \tau } ( y, \(1+ c_{1}/L  \)\ep_{y}). \nonumber
  \end{eqnarray}
  
  Hence  
  \begin{equation}
\bar \mu_{\tau }\(y',\ep_{y'}\)\geq   {1 \over \pi}  t_{L}\left(z \right) \label{timet.ex}
  \end{equation} 
   implies that
  \begin{equation}
\bar\mu_{ \tau } ( y, \(1+ c_{1}/L  \)\ep_{y})\geq   \frac{\om_{\ep_{y'}}}{\om_{ \(1+ c_{1}/L  \)\ep_{y}}} {1 \over \pi}  t_{L}\left(z \right).  \label{timet.4}
  \end{equation}
But under our assumptions
 \begin{equation}
 \frac{\om_{\ep_{y'}}}{\om_{ \(1+ c_{1}/L  \)\ep_{y}}} = 1+O\(1/L\).\label{timet.6}
 \end{equation} 
This gives (\ref{timet.9}).\qed

\section{Green's functions and proof of Lemma \ref{lem-clt}}\label{sec-dev}

  Let $  G_{a}(x,y)$ denote the potential density for Brownian motion killed the first time it leaves $B_{e}(0,a)$, that is, the Green's function for $B_{e}(0,a)$. 
  We have, see \cite[Section 2]{DPRZint} or \cite[Chapter 2, (1.1)]{Doob},  
  \begin{equation}
  G_{a}(x,y)=-\frac{1}{\pi}\log |x-y|+\frac{1}{\pi}\log \(\frac{|y|}{a}|x-y_{a}^{\ast}|\), \hspace{.2 in}y\neq 0,\label{gr.1}
  \end{equation}
  where
  \begin{equation}
  y_{a}^{\ast}=\frac{a^{2}y}{|y|^{2}},\label{gr.2}
  \end{equation}
  and
  \begin{equation}
 G_{a}(x,0)=-\frac{1}{\pi}\log |x|+ \frac{1}{\pi}\log a.\label{gre.3}
  \end{equation}

   Let $v$ denote the south pole of $\S^{2}$.
   We claim that in the isothermal coordinates induced by stereographic projection $\si$, the Green's function for $\si\(B_{d}(   v,h(a))\)= B_{e}(   0,a) $ is just $G_{a}(x,y)$. To see this we must show that if $\De_{\S^{2}}$ is the Laplacian for $\S^{2}$ in isothermal coordinates and $dV(y)$
 is the volume measure, then
 \begin{equation}
\frac{1}{2}\De_{\S^{2}}\int G_{a}(x,y)f(y)\,dV(y) =-f(x)\label{13.1}
 \end{equation}
 for all continuous $f$ compactly supported in $B_{e}(   0,a) $. 
 
 For $x=\(   x_{1}, x_{2}\)$, let
\begin{equation}
g(x) =\frac1{(1+{\frac14}(x_{1}^{2}+x_{2}^{2}))^{2}}.\label{eu.2}
\end{equation}
As shown in \cite[Chapter 7, p. 6-9]{spivak}, the stereographic projection $\si$
 is an isometry if we give $R^{2}$ the metric
 \begin{equation}
g(x) \(dx_{1}\otimes dx_{1}+dx_{2}\otimes dx_{2}\). \label{eu.3}
 \end{equation} 

Because of (\ref{eu.3}) the Laplace-Beltrami operator takes the form
\begin{equation}
{1 \over g(x)}\(\partial^{2}_{x_{1}}+\partial^{2}_{x_{2}}\).\label{eu.1x}
\end{equation}
Thus, $\De_{\S^{2}}=\frac{1}{g}\De$ and $dV(y)=g(y)\,dy$, so that (\ref{13.1}) holds. 
 
{\bf  Proof of   Lemma \ref{lem-clt}: }  Let $\ep=h(\al)$ so that  $h(\al)\leq h_{k}$. lf $\tau_{h_{k-1} }$ is the first exit time of $B_{d}(v,h_{k-1} )$ and   $\rho_{h_{k}  }$ is uniform measure on $\partial B_{d}(   v,h_{k} )$, then by symmetry, for   any $z\in \partial B_{d}(   v, h_{k})$ 
\bea
J_{1}&=:&
  \Ebm^{z} \(   \int_{0}^{ \tau_{h_{k-1} } }1_{\{   B_{d}(   v, \ep)\}}(  X_{t})\,dt \)\label{gre.4}\\
  &=& \Ebm^{\rho_{h_{k}  }} \(   \int_{0}^{ \tau_{h_{k-1} } }1_{\{   B_{d}(   v,\ep)\}}(  X_{t})\,dt         \).\nn
\eea
 
 Since uniform measure $\rho_{h_{k} }$ on $\partial B_{d}(   v,h_{k} )$ goes over to uniform measure $\ga_{r_{k}}$ on $\partial B_{e}(   0,r_{k}) $, using the discussion at the beginning of this section  we have
 \begin{equation}
  J_{1}= \int_{B_{e}(0, \al )}\int_{\partial B_{e}(   0,r_{k})} G_{r_{k-1}}(x,y) \,d\ga_{r_{k}} (x)    g(y)\,dy.  \label{gre.4a}
 \end{equation}

We recall, \cite[Chapter 2, Prop. 4.9]{PS} or \cite[Chapter 1, (5.4), (5.5)]{Doob}, that
 \begin{equation}
 \int_{ \partial   B_{e}(  0,b) } \log \(|x-y|\) \,d\ga_{b} (x)=\log \(b\vee | y|\).\label{gr.22}
 \end{equation}
  This shows that for $y\in B_{e}(  0,r_{k}) $  
  \bea
 && \int_{ \partial   B_{e}(  0,r_{k}) } G_{r_{k-1}}(x,y) \,d\ga_{r_{k}}(x) \label{gr.23}\\
 &&=\frac{1}{\pi}\int_{ \partial   B_{e}(  0,r_{k}) } \(- \log |x-y|+\log \(\frac{|y|}{r_{k-1}}|x-y_{r_{k-1}}^{\ast}|\)\) \,d\ga_{r_{k}}(x) \nn\\
 &&=\frac{1}{\pi} \(- \log r_{k}+\log \(\frac{|y|}{r_{k-1}}|y_{r_{k-1}}^{\ast}|\)\)\nn\\
 &&=\frac{1}{\pi}\(-\log r_{k}+\log r_{k-1}\)=\frac{1}{\pi} \log (r_{k-1}/r_{k})= \frac{1}{\pi}. \nn
  \eea
Thus 
  \be
  J_{1}=\frac{1}{\pi}  \int_{B_{e}(0, \al )}g(y)\,dy=\frac{1}{\pi}   \mbox{Area }(B_{d}(  v,h(\al))) =\frac{1}{\pi}\om_{h(\al)}=\frac{1}{\pi}\om_{\ep}.\label{gre.5a} 
  \ee
It follows that for any $z\in \partial B_{d}(   v,h_{k})$
   \be 
  \Ebm^{z}\( \overline \MM_{v,\ep, h_{k}, h_{k-1}}(1)\)=\frac{1}{\pi}.
 \label{kap.130} 
 \ee
  
  By the Kac moment formula, for any $z\in \partial B_{d}(   v, h_{k})$, with $x=\si(z)$
  \begin{eqnarray} 
  && \Ebm^{z } \( \(  \int_{0}^{ \tau_{h_{k-1} } }1_{\{   B_{d}(   v,\ep)\}}(  X_{t})\,dt\)^{n}        \)
  \label{gre.6}\\
  &&=n! \int_{B^{n}_{e}(0, \al)}  G_{r_{k-1}}(x,y_{1})\prod_{j=2}^{n}G_{r_{k-1}}( y_{j-1}, y_{j})   \prod_{i=1}^{n}g(y_{i})\,dy_{i} \nonumber\\
  &&\leq c^{n } n! \int_{B^{n}_{e}(0, \al )}  G_{r_{k-1}}(x,y_{1})\prod_{j=2}^{n}G_{r_{k-1}}( y_{j-1}, y_{j})   \prod_{i=1}^{n} \,dy_{i} \nonumber\\
  &&\leq c^{n } n! \al^{2n}\(\log \(r_{k-1}/\al\)+c_{0}\)^{n} \nonumber
  \end{eqnarray}
where the last inequality follows as in the proof of \cite[Lemma 2.1]{DPRZint}. It follows that for any $z\in \partial B_{d}(   v, h_{k})$
 \be 
  \Ebm^{z}\(\(\overline \MM_{v,\ep, h_{k}, h_{k-1}}(1)\)^{n} \)\leq c^{n } n!  \(\log \(r_{k-1}/\al\)+c_{0}\)^{n}. 
 \label{kap.131} 
 \ee
 
  By (\ref{dr.1c}), our assumption that  $h_{k}/100\leq \ep\leq h_{k}$ implies that $e\leq r_{k-1}/\al\leq 200e$.  
 Using (\ref{kap.130}) and (\ref{kap.131}), our Lemma then follows as in  the proof of \cite[Lemma 2.2]{DPRZcover} 
\qed

  \section{Continuity Estimates}\label{sec-var}
 
 For fixed $u\in \S^{2}$,
let $\tau_{a}$ be the first exit time of $B_{d}(u,a)$ and let $\rho_{m}  $ be uniform measure on $\partial B_{d}(u,m)$.
  \bl \label{lem-varsphere}  
If $  d(u, v), d(u,\wt v)\leq d_{0}h_{L}/2$,  $d_{0}/L\leq \bar d=:d(v,\wt v)/h_{L}\leq d_{0}$,
  and $h_{L}/20\leq \ep\leq h_{L+1}$, then
 \be
 \Ebm^{\rho_{h_{L} }} \(   \int_{0}^{ \tau_{h_{L-1}} }1_{\{   B_{d}(   v,\ep)\}}(  X_{t})\,dt-\int_{0}^{ \tau_{h_{L-1}} }1_{\{   B_{d}(   \wt v,\ep)\}}(  X_{t})\,dt        \)=0,\label{12.1sp}
 \ee

  \be 
  \Ebm^{\rho_{h_{L} }} \(  \( \int_{0}^{ \tau_{h_{L-1}} }1_{\{   B_{d}(   v,\ep)\}}(  X_{t})\,dt-\int_{0}^{ \tau_{h_{L-1}} }1_{\{   B_{d}(\wt v,\ep)\}}(  X_{t})\,dt       \)^{2}\)\leq c\ep^{4}\bar d^{2},\label{12.2sp} 
 \ee 
  and
   \be 
\sup_{x\in \partial B_{d}(u,h_{L} )}  \Ebm^{x} \(  \( \int_{0}^{ \tau_{h_{L-1}} }1_{\{   B_{d}(   v,\ep)\}}(  X_{t})\,dt-\int_{0}^{ \tau_{h_{L-1}} }1_{\{   B_{d}(\wt v,\ep)\}}(  X_{t})\,dt       \)^{2}\)\leq c\ep^{4}\bar d^{2}.\label{12.2spi} 
 \ee 
  \el

   {\bf  Proof of Lemma \ref{lem-varsphere}: } 
 As in (\ref{gre.4})-(\ref{gre.4a}) we have
   \begin{eqnarray}
  J_{2}&=& \Ebm^{\rho_{h_{L}  }} \(   \int_{0}^{ \tau_{h_{L-1} } }1_{\{   B_{d}(   v,\ep)\}}(  X_{t})\,dt-\int_{0}^{ \tau_{h_{L-1} } }1_{\{   B_{d}(   \wt v,\ep)\}}(  X_{t})\,dt        \)
\nn\\
  &=&    \int\int G_{r_{L-1}}(x,y) \,d\ga_{r_{L}} (x)\,d\mu_{v,\wt v}(y),  \label{12.3sp}
  \end{eqnarray}
    where  
  \begin{equation}
 d\mu_{v,\wt v}(y)=\(1_{\{\si( B_{d}(   v,\ep))   \}}-1_{\{ \si( B_{d}(  \wt v,\ep)) \}}\)(y)g(y)\,dy. \label{gr.4sp}
  \end{equation}
  Then by (\ref{gr.23})-(\ref{gre.5a}) we have that 
  \bea
  J_{2}&=&\frac{1}{\pi} \int  d\mu_{v,\wt v}(y)\label{gr.24m}\\
  &=&\frac{1}{\pi}  \(\mbox{Area }(B_{d}(   v,\ep))-\mbox{Area }(B_{d}( \wt  v,\ep))\)=0,\nn
  \eea
 since all balls of radius $\ep$  on the sphere have area $\om_{\ep}=2\pi (   1-\cos \ep)$. This completes the proof of (\ref{12.1sp}). 
  
  We next observe that
  \begin{eqnarray}
  &&     \Ebm^{\rho_{h_{L}  }} \(  \( \int_{0}^{ \tau_{h_{L-1} } }1_{\{   B_{d}(   v,\ep)\}}(  X_{t})\,dt-\int_{0}^{ \tau_{h_{L-1} } }1_{\{   B_{d}(\wt v,\ep)\}}(  X_{t})\,dt       \)^{2}\)
  \label{12.5sp}\\
  &&  =2 \int\int\int G_{r_{L-1}}(x,y) G_{r_{L-1}}(y,z)\,d\ga_{r_{L}} (x)\,d\mu_{v,\wt v}(y)\,d\mu_{v,\wt v}(z)\nonumber\\
  &&  =\frac{2}{\pi}  \int\int  G_{r_{L-1}}(y,z) \,d\mu_{v,\wt v}(y)\,d\mu_{v,\wt v}(z)\nonumber
  \end{eqnarray}
  as above. 
  
  We note that  for $b<a$
   \begin{equation}
  G_{a}(bx,by)=-\frac{1}{\pi}\log \( b\,|x-y|\)+\frac{1}{\pi}\log \(b\,\frac{|y|}{a/b}|x-y_{a/b}^{\ast}|\)=G_{a/b}(x,y),\label{gr.5}
  \end{equation}
  since 
    \begin{equation}
  (by)_{a}^{\ast}=\frac{a^{2}by}{b^{2}|y|^{2}}=b  y_{a/b}^{\ast}.\label{gr.6}
  \end{equation}

  Using this to scale by $r_{L}$  we see that 
    \begin{equation}
  \int\int  G_{r_{L-1}}(y,z) \,d\mu_{v,\wt v}(y)\,d\mu_{v,\wt v}(z)=r_{L}^{4}\int\int  G_{e}(y,z)  \,d\mu_{L, v,\wt v}(y)\,d\mu_{L,v,\wt v}(z),\label{12.6sp}
  \end{equation}
    where
    \bea
 d\mu_{L,v,\wt v}(y)&=&\(1_{\{\si( B_{d}(   v,\ep))   \}}-1_{\{ \si( B_{d}(  \wt v,\ep)) \}}\)(r_{L} y)g(r_{L} y)\,dy\nn\\
 &=&\(1_{\{\frac{1}{r_{L}}\si( B_{d}(   v,\ep))   \}}-1_{\{ \frac{1}{r_{L}}\si( B_{d}(  \wt v,\ep)) \}}\)( y)g(r_{L} y)\,dy. \label{gr.4spa}
  \eea
For $y$ in our range we have $g(r_{L} y)=1+O(\ep)$, and it is easy to check that up to errors of order $\ep$, $\frac{1}{r_{L}}\si( B_{d}(   v,\ep)) $ and $\frac{1}{r_{L}}\si( B_{d}( \wt  v,\ep)) $ can be replaced by  $ B_{e}(   v',\ep/r_{L})$ and $  B_{e}(  v'-(0,\de\ep/r_{L} ), \ep/r_{L})$ for some $  v' $ with $  |v'|\leq d_{0}$ and $0<\de\leq c_{2}\bar d$.  

Hence with
    \begin{equation}
 d\nu_{v'}(y)=\(1_{\{   B_{e}(   v',\ep/r_{L})\}}-1_{\{   B_{e}(   v'-(0,\de \ep/r_{L}),\ep/r_{L})\}}\)(y)\,dy, \label{gr.9}
  \end{equation}
  it remains to show that
  \begin{equation}
  \int\int  G_{e}(y,z) \,d\mu_{v'}(y)\,d\mu_{v'}(z)\leq C\de^{2}.\label{gr.30}
  \end{equation}
  
   The symmetric difference of $ B_{e}(   v',s)$ and $  B_{e}(  v'-(0,\de s), s)$
  consist of two disjoint  pieces we denote by $A,B$. They have the same area
  \begin{equation}
  \mbox{Area}\(A\)=2 s^{2}\arcsin \({\de \over 2 }\)+{\de s\over 2 }\sqrt{(4-\de^{2}) s^{2}}\asymp \de s^{2}.
  \label{12.3x}
  \end{equation}
  
  We observe   that
  \begin{equation}
 \lim_{y\to 0}\frac{|y|}{a}|z-y_{a}^{\ast}|= \lim_{y\to 0}\frac{|y|}{a}\(| y_{a}^{\ast}|+O(1)\)= \lim_{y\to 0}\frac{|y|}{a}\(\frac{a^{2}|y|}{|y|^{2}}+O(1)\)=a.\label{gr.4a}
  \end{equation}
 It follows that for $y,z$ in our range, $\log \(\frac{|y|}{e}|z-y_{e}^{\ast}|\)$ is bounded, hence to prove (\ref{gr.30}) it suffices to show that
  \begin{equation}
  \int\int \big |\log |y-z|\big |\,d\mu_{v'}(y)\,d\mu_{v'}(z)\leq C\de^{2}.\label{gr.31}
  \end{equation}
It is then easy to see that we need only show that
  \begin{equation}
  \int_{A}\int_{A}  \big |\log |y-z|\big |\,d y\,d z\leq C\de^{2}.\label{gr.32}
  \end{equation}
 It is also clear that we only need to consider $|y-z|\leq 1/2$. Writing $y=(y_{1},y_{2}), z=(z_{1},z_{2})$ we see that
 \begin{eqnarray}
 &&\int_{A}\int_{A}  \big |\log |y-z|\big | 1_{\{|y-z|\leq 1/2\}}\,d y\,d z
 \label{gr.33}\\
 &&  \leq  \int_{[0,1]\times [0,\de]}\int_{[0,1]\times [0,\de]}\big |\log |y-z|\big | 1_{\{|y-z|\leq 1/2\}}\,d y_{1}\,dy_{2}\, d z_{1}\,dz_{2}\nonumber\\
 &&  \leq  \int_{[0,1]\times [0,\de]}\int_{[0,1]\times [0,\de]} \big |\log |y_{1}-z_{1}|\big |  \,d y_{1}\,dy_{2}\, d z_{1}\,dz_{2}\leq C\de^{2},\nonumber
 \end{eqnarray}
 which completes the proof of (\ref{gr.32}).
 
 To obtain (\ref{12.2spi}), arguing as before we need to show that
 \begin{equation}
K_{1}=\sup_{x\in \partial B_{e}(0,r_{L})} \int\int  G_{r_{L-1}}(x,y) G_{r_{L-1}}(y,z) \,d\mu_{v,\wt v}(y)\,d\mu_{v,\wt v}(z)\leq c\ep^{4}\de^{2}. \label{gr.40}
 \end{equation}
 Scaling in $r_{L} $ as before shows that
 \begin{equation}
 K_{1}=r_{L} ^{4}   \sup_{x\in \partial B_{e}(0,1)} \int\int  G_{e}(x,y) G_{e}(y,z) \,d\mu_{L, v,\wt v}(y)\,d\mu_{L, v,\wt v}(z)  \label{gr.41}
 \end{equation}
 But for $y$ in our range, $G_{e}(x,y)$ is bounded uniformly in $x\in \partial B_{e}(0,1)$, so that (\ref{12.2spi}) follows as before. 
  \qed

The same proof shows that
   \bea 
&&\sup_{x\in \partial B_{d}(u,h_{L})}  \Ebm^{x} \(  \( \int_{0}^{ \tau_{h_{L-1}} }1_{\{   B_{d}(   v,\ep)\}}(  X_{t})\,dt-\int_{0}^{ \tau_{h_{L-1}} }1_{\{   B_{d}(\wt v,\ep)\}}(  X_{t})\,dt       \)^{2n}\)\nn\\
&&\hspace{3.5 in} \leq (2n)!c_{1}^{2n}\ep^{4n}\bar d^{2n}.\label{12.2nspg} 
 \eea 
 and hence by the Cauchy-Schwarz inequality
    \bea 
&&\hspace{-.25 in}\sup_{x\in \partial B_{d}(u,h_{L})}  \Ebm^{x} \(  \Big |\frac{1}{\om_{\ep}} \int_{0}^{ \tau_{h_{L-1}} }1_{\{   B_{d}(   v,\ep)\}}(  X_{t})\,dt-\frac{1}{\om_{\ep}}\int_{0}^{ \tau_{h_{L-1}} }1_{\{   B_{d}(\wt v,\ep)\}}(  X_{t})\,dt     \Big |^{n}\)\nn\\
&&\hspace{3.5 in} \leq n!c_{2}^{n} \bar d^{n}.\label{12.2nsp} 
 \eea
 
 Recall (\ref{unif.1a}) and set
 \begin{equation}
Y^{(n)}_{y}= \overline \MM_{y, \ep, u, h_{L},h_{L-1}}\(n \) \nn
 \end{equation}

\bl\label{lem-conta}
For some $d_{0}>0$ we can find  $C_{0}>0$ such that, if $  d(u, v), d(u,\wt v)\leq d_{0}h_{L}/2$,  $d_{0}/L\leq \bar d(v,\wt v)=:d(v,\wt v)/h_{L}\leq d_{0}$,
  $h_{L}/20\leq \ep\leq h_{L+1}$,  and
 $\th\leq \sqrt{\bar d(v,\wt v) n}/2$, then  
\begin{equation}
\Pbm \(|Y^{(n)}_{v}-Y^{(n)}_{\wt v}|
\geq \th\sqrt{n}\)\leq e^{-C_{0}\th^{2}/\bar d^{1/2}(v,\wt v)}.\label{cont.1}
\end{equation}
\el

{\bf  Proof of Lemma \ref{lem-conta}: }We follow the proof of \cite[Lemma 5.1]{BRZ}. 

Let   $T_{i}  $ denote the successive excursion times $T_{\partial B_{d}(u,h_{L})}\circ \th_{T_{\partial B_{d}(u,h_{L-1})}}$ from $\partial B_{d}(u,h_{L-1})$ to $\partial B_{d}(u,h_{L})$ and set
\begin{equation}
Y_{v,i}=  \frac{1}{\om_{\ep}} \int_{0}^{ \tau_{h_{L-1}} }1_{\{   B_{d}(   v,\ep)\}}(  X_{t+T_{i}})\,dt,\label{cont.2}
\end{equation}
so that
\begin{equation}
Y^{(n)}_{v}=\sum_{i=1}^{n}Y_{v,i}.\label{cont.3}
\end{equation}
Let $J$ be a geometric random variable with success parameter $p_{3}>0$,
independent of $\{Y_{v,i},Y_{\wt v,i}\}$. It follows from (\ref{12.2nsp}) and the proof of \cite[Corollary 5.3]{BRZ} that, abbreviating  $\bar  d=\bar d(v,\wt v)$, if
 $c_{2}\bar d \lambda\le  p_{3}/2  $ then  for some $c_{4}$
\bea
&&\sup_{x\in \partial B_{d}(u,h_{L})}\mathbb{E}^{x}\left(\exp\left(\lambda\sum_{i=1}^{J-1}\left|\left(Y_{v,i}-Y_{\wt v,i}\right)\right|\right)\right)\leq e^{c_{4}\bar d\lambda /p_3  },\label{basicexp1}
\eea
and from (\ref{12.2nsp}) together with (\ref{12.1sp}) and the proof and notation of \cite[Lemma 5.5]{BRZ} it follows that, after perhaps enlarging $c_{4}$
\bea
&& \mathbb{E}\left(\exp\left(\lambda\sum_{i=J_{1}}^{J_{2}-1} \left(Y_{v,i}-Y_{\wt v,i}\right) \(X_{\cdot}^{i}\)\right)\right)\leq e^{c_{4}\(\bar d\lambda /p_3\)^{2}  }.\label{basicexp2}
\eea
 Instead of \cite[(5.33)]{BRZ} we set
 \begin{equation}
 \de=\th/\sqrt{\bar dn}\leq 1/2.\nn
 \end{equation}
It then follows from the proof of \cite[Lemma 5.1]{BRZ} that for $c_{2}\bar d\la\leq p_{3}/2$
\bea
&&
\Pbm \(|Y^{(n)}_{v}-Y^{(n)}_{\wt v}|
\geq \th\sqrt{n}\)\label{cont.4}\\
&&\leq e^{-\bar c \th^{2}/\bar d}+\exp\(c_{4}\la^{2}\bar d^{2}n/p_{3}+2c_{4}\la   \th\sqrt{\bar d n}-\la\th \sqrt{n}  \)\nn
\eea
By taking  $d_{0}$ sufficiently small we can be sure that $c_{2}\bar d\leq 1$, so the above holds for any $ \la\leq p_{3}/2$.
If we set
\begin{equation}
\la=p_{3}\th/\sqrt{\bar dn}\le  p_{3}/2\nn
\end{equation}
we see that
\begin{eqnarray}
&&\exp\(c_{4}\la^{2}\bar d^{2}n/p_{3}+2c_{4}\la   \th\sqrt{\bar d n}-\la\th \sqrt{n}  \)
\nn\\
&& = \exp\(c_{4}p_{3}\th^{2}\bar d+2c_{4}p_{3} \th^{2}-p_{3} \th^{2}/\sqrt{\bar d}  \), \nonumber
\end{eqnarray}
which completes the proof of (\ref{cont.1}) for $d_{0}$ sufficiently small.\qed

\section{  From the sphere to the plane, and back}\label{sec-euclid}

Using (\ref{eu.1x}) it follows from  \cite[Chapter 5, Theorem 1.9]{RY}, that we can find a planar Brownian motion $W_{t}$
such that in  the isothermal coordinates induced by stereographic projection,
\begin{equation}
 X_{t}=W_{U_{t}}, \hspace{.1 in}\mbox{ where  }\hspace{.1 in}U_{t}=\int_{0}^{t}{1 \over g(X_{s})}\,ds.\label{eu.1y}
\end{equation}

We take the $v$ of this paper to be $v=(0,0,0 )$.
Let 
\be
D_{\ast} =\si \(   B_{d}\(   v,r^{ \ast}\)\)=B_{e}\( ( 0,0),2 \tan (r^{ \ast}/2)\).\label{dast}\ee
For the last equality see \cite[(2.4)]{BRZ}.
If $\th$ is the first hitting time of $\partial D_{\ast} $ by $W_{t}$, then under the coupling (\ref{eu.1y}) we see that $\th=U_{\tau}$. Set
\begin{equation}
\mu_{ \th}(x,\ep)=\frac{1}{\pi \ep^{ 2} }\int_{0}^{  \th }1_{\{   B_{e}(   x,\ep)\}}(  W_{t})\,dt.\label{eu.1}
\end{equation}
 
 \bl\label{lem-disc}
 For some 
$-\ff<d_{1}, d_{2}, d_{3},  d_{4}<\ff$, all $x\in D_{\ast} $ and all $\ep$ sufficiently small
\begin{equation}
\mu_{ \th}(x,\ep)\leq ( 1+d_{1}   \,\,\ep) \,\, \bar \mu_{\tau}(x,g^{1/2}(x) \ep ( 1+d_{2}\ep)),\label{eu.6a}
\end{equation}
and
\begin{equation}
\mu_{ \th}(
x,\ep)\geq ( 1+d_{3}  \,\,\ep) \,\, \bar \mu_{\tau}(x,g^{1/2}(x) \ep ( 1+d_{4}\ep)).\label{eu.6b}
\end{equation}
 \el

{\bf  Proof of Lemma \ref{lem-disc}: } We first note that for $\ep$ sufficiently small, we can find $c_{1}
< c_{2}$ such that uniformly in $x'\in  B_{e}(   x,2\ep)$ and $x\in D_{\ast} $
\begin{equation}
g(x)(1+c_{1}\ep)\leq g(x')\leq g(x)(1+c_{2}\ep).\label{eu.0}
\end{equation}
 For $x'\in  B_{e}(   x,\ep)$, with $x_{t}=x+t(x'-x)$
 \begin{equation}
d(x,x')\leq \int_{0}^{1}g^{1/2}(x_{t})|x'-x|\,dt\leq g^{1/2}(x)|x'-x   |( 1+c_{3}\ep). \label{eu.3a}
 \end{equation}
 Hence 
 \begin{equation}
 B_{e}(   x,\ep)\subseteq B_{d}(   x, g^{1/2}(x) \ep ( 1+c_{3}\ep)).\label{eu.4}
\end{equation}
 Similarly, for some  $c_{4}
< c_{3}$
\begin{equation}
 B_{e}(   x,\ep)\supseteq B_{d}(   x, g^{1/2}(x) \ep ( 1+c_{4}\ep)).\label{eu.01x}
\end{equation}
 
Consider 
\begin{equation}
\int_{0}^{ \tau }1_{\{   B_{d}(   x, g^{1/2}(x) \ep ( 1+c_{3}\ep))\}}( W_{U_{t}})\,dt.\label{eu.20}
\end{equation}
By the nature of $U_{t}$ in (\ref{eu.1y}) it follows that whenever the path $W_{U_{t}}$ enters 
$ B_{d}(   x, g^{1/2}(x) \ep ( 1+c_{3}\ep))$ it is slowed by a variable factor between ${1 \over g  (x)  ( 1+c_{5}\ep)} $
and ${1 \over g  (x)  ( 1+c_{6}\ep)} $. Hence the amount of time spent in $B_{d}(   x, g^{1/2}(x) \ep ( 1+c_{3}\ep))$ during each incursion is multiplied by a variable factor between $ g  (x)  ( 1+c_{7}\ep)$
and $  g  (x)  ( 1+c_{8}\ep)$. Thus 
\begin{equation}
\int_{0}^{ \tau }1_{\{   B_{d}(   x, g^{1/2}(x) \ep ( 1+c_{3}\ep))\}}( W_{U_{t}})\,dt\geq   g  (x)  ( 1+c_{7}\ep)  \int_{0}^{ \th }1_{\{   B_{d}(   x, g^{1/2}(x) \ep ( 1+c_{3}\ep))\}}( W_{t})\,dt,\label{eu.20a}
\end{equation}
and 
\begin{equation}
\int_{0}^{ \tau }1_{\{   B_{d}(   x, g^{1/2}(x) \ep ( 1+c_{3}\ep))\}}( W_{U_{t}})\,dt\leq   g  (x)  ( 1+c_{8}\ep)  \int_{0}^{ \th }1_{\{   B_{d}(   x, g^{1/2}(x) \ep ( 1+c_{3}\ep))\}}( W_{t})\,dt,\label{eu.20b}
\end{equation}

It follows from (\ref{eu.4}) and  (\ref{eu.20a}) that
\begin{eqnarray}
&&\int_{0}^{  \th }1_{\{   B_{e}(   x,\ep)\}}(  W_{t})\,dt
\label{eu.5}\\
&&\leq   \int_{0}^{ \th }1_{\{   B_{d}(   x, g^{1/2}(x) \ep ( 1+c_{3}\ep))\}}( W_{t})\,dt \nonumber\\
&&\leq  {1 \over g  (x)  ( 1+c_{7}\ep)}  \int_{0}^{\tau }1_{\{   B_{d}(   x,g^{1/2}(z) \ep ( 1+c_{3}\ep))\}}( W_{U_{t}})\,dt.  \nonumber
\end{eqnarray}
Since $\om_{\de}=2\pi\(   1-\cos \(   \de\)\) $, we see that if we set $\de_{x}= g^{1/2}(x) \ep ( 1+c_{3}\ep) $, then uniformly in $x\in D_{\ast}^{2}$  and sufficiently small $\ep$
\begin{equation}
 (1+f'_{0}\ep)     \leq  \frac{\om_{\de_{x}}}{\pi g\(   x\) \ep^{2}}\leq (1+f_{0}\ep),\nn
\end{equation}
  so that by (\ref{eu.4}) and  (\ref{eu.5})
\begin{eqnarray}
\mu_{ \th}(x,\ep)&=&\frac{1}{\pi \ep^{ 2} }\int_{0}^{  \th }1_{\{   B_{e}(  x,\ep)\}}(  W_{t})\,dt
\nn\\
&\leq &   \frac{1}{\pi \ep^{ 2} g  (x)  ( 1+c_{7}\ep)}\int_{0}^{\tau }1_{\{   B_{d}(   x,g^{1/2}(z) \ep ( 1+c_{3}\ep))\}}( W_{U_{t}})\,dt\nonumber\\
&= &   \frac{1}{  ( 1+c_{7}\ep)}\,\,\frac{\om_{\de_{x}}}{\pi g\(   x\) \ep^{2}}\,\,\bar \mu_{\tau}(x,g^{1/2}(x) \ep ( 1+c_{3}\ep))\nonumber\\
&\leq &( 1+\wh d   \,\,\ep) \,\, \bar \mu_{\tau}(x,g^{1/2}(x) \ep ( 1+c_{3}\ep)),\label{eu.6}
\end{eqnarray}
where
\begin{equation}
( 1+\wh d \,\,\ep)={1+f_{0}\ep \over 1+c_{7}\ep}.\nn
\end{equation}
This completes the proof of (\ref{eu.6a}).

 The lower bound (\ref{eu.6b}) is proven similarly using (\ref{eu.01x}) and (\ref{eu.20b}).  
\qed

 \bl\label{lem-timeu}
We can find $0<c,c',z_{0}<\ff$ such that for $L $ large and  all $\frac{1}{12} h_{L }\leq \ep\leq \frac{1}{3}h_{L }$ and $z_{0}\leq z\leq \log L$, 
\begin{equation}
c' ze^{-c^{\ast}\sqrt{\pi}z}\leq \Pbm\(\sqrt{\sup_{y }\mu_{ \th } (   y, \ep) }\geq m_{\ep}+ z \)\leq 
cze^{-c^{\ast}\sqrt{\pi}z}.\label{timeu.1}
\end{equation}
\el

{\bf  Proof of Lemma \ref{lem-timeu}: }We consider the upper bound. By (\ref{eu.6a}) 
it suffices to show that
\begin{equation}
 \Pbm\(\sqrt{\sup_{y }\bar \mu_{\tau}(y,g^{1/2}(y) \ep ( 1+d_{2}\ep))}\geq m_{\ep}+ z \)\leq 
cze^{-c^{\ast}\sqrt{\pi}z}.\label{timeu.2}
\end{equation}

Since for $\ep$ in our range
\begin{equation}
\(   m_{\ep}+z\)^{ 2}={1 \over \pi}  t_{L}\left(\sqrt{2\pi}\,\,z +O(1)\right),\nn
\end{equation}
(compare (\ref{con1.9})), (\ref{timeu.2}) follows from Lemma \ref{lem-unif} once we verify the condition that 
$|\ep_{y}-\ep_{y'}|\leq C\,d(y,y')/L$, where now $\ep_{y}=g^{1/2}(y) \ep ( 1+d_{2}\ep))$. This follows easily since $g$ is smooth and we can assume that $\frac{4}{5}\leq g^{1/2}(y)\leq 1$. We also point out that for $\frac{1}{12} h_{L }\leq \ep\leq \frac{1}{3}h_{L }$ and $L$  large we have $\frac{1}{20} h_{L }\leq \ep_{y}\leq  h_{L +1}$.

The lower bound is similar.
 \qed
 
 We note that $\frac{1}{3}h_{L+1 }=\frac{1}{3e}h_{L }\geq \frac{1}{12} h_{L }$, so all $\ep$ are covered by Lemma \ref{lem-timeu}. 
 
 Lemma \ref{lem-timeu} is the analog of Theorem \ref{theo-planetight}, but where now $\th$  is the first hitting time of $\partial D_{\ast} $, see (\ref{dast}). Theorem \ref{theo-planetight} then follows by Brownian scaling. To spell this out for later use, let  $\th_{a}$  be the first hitting time of $\partial B_{e} (0,a)$ and set 
\begin{equation}
 \mu_{a } (   x,\ep) =\frac{1}{\pi\ep^{ 2} }\int_{0}^{  \th_{a} }1_{\{   B(   x,\ep)\}}(  W_{t})\,dt.\label{plane.1h}
\end{equation}
Then it follows from  Brownian scaling that for any $a,b>0$,
\begin{equation}
\{ \mu_{a } (   x,\ep_{x});\,x, \ep_{x}\}\overset{\mbox{law}}{=}\{ \mu_{ba } (   bx,b\ep_{x});\,x, \ep_{x}\}.\label{pi.1}
\end{equation}

The left tail and then Theorem \ref{theo-tight0plane} can be proven similarly.

\subsection{$0<r^{\ast}<\pi$}

We first note the following extension of Lemma \ref{lem-timeu}.
 \bl\label{lem-timeua}
We can find $0<c,c',z_{0}<\ff$ such that for $L $ large and  all $\frac{1}{12} h_{L }\leq \ep_{y}\leq \frac{1}{3}h_{L }$ with $|\ep_{y}-\ep_{y'}|\leq C |y-y'|/L$ and $z_{0}\leq z\leq \log L$, 
\begin{equation}
c' ze^{-2z}\leq \Pbm\( \sup_{y }\mu_{ \th } (   y, \ep_{y})  \geq {1 \over \pi}  t_{L}\left(z \right) \)\leq 
cze^{-2z}.\label{timeu.1a}
\end{equation}
\el
This follows as in  the proof of Lemma \ref{lem-timeu}, once we observe that in Lemma \ref{lem-disc} we can allow the $\ep$ to depend on $x$.

It follows from (\ref{pi.1}) that  for any fixed $a>0$,  Lemma \ref{lem-timeua} holds with $\th$ replaced by $\th_{a}$.

We now show that Theorem \ref{theo-spheretight} holds for any $0<r^{\ast}<\pi$. This is done by using  Lemma \ref{lem-disc}.  That is, with $a=2 \tan (r^{ \ast}/2)$ we have that
for some 
$-\ff<d_{1}, d_{2}, d_{3},  d_{4}<\ff$, all $x\in D_{a} $ and all $\ep$ sufficiently small
\begin{equation}
\bar \mu_{\tau}(x,\ep)\leq ( 1+d_{1}   \,\,\ep) \,\, \mu_{ \th_{a}}(x,g^{-1/2}(x) \ep ( 1+d_{2}\ep)),\label{eu.6am}
\end{equation}
and
\begin{equation}
 \bar \mu_{\tau}(
x,\ep)\geq ( 1+d_{3}  \,\,\ep) \,\,\mu_{ \th_{a}}(x,g^{-1/2}(x) \ep ( 1+d_{4}\ep)),\label{eu.6bm}
\end{equation}
Theorem \ref{theo-spheretight} then follows from Lemma \ref{lem-timeua} just as  Lemma \ref{lem-timeu} followed from  Lemma \ref{lem-unif}. Theorem \ref{theo-tight0} can be proven similarly.

\section{Appendix I: Barrier estimates}\label{sec:BoundaryCrossing}

In what follows, we use the notation $H_{y,\delta}=[y,y+\delta]$ from
\cite{BRZ2}. The following is a variant of \cite[Theorem 1.1]{BRZ2}, which can be proven similarly. Re set
\begin{equation}
f_{a, b}\left(l;L\right)=a+(b-a)\frac{l}{L}.\label{fdef.8}
\end{equation}

\bt
\label{thm: GW Barrier}
 a)  For all fixed $\delta>0, C\geq 0$,  $\eta>1$  and 
$\varepsilon\in\left(0,\frac{1}{2}\right)$ we have,  
uniformly
in  $\sqrt{2}\le x, y\le  \eta  L$ 
 such that ${x^{2}}/{2}\in\mathbb{N}$, any $0\le x\le a, \,\,0\le y\le b$, 
that
\begin{eqnarray}
&&\Pgw_{x^{2}/2}\left(\sqrt{2T_{l}}\le f_{a, b}\left(l;L\right)+Cl_{L}^{\frac{1}{2}-\varepsilon},l=1,\ldots,L-1,\sqrt{2T_{L}}\in
 H_{y,\delta}\right) \nonumber\\
&&\hspace{2 in}\;\le c\frac{\left(1+a-x\right)\left(1+b-y\right)}{L} \sqrt{
 \frac{x}{yL}}e^{-\frac{\left(x-y\right)^{2}}{2L}}.
\label{eq: GW curved barrier upper bound}
\end{eqnarray}

 b) Let $\mbox{Tube}_{C,\tilde{C} }(   l;L)= [f_{x, y}\left(l;L\right)-\tilde{C}l_{L}^{\frac{1}{2}+\varepsilon}, f_{a, b}\left(l;L\right)-Cl_{L}^{\frac{1}{2}-\varepsilon} ]\ $. If,  in addition to the conditions in part a),  
we also have  $\left(1+a-x\right)\left(1+b-y\right)\le  \eta L$, $\max (xy, |y-x|)\geq    L/\eta$
and $\left[y,y+\delta\right]\cap\sqrt{2\mathbb{Z}}\ne\emptyset$, and  
 $\mbox{Tube}_{C,\tilde{C} }(   l;L)\cap\sqrt{2\mathbb{N}}\ne\emptyset$   for all l=1,\ldots,L-1,  then 
\bea
&&\Pgw_{x^{2}/2}\left( 
\sqrt{2T_{l}}\in \mbox{Tube}_{C,\tilde{C} }(   l;L), l=1,\ldots,L-1,  
\sqrt{2T_{L}}\in H_{y,\de} \right)\nn\\
&&\hspace{ 1in} \;\ge c\frac{\left(1+a-x\right)\left(1+b-y\right)}{L}\times\left(\sqrt{
  \frac{x}{yL}}\wedge1\right) e^{-\frac{\left(x-y\right)^{2}}{2L}},
\label{eq: gw lower bound}
\eea
 and the estimate is uniform in such $x,y,a,b$ and all $L$. 
 
 Similar results hold if we delete the barrier condition on some fixed finite interval.
\et

For the last statement, we simply note that following the proof of
 \cite[Lemma 2.3]{BRZ2}
we can show that the analogue of \cite[Theorem 1.1]{BRZ2} holds where we skip some fixed finite interval.

Recall that 
\begin{equation}
\rho_{L} = 2-\frac{\log L}{L}, \hspace{.2 in}\alpha_{  z, \pm}\left(l\right)=\alpha\left(l,L, z\right)= \rho_{L}l +z \pm l_{L}^{ 1/4}.\label{18.01}
\end{equation}
 

\bl\label{prop:BarrierSecGWPropUB}
Let $m=k_{y}+1\leq \log L$.
For any   $k\geq L-\(   \log L\)^{ 4}$, $0\leq j\leq \al_{z, +}(k)/2$ and $z\leq \log L$   
\begin{eqnarray}
&&
\Pbm\left[\sqrt{2T_{y,l}^{\tau}} \le  \al_{z, +}(   l),\,l=m,\ldots,k-1;\,  \sqrt{2T_{y,k}^{\tau}}\in I_{\al_{z, +}(k)-j}\right]\label{18.20}\\
&&\nn
\leq  ce^{-2k-2z-2k^{1/4}_{L}+2j } \times m^{2}\left(1+z+m + k^{ 1/4}_{L}\right)\left(1+j  \right).\end{eqnarray}
\el

{\bf   Proof of Lemma (\ref{prop:BarrierSecGWPropUB}): }  Using the Markov property, the probability in (\ref{18.20}) is bounded by
\begin{eqnarray}
&& \sum_{s=0}^{\al_{z, +}(m)}\Pbm\left[\sqrt{2T_{y,m}^{\tau}} \in I_{s}\right] \nn\\
&&\hspace{1 in}   \times  \sup_{x\in I_{s}}  \Pbm\left[\sqrt{2T_{y,l}^{m,x^{2}/2}} \le  \al_{z, +}(   l),\,l=m+1,\ldots,k-1;\,\right.
\nn\\
&&\left.\hspace{2.8 in} \sqrt{2T_{y,k}^{m,x^{2}/2}}\in I_{\al_{z, +}(k)-j}\right], \label{18.20f}
\end{eqnarray}
and using the fact that $T_{y,m}^{\tau}\leq T_{y,m}^{m,0}$ and (\ref{1.1}), we see that (\ref{18.20f}) is bounded by
\begin{eqnarray}
&& \sum_{s=0}^{\al_{z, +}(m)}e^{-s^{2}/2m}   \sup_{x\in I_{s}}  \Pbm\left[\sqrt{2T_{y,l}^{m,x^{2}/2}} \le  \al_{z, +}(   l),\,l=m+1,\ldots,k-1;\, \right.
\nn\\
&&\left.\hspace{2.8 in} \sqrt{2T_{y,k}^{m,x^{2}/2}}\in I_{\al_{z, +}(k)-j}\right].\label{18.20g}
\end{eqnarray}

Recall that $\al_{z, +}(   l)= \rho_{L} l+z + l^{1/4}_{L}$. 
 Using this we can write the last  probability as 
\bea
&&
K_{1,s}:=\Pbm \left[ \sqrt{2T_{y,l}^{m,x^{2}/2}} \le  \rho_{L} l+z + l^{1/4}_{L} \mbox{\,\ for }l=m+1,\ldots, k-1;\right.\nn\\ 
&&\left.\hspace{2.5in}\sqrt{2T_{y,k}^{m,x^{2}/2}}\in 
I_{  \rho_{L} k+z + k^{1/4}_{L}- j } \right].\label{18.21}
\eea

Using the fact that for all $1\leq l\leq k-1$
\begin{equation}
 l_{L}^{1/4}\leq l_{k}^{1/4}+ k_{L}^{1/4},\label{bar.6}
\end{equation}
see \cite{BRZ}, 
 it follows that   
 \bea
&&
K_{1,s}\leq \Pgw_{x^{2}/2} \left[ \sqrt{2T_{l} }  \le  \rho_{L}\(m+l\)+ z+ \(m+l\)^{1/4}_{k}+  k_{L}^{1/4},\right.\nn\\ 
  &&\left.\hspace{.2in}\mbox{\,\ for }l=1,\ldots, k-m-1;\sqrt{2T_{k-m} }\in  I_{\rho_{L}k+ z+ k^{ 1/4}_{L}- j } \right].\label{18.22}
\eea

Thus using (\ref{eq: GW curved barrier upper bound}),   with $a=\rho_{L}m+ z+ (m^{1/4}_{k}+k_{L}^{1/4}) $ and $b =\rho_{L}k+ z+ k^{ 1/4}_{L}$, $y=\rho_{L}k+ z+ k^{ 1/4}_{L}- j $,
\be
K_{1,s}\leq c\frac{\left(1 + a -x\right)\(   1+ j\) }{k-m}\sqrt{\frac{x }{y\(k-m\)}}e^{-\frac{\left(\rho_{L}k+ z+ k^{ 1/4}_{L}- j -x\right)^{2}}{2\(k-m\)}}.\nn
\ee
We have
\bea
 e^{-\frac{\left(\rho_{L}k+ z+ k^{ 1/4}_{L}- j  -x\right)^{2}}{2\(k-m\)}}&\leq  &ce^{-\frac{\left(\rho_{L}k\right)^{2}}{2k(1-m/k)}}\,e^{ -2\(   z +  k_{L}^{1/4} - j-x \)}\nn\\
&\le &c e^{2k\frac{\log L}{L}}e^{-2(k+m)-2z-2 k_{L}^{1/4}+2j+2x},
\nn\eea
  $\frac{ 1 }{k-m} \sqrt{\frac{x }{y\(k-m\)}}\,\, e^{2k\frac{\log L}{L}}\asymp \sqrt{x}$, and by assumption, \[a-x\leq c\left(k_{L}^{1/4}+m+z\right).\] Hence we can bound (\ref{18.20g}) by 
 \be
c\left(1 +  k_{L}^{1/4}+m+z\right)\(   1+ j\) e^{-2k-2z-2  k_{L}^{1/4}+2j }\sum_{s=0}^{\al_{z, +}(m)}\sqrt{s}e^{-(s-2m)^{2}/2m}.
\nn
\ee
 
 Our Lemma follows.
\qed

 \bl
\label{prop:BarrierSecGWProp} For all $L $ sufficiently  large,   and all $0\leq   z\leq \log L$,  
\bea 
&&\Pbm\left[   \sqrt{2T_{y, l}^{1}}\le \al_{z,-}\left(l\right)\mbox{\,\ for }l=1,\ldots,L-1;\,\sqrt{2T_{y,L}^{1}}\geq \rho_{L} +z \right]\nn\\
&&\leq \Pbm \left[   \sqrt{2T_{y, l}^{1}}\le  \rho_{L}l +z\mbox{\,\ for }l=1,\ldots,L-1;\,\sqrt{2T_{y,L}^{1}}\geq \rho_{L}L +z \right]\nn\\
&&\hspace{2.5 in}  \leq c (1+z) e^{ -2L-2z-z^{2}/4L}.\label{14.6}
\eea
and
\bea 
&&\Pbm \left[   \sqrt{2T_{y, l}^{1}}\le \al_{z,-}\left(l\right)\mbox{\,\ for }l=1,\ldots,L-1;\,\sqrt{2T_{y, L}^{1}}\in I_{\rho_{L}L +z }\right]\nn\\
&&\hspace{2 in}  \geq c (1+z)   e^{ -2L-2z-z^{2}/4L}.\label{14.6a}
\eea

Similar results hold if we delete the barrier condition on some fixed finite interval.
\el

{\bf   Proof of Lemma \ref{prop:BarrierSecGWProp}: }The first inequality in (\ref{14.6}) is obvious. Theorem \ref{thm: GW Barrier} requires that $ y\leq b$ which we will not have if we go all the way to $L$. Instead, using the Markov property at $l=L-1$ and (\ref{eq: tail bound}) we bound 
\bea 
&& \Pbm \left[ \sqrt{2T_{y, l}^{1}}\le \rho_{L}l +z\mbox{\,\ for }l=1,\ldots,L-1;\,\sqrt{2T_{y, L}^{1}}\geq \rho_{L} +z \right]\nn \\
&&\leq c\sum_{j=1}^{\rho_{L}(L-1)+z}\Pbm \left[   \sqrt{2T_{y, l}^{1}}\le  \rho_{L}l +z\mbox{\,\ for }l=1,\ldots,L-2;\right.\nn\\
&&\left.\hspace{1 in}\sqrt{2T_{y, L-1}^{1}}\in I_{ \rho_{L}(L-1)+z-j} \right]e^{-j^{2}/2}.\label{14.6ub1}
\eea

If $j\geq L/2$, then $e^{-j^{2}/2}\leq e^{-L^{2}/8}$ so we get a bound much smaller than (\ref{14.6}). Thus we need only bound the sum over $1\leq j\leq L/2$.
  
It follows  from (\ref{eq: GW curved barrier upper bound}), with $a=  z$, $ b=\rho_{L}(L-1)+z,\,y=\rho_{L}(L-1)+z-j$ that the last probability is bounded by
\bea
&&
c\frac{(1+z) \(   1+j\) }{L}\sqrt{\frac{1}{L^{ 2}}}e^{ -\(\rho_{L}(L-1)+z-j\)^{ 2}/2(L-1)}\nn\\
&&\hspace{1 in}\leq c(1+z) \(   1+j\) e^{ -2L-2(z-j)-(z-j)^{2}/4L},\nonumber 
\eea
and our upper bound follows after summing over $j$. 

The lower bound follows similarly  using (\ref{eq: gw lower bound}). The last statement in our Lemma  comes from the last statement in Theorem  \ref{thm: GW Barrier}.
\qed


\bl\label{prop:BarrierSecGWgenlow}
 If   $k\geq L/2, \,\,0\leq   z\leq \log L$ and $L $ is sufficiently large,  then uniformly in  $ 0\leq  p\leq k $,
\begin{eqnarray}
&&\Pbm \left[\sqrt{2\mathcal{T}_{y,l}^{1}}\le  \rho_{L}l +z \mbox{\,\ for }l=1,\ldots,k-1;\,\, \sqrt{2\mathcal{T}_{y,k}^{ 1}}\in I_{\rho_{L}k+z-p}\right] 
\label{genlow.1}\\
&& \leq C(1+z)\(   1+p\)  e^{ -2k-2(z -p)-(z -p)^{2}/4k  }. \nonumber
\end{eqnarray}

\el

{\bf   Proof of Lemma \ref{prop:BarrierSecGWgenlow}: }  
Using Theorem \ref{thm: GW Barrier} with $a=  z$, $y=\rho_{L}k+z-p, b=\rho_{L}k+z$ this is bounded by 
\bea
&&
 c\frac{(1+z)\(   1+p\) }{k}\sqrt{\frac{1}{k^{2}}}e^{ -\( \rho_{L}k+z-p\)^{ 2}/2k}\label{genlow.4}\\
 &&\leq C\frac{(1+z)\(   1+p\) }{k^{2}} e^{2\log (L)k/L}e^{ -2k-2(z -p)-(z -p)^{2}/4k }.\nn
\eea
   (\ref{genlow.1})  follows  since  by the convexity of $\log$ we have  $  e^{2\log (L)k/L} \leq e^{2\log (k)}=k^{2}$.
\qed

\bl\label{prop:BarrierSecGWgenhigh}
 If   $k\leq \log^{5}L,\,\,0\leq   z\leq \log L$ and $L $ is sufficiently large,  then uniformly in $  m=  \rho_{L}k+z-t $, 
  \begin{eqnarray}
 && \Pbm \left[ \sqrt{2\mathcal{T}_{y,l}^{k,m^{2}/2}}\le \rho_{L}l +z\mbox{\,\ for }l=k+1,\ldots,L-1;\,\, \sqrt{2\mathcal{T}_{y,L}^{ k,m^{2}}}\geq\rho_{L}L+z \right]
\nn\\
 &&\hspace{1 in}\leq c(1+ t)m^{1/2} e^{- 2(L-k)-2t-t^{2}/4(L-k)}.  \label{14.8}
 \end{eqnarray} 

 
\el

{\bf   Proof of Lemma \ref{prop:BarrierSecGWgenhigh}: } As before, we can bound  the probability by
\bea
&&\sum_{j=0}^{ \rho_{L}(L-1)+z}\Pbm \left[ \sqrt{2\mathcal{T}_{y,l}^{k,m^{2}/2}}\le  \rho_{L}l +z\mbox{\,\ for }l=k+1,\ldots,L-2;\right.\\
&&\hspace{1.5 in}\left.\, \,\sqrt{2T_{y, L-1}^{k,m^{2}/2}}  \in I_{ \rho_{L}(L-1)+z-j}\right]e^{-j^{2}/2}.\nn\label{genhigh.2}
\eea 
Also, as before,    we need only bound the sum over $1\leq j\leq L/2$.

It follows  from (\ref{eq: GW curved barrier upper bound}),   with the $x$ of that estimate given by  $  m= a-t$  and $a= \rho_{L}k+z,      b=\rho_{L}(L-1)+z, y= \rho_{L}(L-1)+z-j$  that the  last probability is bounded by
\be 
c\frac{\left(1+ t\right) \left(1+ j\right)   }{L-k}\sqrt{\frac{m}{L(L-k)}}e^{-\frac{\left(\rho_{L}\left(L-k-1\right)+ t- j \right)^{2}}{2\left(L-k-1\right)}},\label{genhigh.4}
\ee
and
\bea
&&
e^{-\frac{\left(\rho_{L}\left(L-k-1\right)+ t- j \right)^{2}}{2\left(L-k-1\right)}}\nn\\
&&\leq c e^{-\frac{\left(\rho_{L}\left(L-k-1\right)\right)^{2}}{2(L-k-1)}-2(t-j)-(t-j)^{2}/2(L-k)}\nn\\
&&\leq Ce^{2\log (L)(L-k)/L}e^{-2\left(L-k\right)-2(t-j)-(t-j)^{2}/2(L-k)}.\label{genhigh.5}
\eea
This gives  
  \begin{eqnarray}
 &&\hspace{-.2 in} \Pbm \left[ \sqrt{2\mathcal{T}_{y,l}^{k,m^{2}/2}}\le \rho_{L}l +z\mbox{\,\ for }l=k+1,\ldots,L-1;\,\, \sqrt{2\mathcal{T}_{y,L}^{ k,m^{2}/2}}\geq\rho_{L}L+z \right]
 \label{genhigh.1}\nn\\
 &&\leq C\sum_{j=0}^{ \rho_{L}(L-1)+z} \frac{\left(1+ t\right) \left(1+ j\right)    }{L-k}\sqrt{\frac{m}{L(L-k)}}L^{\frac{2(L-k)}{L}} \nn\\
 &&\hspace{2 in}e^{-2\left(L-k\right)-2(t-j)-(t-j)^{2}/2(L-k)}e^{-j^{2}/2}.  \nonumber
 \end{eqnarray}



(\ref{14.8})   follows  since  by the convexity of $\log$,  $  e^{2\log (L)(L-k)/L} \leq e^{2\log ((L-k))}=(L-k)^{2}$.        
\qed

The following Lemma states the barrier estimates needed for the proof of the left tail estimates in Theorem \ref{theo-tight0}. For notation see Section \ref{sec-lowertail}. The proof of this Lemma is similar to the proofs of Lemmas \ref{prop:BarrierSecGWProp} -\ref{prop:BarrierSecGWgenhigh}.
 
 \bl\label{lem-NegBarrier}  For all $L $ sufficiently  large,   and all $|z|\leq \log L$,  
\bea 
&&\Pbm\left[  \sqrt{2T_{y, l}^{1}}\le \wh\gamma_{z,-}\left(l\right)\mbox{\,\ for }l=1,\ldots,L-1;\,\sqrt{2T_{y,L}^{1}}\geq \rho_{L}L +z \right]\nn\\
&&\leq \Pbm \left[   \sqrt{2T_{y, l}^{1}}\le  \wh\bb_{ z}(l)\mbox{\,\ for }l=1,\ldots,L-1;\,\sqrt{2T_{y,L}^{1}}\geq \rho_{L}L +z \right]\nn\\
&&\hspace{3 in}  \leq c  e^{ -2L-2z-z^{2}/4L}.\label{14.6neg}
\eea
and
\bea 
&&\Pbm \left[   \sqrt{2T_{y, l}^{1}}\le \wh\gamma_{z,-}\left(l\right)\mbox{\,\ for }l=1,\ldots,L-1;\,\sqrt{2T_{y, L}^{1}}\in I_{\rho_{L}L +z }\right]\nn\\
&&\hspace{3 in}  \geq c    e^{ -2L-2z-z^{2}/4L}.\label{14.6aneg}
\eea

Similar results hold if we delete the barrier condition on some fixed finite interval.

If $k\geq L/2,\,\, |z|\leq \log L$, and $L$ is sufficiently large, then uniformly in  $  p\leq  k$,
\begin{eqnarray}
&&\Pbm \left[\sqrt{2\mathcal{T}_{y,l}^{1}}\le \wh \bb_{ z}(l) \mbox{\,\ for }l=1,\ldots,k-1;\,\, \sqrt{2\mathcal{T}_{y,k}^{ 1}}\in I_{\wh  \bb_{z}(k)-p}\right] 
\label{genlow.1neg}\\
&& \leq C \(   1+p\)  e^{ -2k-2(z -p)- (z -p)^{2}/4k}. \nonumber
\end{eqnarray}

If $k\leq \log^{5} L ,\,\, |z|\leq \log L$, and $L$ is sufficiently large, then uniformly in $  m=  \wh\bb_{ z}(k)-t $,
  \begin{eqnarray}
 && \Pbm \left[ \sqrt{2\mathcal{T}_{y,l}^{k,m^{2}}}\le  \wh\bb_{ z}(l)\mbox{\,\ for }l=k+1,\ldots,L-1;\,\, \sqrt{2\mathcal{T}_{y,L}^{ k,m^{2}}}\geq\rho_{L}+z \right]
 \label{14.8neg}\\
 &&\leq c(1+ t)m^{1/2} e^{- 2(L-k)- 2t-t^{2}/4(L-k) }.  \nonumber
 \end{eqnarray} 
 \el

\bigskip
\noindent
\begin{tabular}{lll} 
& Jay Rosen\\
& Department of Mathematics\\
&  College of Staten Island, CUNY\\
& Staten Island, NY 10314 \\
& jrosen30@optimum.net\\
\end{tabular}

\end{document}